\documentclass[twocolumn,final]{svjour3} 
\usepackage{amsmath}
\usepackage{amssymb}
\usepackage{graphicx}
\usepackage{color}
\usepackage{subfigure}
\usepackage{ifthen}  % Conditional 
\usepackage[utf8]{inputenc} %unicode support
\usepackage{bm}
\usepackage{float}
\smartqed

\begin{document}

\title{Generalised approach to modelling a three-tiered microbial food-web}
\author{Tewfik Sari \and Matthew J. Wade}

\institute{T. Sari \at Irstea, UMR Itap, Montpellier, France \& Universit\'e de Haute Alsace, Laboratoire de Math\'ematiques, Mulhouse, France. \\\\email{: tewfik.sari@irstea.fr}
\and M. J. Wade \at School of Civil Engineering of Geosciences, Newcastle University, Newcastle-upon-Tyne NE1 7RU, United Kingdom.\\\\email{: matthew.wade@newcastle.ac.uk}}
\maketitle

\begin{abstract}
Ecological modelling of increasingly more complex microbial populations is necessary to reflect the highly functional and diverse behaviour inherent to many systems found in reality. Anaerobic digestion is one such process that has benefitted from the application of mathematical analysis not only for characterising the biological dynamics, but also to investigate emergent behaviour not apparent by simulation alone. Nevertheless, the standard modelling approach has been to describe biological systems using sets of differential equations whose kinetics are generally described by some empirically derived function of growth. The drawbacks of this are two-fold; the growth functions are derived from empirical studies that may not be representative of the system to be modelled and whose parameters may not have a mechanistic meaning, and mathematical analysis is restricted by a conformity to an assumption of the dynamics. Here, we attempt to address these challenges by investigating a generalised form of a three-tier chlorophenol mineralising food-web previously only analysed numerically. We examine the existence and stability of the identified steady-states and find that, without a decay term, the system may be characterised analytically. However, it is necessary to perform numerical analysis for the case when maintenance is included, but in both cases we verify the discovery of two important phenomena; i) the washout steady-state is always stable, and ii) the two other steady-states can be unstable according to the initial conditions and operating parameters.
\end{abstract}


\section{Introduction}\label{Intro}

The mathematical modelling of engineered biological systems has entered a new era in recent years with the expansion and standardisation of existing models aimed at collating disparate components of these processes and provide scientists, engineers and practitioners with the tools to better predict, control and optimise them. These forms of mechanistic models emerged initially with the Activated Sludge Models~\cite{henze87,henze99} for wastewater treatment processes, followed by the Anaerobic Digestion Model No. 1 (ADM1)~\cite{batstone02} a few years later. The development of ADM1 was enabled largely due to the possibilities for better identification and characterisation of functional groups responsible for the discrete degradation steps operating in series within anaerobic digesters. It describes a set of fairly complex stoichiometric and kinetic functions representing the standard anaerobic process, remaining the scientific benchmark to the present day, despite an general understanding of its limitations in describing all necessary biochemical transformations. Indeed, there has been a growing argument that the model should take advantage of improved empirical understanding and extension of biochemical processes included in its structure, to acquire a better trade-off between model realism and complexity~\cite{jeppsson13}.

It has previously been shown that the simplification or reduction in model complexity can preserve biological meaning whilst reducing the computational effort required to find mathematical solutions of the model equations~\cite{garcia13,bornhoft13,weedermann15}. Whilst simpler models are approximations of real systems, it can be beneficial to consider a reduced model to better understand biological phenomena of sub-processes without the need to consider extraneous system parameters and variables, which tend to make mathematical analysis intractable and cumbersome. Nevertheless, even with gross simplification of a biological system based on a set of ordinary differential equations (ODEs) of relatively low dimensionality, analytical techniques are unable to provide general solutions for the system and numerical methods must suffice. 

As an example of this, for anaerobic digestion, a previous study investigated the effect of maintenance on the stability of a two-tiered `food-chain' comprising two species and two substrates~\cite{xu11}. Although the authors were not able to determine the general conditions under which this four dimensional syntrophic consortium was stable, further work has shown that a model with generality can be used to answer the question posed, determining that the two-tiered food-chain is always stable when maintenance is included~\cite{sari14b}.

In a more recent example, the model described by~\cite{xu11} was extended by the addition of a third organism and substrate to create a three-tiered `food-web'~\cite{wade15}. In this model, the stability of some steady-states could be determined analytically, but due to the complexity of the Jacobian matrix for certain steady-states, local solutions were necessary using numerical analysis, when considering the full system behaviour. Although the results were important in revealing emergent properties of this extended model, the motivation of this work is to determine whether the approach carried out in~\cite{sari14b}, can be applied to the three-tiered model from~\cite{wade15}, to provide some general properties of that system.

The paper is organised as follows. In Section~\ref{sec:model}, we present a description of the model to be investigated, before providing an alternative reduction of its structure than that given by~\cite{wade15}, in Section~\ref{mui}. With Section~\ref{sec:ES} we demonstrate the existence of the three steady-states and define four interesting cases for specific parameter values that are investigated using the analytical solutions, whilst also indicating the regions of existence of the steady-states for the operating parameter values (dilution rate and substrate input concentration). In Section~\ref{sec:SS} we perform local stability analysis of the steady-states without maintenance, and in Section~\ref{numanaly}, perform a comprehensive numerical stability analysis of the four cases for both the model with and without a decay constant. We show that our approach leads to the discovery of five operating regions, in which one leads to the possibility of instability of the positive steady state, where all three organisms exist, a fact that has not be reported by~\cite{wade15}. Indeed, we prove that a stable limit-cycle can occur in this region. Finally, in Section~\ref{kinetic}, we make comment on the role of the kinetic parameters used in the four example cases, in maintaining stability, which points to the importance of the relative aptitude of the two hydrogen consumers in sustaining a viable chlorophenol mineralising community. In the Appendix we describe the numerical method used in Section~\ref{numanaly} and we give the proofs of the results.

%============================================
\section{The model}\label{sec:model}
%============================================
The model developed in~\cite{wade15} has six components, three substrate (chlorophenol, phenol and hydrogen) and three biomass (chlorophenol, phenol and hydrogen degraders) variables. The substrate and biomass concentrations evolve according to the six-dimensional dynamical of ODEs

{\small
\begin{align} 
\frac{{\rm d}X_{\rm ch}}{{\rm d}t} &=  -D X_{\rm ch} + Y_{\rm ch}f_0\left(S_{\rm ch},S_{\rm H_2}\right) X_{\rm ch}- 
k_{\rm dec, ch}X_{\rm ch}\label{model0a}\\
%[2mm]
\frac{{\rm d}X_{\rm ph}}{{\rm d}t} &=  -D X_{\rm ph} + Y_{\rm ph}f_1\left(S_{\rm ph},S_{\rm H_2}\right) X_{\rm ph}- 
k_{\rm dec, ph}X_{\rm ph}\label{model0b}\\
%[2mm]
\frac{{\rm d}X_{\rm H_2}}{{\rm d}t} &=  -D X_{\rm H_2} + Y_{\rm H_2}f_2\left(S_{\rm H_2}\right) X_{\rm H_2}- 
k_{\rm dec, H_2}X_{\rm H_2}\label{model0c}\\
%[2mm]
\frac{{\rm d}S_{\rm ch}}{{\rm d}t} &=  D \left(S_{\rm ch,in} - S_{\rm ch}\right) - 
f_0\left(S_{\rm ch},S_{\rm H_2}\right) X_{\rm ch}\label{model0d}\\
%[2mm]
%\end{align}
%\begin{align}
\frac{{\rm d}S_{\rm ph}}{{\rm d}t} &=  D \left(S_{\rm ph,in} - S_{\rm ph}\right)
+\frac{224}{208}\left(1-Y_{\rm ch}\right)f_0\left(S_{\rm ch},S_{\rm H_2}\right) X_{\rm ch} \nonumber\\
%[2mm]
& - f_1\left(S_{\rm ph},S_{\rm H_2}\right) X_{\rm ph}\label{model0e}\\
%[2mm]
\frac{{\rm d}S_{\rm H_2}}{{\rm d}t} &=  \left(S_{\rm H_2,in}- S_{\rm H_2}\right) 
+ \frac{32}{224}\left(1-Y_{\rm ph}\right) f_1\left(S_{\rm ph},S_{\rm H_2}\right) X_{\rm ph}\nonumber\\
%[2mm]
&-\frac{16}{208}f_0\left(S_{\rm ch},S_{\rm H_2}\right) X_{\rm ch}- f_2\left(S_{\rm H_2}\right) X_{\rm H_2} \label{model0f}
\end{align}}%

\noindent where $S_{\rm ch}$ and $X_{\rm ch}$ are the chlorophenol substrate and biomass concentrations, $S_{\rm ph}$ and $X_{\rm ph}$ those for phenol and $S_{\rm H_2}$ and $X_{\rm H_2}$ those for hydrogen; $Y_{\rm ch}$, $Y_{\rm ph}$ and $Y_{\rm H2}$ are the yield coefficients, $224/208\left(1-Y_{\rm ch}\right)$ represents the part of chlorophenol degraded to phenol, and $32/224\left(1-Y_{\rm ph}\right)$ represents the part of phenol that is transformed to hydrogen. Growth functions take Monod form with hydrogen inhibition acting on the phenol degrader and represented in $f_{1}$ (see Eq.~\ref{Monod}) as a product inhibition term.

\begin{equation}
\begin{array}{l}
f_0\left(S_{\rm ch},S_{\rm H_2}\right)=\frac{k_{m,\rm ch}S_{\rm ch}}{K_{S,ch}+S_{\rm ch}}\frac{S_{H_2}}{K_{S,\rm H_2,c}+S_{\rm H_2}}
\\[2mm]
f_1\left(S_{\rm ph},S_{\rm H_2}\right)=\frac{k_{m,\rm ph}S_{\rm ph}}{K_{S,\rm ph}+S_{\rm ph}}\frac{1}{1+
\frac{S_{\rm H_2}}{K_{i,\rm H_2}}}
\\
f_2\left(S_{\rm H_2}\right)=\frac{k_{m,\rm H_2}S_{\rm H_2}}{K_{S,\rm H_2}+S_{\rm H_2}}
\end{array}
\label{Monod}
\end{equation}

Here, apart from the four operating (or control) parameters, which are the inflowing concentrations $S_{\rm ch,in}$, $S_{\rm ph,in}$, $S_{\rm H_2,in}$ and the dilution rate $D$, that can vary, all others have biological meaning and are fixed depending on the organisms and substrate considered.
We use the following simplified notations in (Eqs.~\ref{model0a}-\ref{model0f})

{\begin{align}
&
X_0=X_{\rm ch},
\quad
X_1=X_{\rm ph},
\quad
X_2=X_{\rm H_2}
\nonumber\\
&
S_0=S_{\rm ch},
\quad
S_{1}=S_{\rm ph},
\quad
S_2=S_{\rm H_2}
\nonumber\\
&
S_0^{\rm in}=S_{\rm ch,in},
\quad
S_{1}^{\rm in}=S_{\rm ph,in},
\quad
S_2^{\rm in}=S_{\rm H_2,in}
\nonumber\\
&
Y_0=Y_{\rm ch},
\quad
Y_1=Y_{\rm ph},
\quad
Y_2=Y_{\rm H_2}
\nonumber\\
&
Y_3=\frac{224}{208}\left(1-Y_{\rm ch}\right),
\quad
Y_4=\frac{32}{224}\left(1-Y_{\rm ph}\right),
\quad
Y_5=\frac{16}{208}
\nonumber\\
&
a_0=k_{\rm dec,\rm ch},
\quad
a_1=k_{\rm dec,\rm ph},
\quad
a_2=k_{\rm dec, H_2}
\nonumber
\end{align}}
With these notations Eqs.~\ref{model0a}-\ref{model0f} can be written as follows

{\small
\begin{align} 
\frac{{\rm d}X_{0}}{{\rm d}t}& =  -D X_{0} + Y_{0}f_{0}\left(S_{0},S_{2}\right) X_{0}- a_0X_{0}\label{model1a}\\%[2mm]
\frac{{\rm d}X_{1}}{{\rm d}t}& =  -D X_{1} + Y_{1}f_{1}\left(S_{1},S_{2}\right) X_{1}- a_1X_{1}\label{model1b}\\%[2mm]
\frac{{\rm d}X_{2}}{{\rm d}t}& =  -D X_{2} + Y_{2}f_{2}\left(S_{2}\right) X_{2}- a_2X_{2}\label{model1c}\\%[2mm]
\frac{{\rm d}S_{0}}{{\rm d}t}& =  D \left(S_{0}^{\rm in} - S_{0}\right) - f_{0}\left(S_{0},S_{2}\right) X_{0}\label{model1d}\\%[2mm]
%\end{align}
%\begin{align}
\frac{{\rm d}S_{1}}{{\rm d}t}& =  D \left(S_{1}^{\rm in} - S_{1}\right)
+Y_3f_{0}\left(S_{0},S_{2}\right) X_{0} - f_{1}\left(S_{1},S_{2}\right) X_{1}\label{model1e}\\%[2mm]
\frac{{\rm d}S_{2}}{{\rm d}t}& =  D\left(S_{2}^{\rm in}- S_{2}\right)+ 
Y_4 f_{1}\left(S_{1},S_{2}\right) X_{1}-Y_5f_{0}\left(S_{0},S_{2}\right) X_{0} \nonumber\\
& - f_{2}\left(S_{2}\right) X_{2} \label{model1f}
\end{align}}%
In~\cite{wade15}, this model is reduced to a dimensionless form that significantly reduces the number of parameters describing the dynamics.  In this paper we do not assume that the growth functions $f_0$, $f_1$ and $f_2$ have the specific analytical expression (Eq.~\ref{Monod}). We will only assume that the growth functions 
satisfy properties that are listed in Appendix~\ref{general}. Therefore, we cannot benefit from the dimensionless rescaling used by~\cite{wade15}, because this rescaling uses some kinetics parameters of the specific growth functions (Eq.~\ref{Monod}), while we work with general {\em unspecified} growth functions. In Section~\ref{mui} we consider another rescaling that does not use the kinetics parameters. Furthermore, we restrict our analysis to the case where we only have one substrate addition to the system, such that: $S_0^{in} > 0$, $S_1^{in} = 0$, and $S_2^{in} = 0$.

%============================================
\section{Model reduction}\label{mui}
%============================================
To ease the mathematical analysis, we can rescale the system (Eqs.~\ref{model1a}-\ref{model1f}) using the following change of variables adapted from~\cite{sari14b}:

{\begin{align}
&
x_0 = \frac{Y_3Y_4}{Y_0} X_0,
\quad 
x_1 = \frac{Y_4}{Y_1} X_1,
\quad 
x_2 = \frac{1}{Y_2} X_1
\nonumber\\
&
s_0 = Y_3Y_4 S_0,
\quad 
s_1 = Y_4 S_1,
\quad 
s_2 = S_2
\nonumber
\end{align}}
We obtain the following system

{\small\begin{align}
\frac{{\rm d}x_{0}}{{\rm d}t} &=  -D x_{0} + \mu_{0}\left(s_{0},s_{2}\right) x_{0}- a_0x_{0}\label{modela}\\%[2mm]
%\end{align}
%\begin{align}
\frac{{\rm d}x_{1}}{{\rm d}t} &=  -D x_{1} + \mu_{1}\left(s_{1},s_{2}\right) x_{1}- a_1x_{1}\label{modelb}\\%[2mm]
\frac{{\rm d}x_{2}}{{\rm d}t} &=  -D x_{2} + \mu_{2}\left(s_{2}\right) x_{2}- a_2x_{2}\label{modelc}\\%[2mm]
\frac{{\rm d}s_{0}}{{\rm d}t} &= D \left(s_{0}^{{\rm in}} - s_{0}\right) - \mu_{0}\left(s_{0},s_{2}\right) x_{0}\label{modeld}\\%[2mm]
\frac{{\rm d}s_{1}}{{\rm d}t} &= -Ds_{1}+\mu_{0}\left(s_{0},s_{2}\right) x_{0}- 
\mu_{1}\left(s_{1},s_{2}\right) x_{1}\label{modele}\\%[2mm]
\frac{{\rm d}s_{2}}{{\rm d}t} &=  -D s_{2} + \mu_{1}\left(s_{1},s_{2}\right) x_{1}-\omega \mu_{0}\left(s_{0},s_{2}\right) x_{0}%\nonumber\\&
- \mu_2\left(s_{2}\right) x_{2} \label{modelf}
\end{align}}
where the inflowing concentration is

{\begin{equation}
s_0^{{\rm in}} = Y_3Y_4 S_0^{{\rm in}},
\label{eqS0in}
\end{equation}}
%=======================
the growth functions are

{\begin{equation}
\begin{array}{l}
\mu_0(s_0,s_2) = {Y_0}f_0\left(\frac{s_0}{Y_3Y_4},s_2\right)\\[2mm]
\mu_1(s_1,s_2) = {Y_1}f_1\left(\frac{s_1}{Y_4},s_2\right)\\[2mm]
\mu_2(s_2) = Y_2f_2(s_2)
\end{array}
\label{rescaling}
\end{equation}}
and
{\begin{equation}
\omega=\frac{Y_5}{Y_3Y_4}=\frac{1}{2(1-Y_0)(1-Y_1)}
\label{eqomega}
\end{equation}}

The benefit of our rescaling is that it permits to fix in Eqs.~\ref{modela}-\ref{modelf} all yield coefficients to one except that denoted by $\omega$ 
{and defined by (Eq.~\ref{eqomega}), and to discuss the existence and stability with respect to this sole parameter.

Using Eq.~\ref{rescaling} and the growth functions (Eq.~\ref{Monod}), we obtain the model (Eqs.~\ref{modela}-\ref{modelf}) with the following Monod-type growth functions

\begin{equation}
\begin{array}{l}
\mu_0\left(s_0,s_2\right)=\frac{m_0s_0}{K_{0}+s_{0}}\frac{s_2}{L_0+s_{2}}\\[2mm]
\mu_1\left(s_1,s_2\right)=\frac{m_1s_1}{K_{1}+s_{1}}\frac{1}{1+s_{2}/K_{i}}\\[2mm]
\mu_2\left(s_{2}\right)=\frac{m_{2}s_{2}}{K_{2}+s_{2}}
\label{monMonod}
\end{array}
\end{equation}
where

{\begin{align}
&
m_0=Y_0k_{m,\rm ch},\quad K_0 = Y_3Y_4K_{s,\rm ch},\quad L_0 = K_{S,\rm H_2,c}
\nonumber\\
&
m_1 = Y_1k_{m,\rm ph},\quad  K_1 = Y_4K_{s,\rm ph},\quad K_i = K_{i,\rm H_2}
\label{cc}\\	
&
m_2 = Y_2k_{m,\rm H_2},\quad K_2 = K_{S,\rm H_2}
\nonumber
\end{align}}
For the numerical simulations we will use the nominal values in Table~\ref{table0} given in~\cite{wade15}. 

\begin{table}[ht]
\centering
\begin{tabular}{lll}
\hline
Parameters & Nominal values     & Units                \\ \hline
$k_{m,\mathrm{ch}}$          & 29                         & $\mathrm{kgCOD_{S}/kgCOD_{X}/d}$                     \\
$K_{S,\mathrm{ch}}$          & 0.053                       & $\mathrm{kgCOD/m^{3}}$ \\
$Y_{\mathrm{ch}}$            & 0.019                       & $\mathrm{kgCOD_X/kgCOD_S}$                       \\
$k_{m,\mathrm{ph}}$          & 26                          & $\mathrm{kgCOD_{S}/kgCOD_{X}/d}$                     \\
$K_{S,\mathrm{ph}}$          & 0.302                        & $\mathrm{kgCOD/m^{3}}$ \\
$Y_{\mathrm{ph}}$            & 0.04                        & $\mathrm{kgCOD_X/kgCOD_S}$                       \\
$k_{m,\mathrm{H_{2}}}$    & 35                          & $\mathrm{kgCOD_{S}/kgCOD_{X}/d}$                     \\
$K_{S,\mathrm{H_{2}}}$       & 2.5$\times10^{-5}$ & $\mathrm{kgCOD/m^{3}}$ \\
$K_{S,\mathrm{H_{2},c}}$  & 1.0$\times10^{-6}$ & $\mathrm{kgCOD/m^{3}}$ \\
$Y_{\mathrm{H_{2}}}$        & 0.06                        & $\mathrm{kgCOD_X/kgCOD_S}$                       \\
$k_{\mathrm{dec,i}}$        & 0.02                        & $\mathrm{d^{-1}}$      \\
$K_{I,\mathrm{H_{2}}}$      & 3.5$\times10^{-6}$ & $\mathrm{kgCOD/m^{3}}$ \\ \hline
\end{tabular}
\caption{Nominal parameter values.}
\label{table0}
\end{table}

%============================================
\section{Existence of steady-states }\label{sec:ES}
%============================================

A steady-state of Eqs.~\ref{modela}-\ref{modelf} is obtained by setting the right-hand sides equal to zero:

{\small\begin{align}
\left[\mu_{0}\left(s_{0},s_{2}\right)-D-a_0\right] x_{0} &=0 \label{m1}\\[2mm]
\left[\mu_{1}\left(s_{1},s_{2}\right)-D-a_1\right] x_{1} &=0 \label{m2}\\[2mm]
\left[\mu_{2}\left(s_{2}\right) -D-a_2\right]x_{2}       &=0 \label{m3}\\[2mm]
D \left(s_{0}^{{\rm in}} - s_{0}\right) - \mu_{0}\left(s_{0},s_{2}\right) x_{0} &=0 \label{m4}\\[2mm]
-D s_{1}+\mu_{0}\left(s_{0},s_{2}\right) x_{0}- \mu_{1}\left(s_{1},s_{2}\right) x_{1}&=0 \label{m5}\\[2mm]
- Ds_{2} + \mu_{1}\left(s_{1},s_{2}\right) x_{1}-\omega \mu_{0}\left(s_{0},s_{2}\right) x_{0}- \mu_2\left(s_{2}\right) x_{2} &=0 \label{m6}
\end{align}}
A steady-state exists (or is said to be `meaningful') if, and only if, all its components are non-negative.

\begin{lemma}\label{lemma1}
The only steady-state of Eqs.~\ref{modela}-\ref{modelf}, for which $x_0=0$ or $x_1=0$, is the steady-state
$${\rm SS1}= (x_0=0,x_1=0,x_2=0,s_0=s_0^{{\rm in}},s_1=0,s_2=0)$$ 
where all species are washed out. This steady-state always exists. It is always stable.
\end{lemma}
%================================================
\begin{figure*}[ht]
\setlength{\unitlength}{1.0cm}
\begin{center}
\begin{picture}(10,4)(0,2.2)
\put(-2.5,-3){\rotatebox{0}{\includegraphics[scale=0.35]{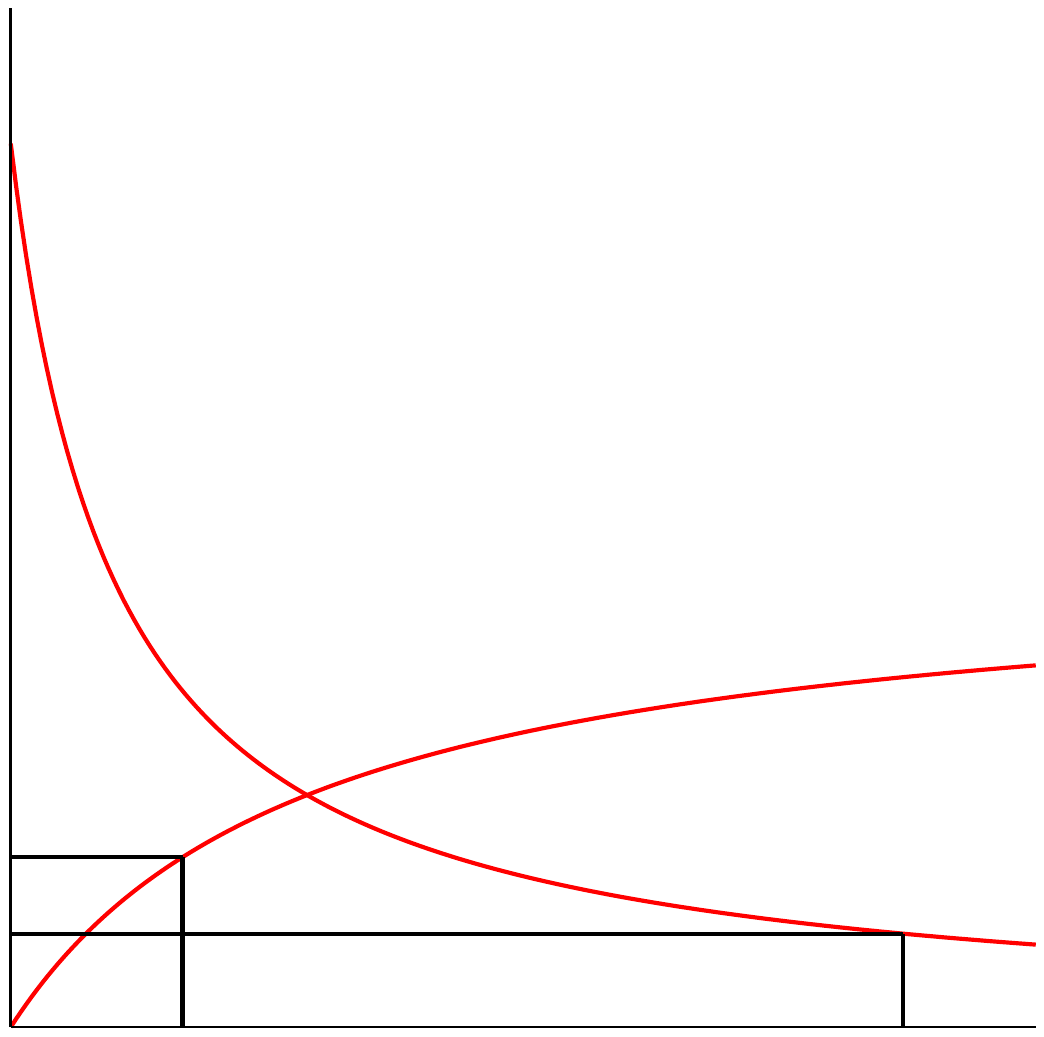}}}
\put(2.5,-3){\rotatebox{0}{\includegraphics[scale=0.35]{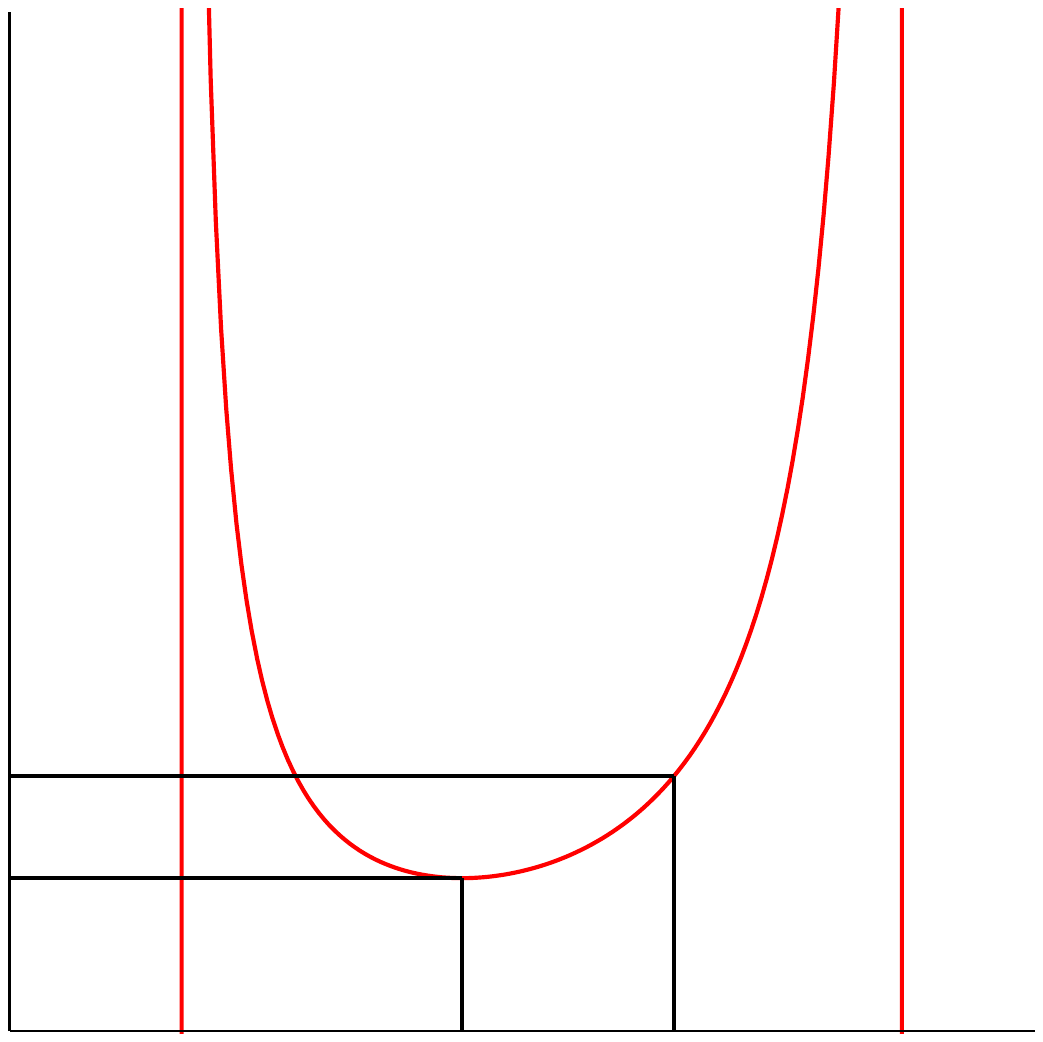}}}
\put(7.5,-3){\rotatebox{0}{\includegraphics[scale=0.35]{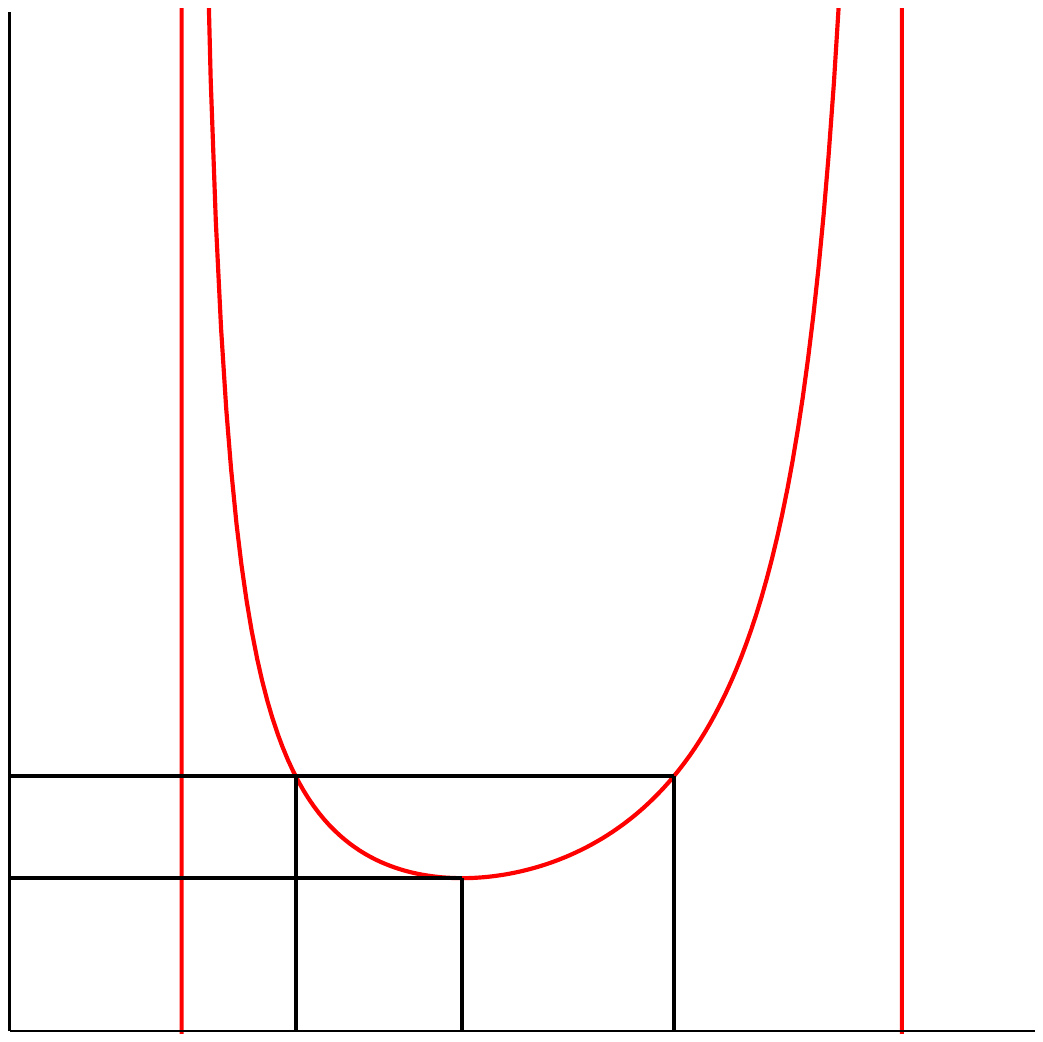}}}
%%%%%%%%%%%%
\put(0,6){(a)}
\put(5,6){(b)}
\put(10,6){(c)}
%===============
\put(2.2,2.3){{$s_2$}}
\put(7.2,2.3){{$s_2$}}
\put(12.2,2.3){{$s_2$}}
%===============
\put(-1.3,5){{$\mu_1(+\infty,s_2)$}}
\put(0.5,3.8){{$\mu_0(+\infty,s_2)$}}
\put(5.5,5){{$\psi(s_2)$}}
\put(10.5,5){{$\psi(s_2)$}}
%================
\put(-2.4,3){{\scriptsize{$D+a_0$}}}
\put(-2.4,2.7){{\scriptsize{$D+a_1$}}}
%==================
\put(2.7,3.3){{\scriptsize{$F_2(D)$}}}
\put(2.7,2.9){{\scriptsize{$F_1(D)$}}}
%===================
\put(8.1,3.3){{\scriptsize{$s_0^{\rm in}$}}}
\put(7.7,2.9){{\scriptsize{$F_1(D)$}}}
%====================
\put(-1,2.1){{\scriptsize{$s_2^0$}}}
\put(1.5,2.1){{\scriptsize{$s_2^1$}}}
\put(4,2.1){{\scriptsize{$s_2^0$}}}
\put(6.5,2.1){{\scriptsize{$s_2^1$}}}
\put(9,2.1){{\scriptsize{$s_2^0$}}}
\put(11.5,2.1){{\scriptsize{$s_2^1$}}}
%====================
\put(5,2.2){{\scriptsize{$\overline{s}_2$}}}
\put(10,2.2){{\scriptsize{$\overline{s}_2$}}}
%=================
\put(5.1,1.8){{\scriptsize{$M_2(D+a_2)$}}}
\put(5.85,2){\vector(0,1){0.4}}
%==================
\put(9.4,2.1){{\scriptsize{$s_2^\flat$}}}
\put(10.8,2.1){{\scriptsize{$s_2^\sharp$}}}
%================
\end{picture}
\end{center}
\caption{Graphical definitions. (a): $s_2^0$ and $s_2^1$. (b) : $\psi(s_2)$, $\overline{s}_2$, $F_1(D)$ and $F_2(D)$. (c): $s_2^\flat$ and $s_2^\sharp$}
\label{fig1}
\end{figure*}

From the previous Lemma we deduce that besides the steady-state SS1, the system can have at most two other steady-states.

SS2: $x_0>0$, $x_1>0$ and $x_2=0$, where species $x_2$ is washed out while species $x_0$ and and $x_1$ exist.

SS3: $x_0>0$, $x_1>0$, and $x_2>0$, where all populations are maintained.

In the following we describe the steady-states SS2 and SS3 of Eqs.~\ref{modela}-\ref{modelf} with the Monod-type growth functions
(Eq.~\ref{monMonod}). The general case with unspecified growth function is provided in Appendix~\ref{general}, with proofs given in Appendix~\ref{proofs}. We use the following notations:

Let $s_2$ be fixed, we define the function $M_0(y,s_2)$ as follows : for all 
$y\in\left[0,\mu_0(+\infty,s_2)=\frac{m_0s_2}{L_0+s_2}\right)$, we let
$$
 M_0(y,s_2)=\frac{K_0y}{\frac{m_0s_2}{L_0+s_2}-y}
$$
Notice that $y\mapsto M_0(y,s_2)$ is the inverse function of the function $s_0\mapsto \mu_0(s_0,s_2)$, that is to say, for all 
$s_0\geq 0$, $s_2\geq 0$ and $y\in[0,\mu_0(+\infty,s_2))$ 
%================
\begin{equation}
s_0=M_0(y,s_2)\Longleftrightarrow y=\mu_0(s_0,s_2)
\label{eqM0explicit}
\end{equation}

Let $s_2$ be fixed, we define the function $M_1(y,s_2)$ as follows : for all 
$y\in\left[0,\mu_1(+\infty,s_2)=\frac{m_1}{1+s_2/K_i}\right)$, we let
$$
M_1(y,s_2)=\frac{K_1y}{\frac{m_1}{1+s_{2}/K_{i}}-y}
$$
Notice that $y\mapsto M_1(y,s_2)$ is the inverse function of the function $s_1\mapsto \mu_1(s_1,s_2)$, that is to say, for all 
$s_1\geq 0$, $s_2\geq 0$ and $y\in[0,\mu_1(+\infty,s_2))$ 
%================
\begin{equation}
s_1=M_1(y,s_2)\Longleftrightarrow y=\mu_1(s_1,s_2)
\label{eqM1explicit}
\end{equation}

We define the function $M_2(s_2)$ as follows : for all 
$y\in\left[0,\mu_2(+\infty)=m_2\right)$, we let
$$
M_2(y)=\frac{K_2y}{m_2-y}
$$
Notice that $y\mapsto M_2(y)$ is the inverse function of the function $s_2\mapsto \mu_2(s_2)$, that is to say, for all 
$s_2\geq 0$ and $y\in[0,\mu_2(+\infty))$ 
%================
\begin{equation}
s_2=M_2(s_2)\Longleftrightarrow y=\mu_2(s_2)
\label{eqM2explicit}
\end{equation}
%================================================
\begin{figure*}[ht]
\setlength{\unitlength}{1.0cm}
\begin{center}
\begin{picture}(10,4)(0,2.8)
\put(-4,-3){\rotatebox{0}{\includegraphics[scale=0.35]{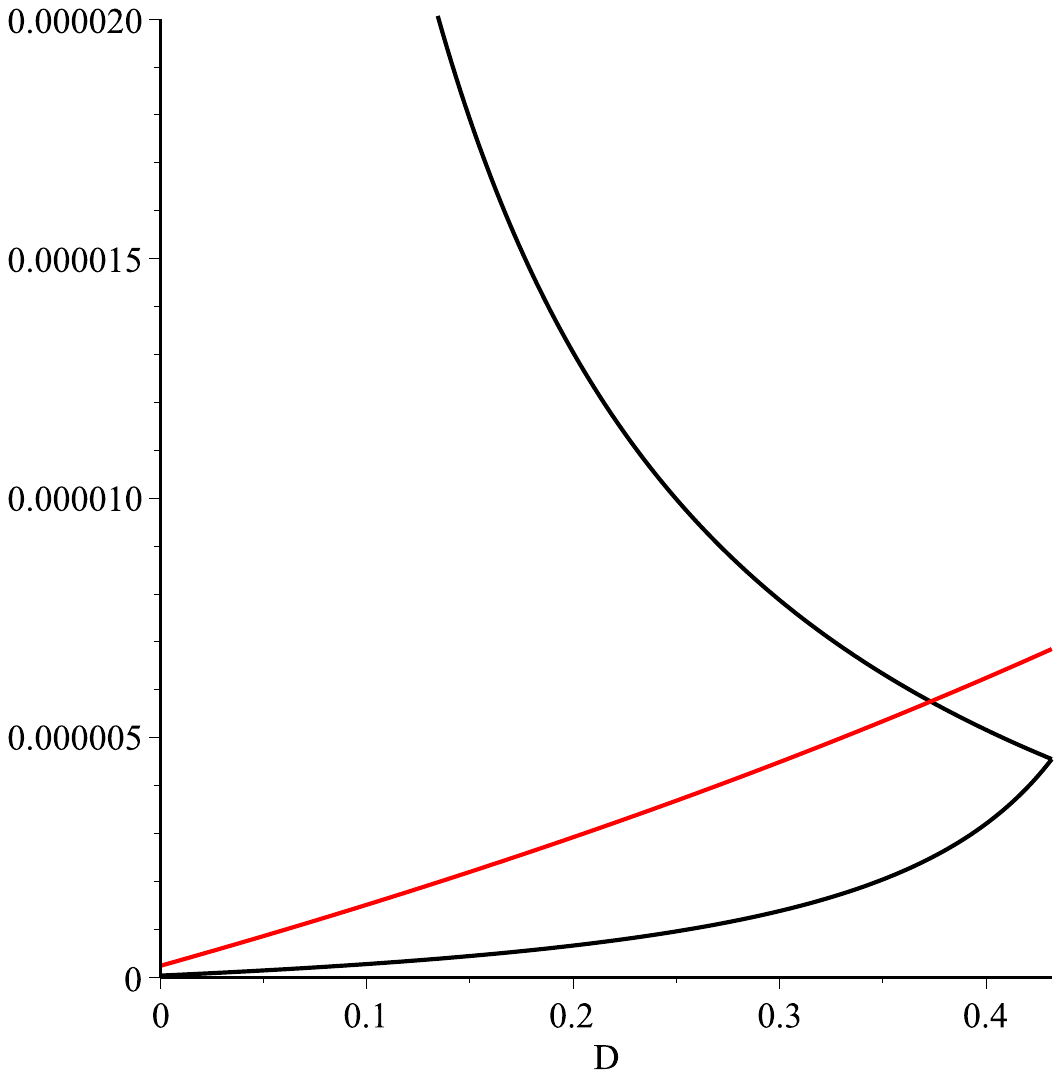}}}
\put(0,-3){\rotatebox{0}{\includegraphics[scale=0.35]{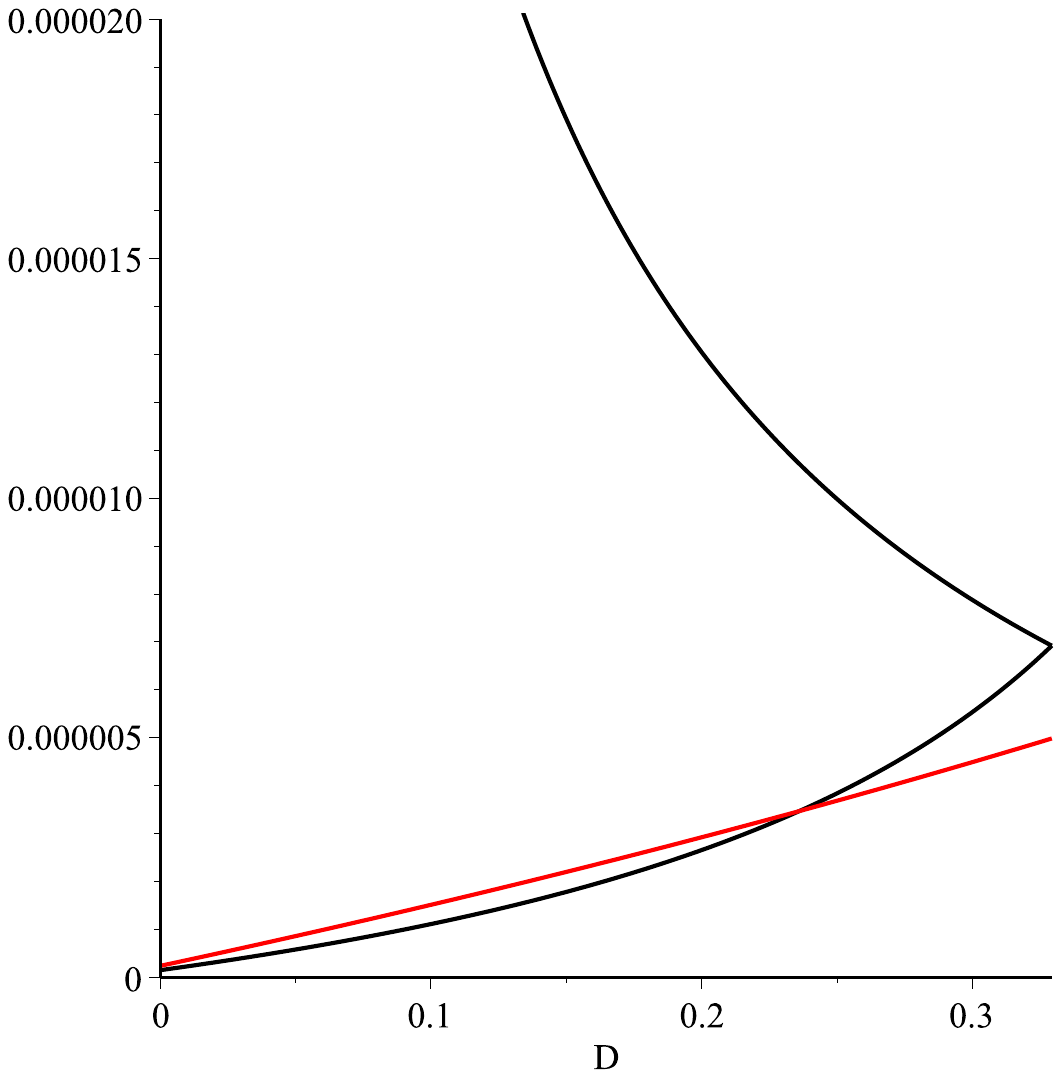}}}
\put(4,-3){\rotatebox{0}{\includegraphics[scale=0.35]{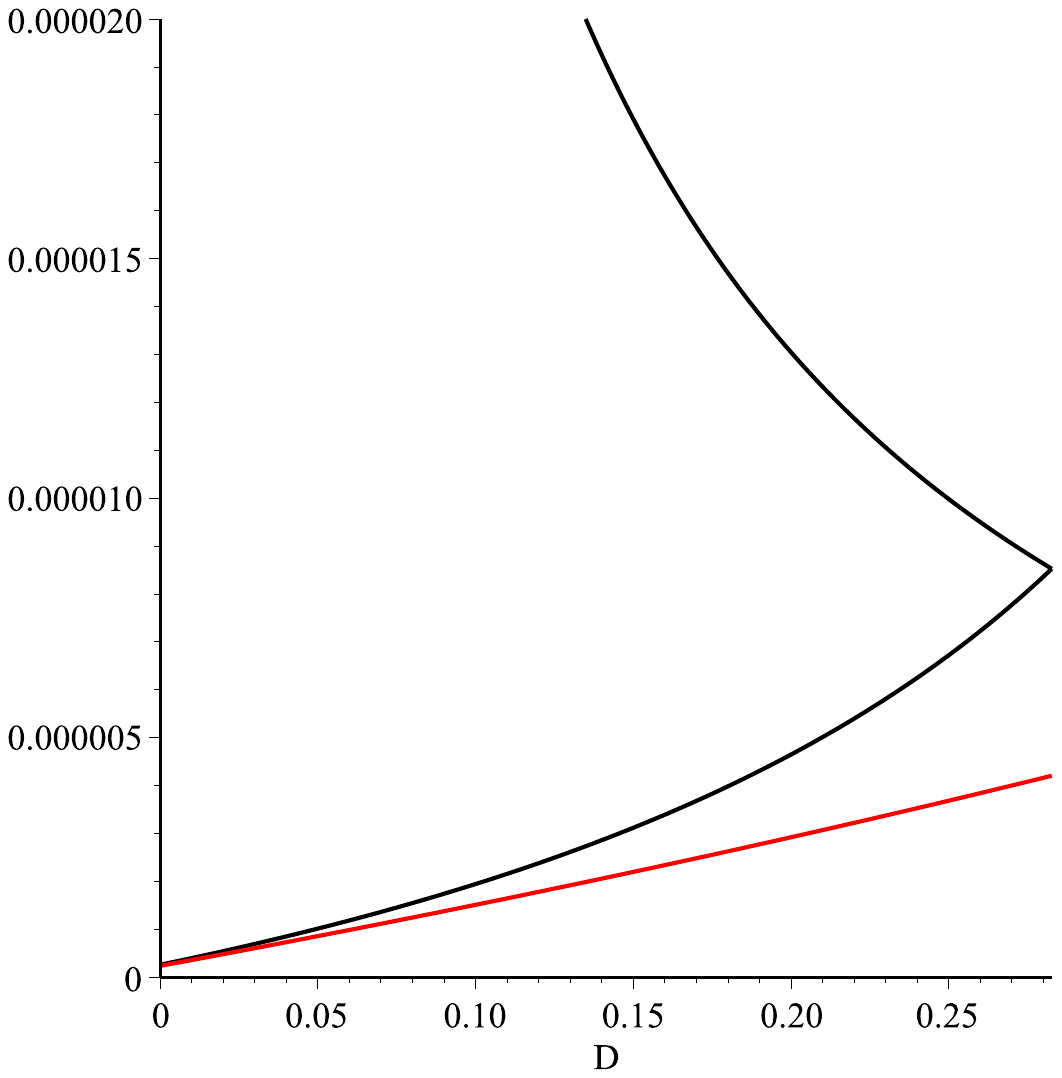}}}
\put(8,-3){\rotatebox{0}{\includegraphics[scale=0.35]{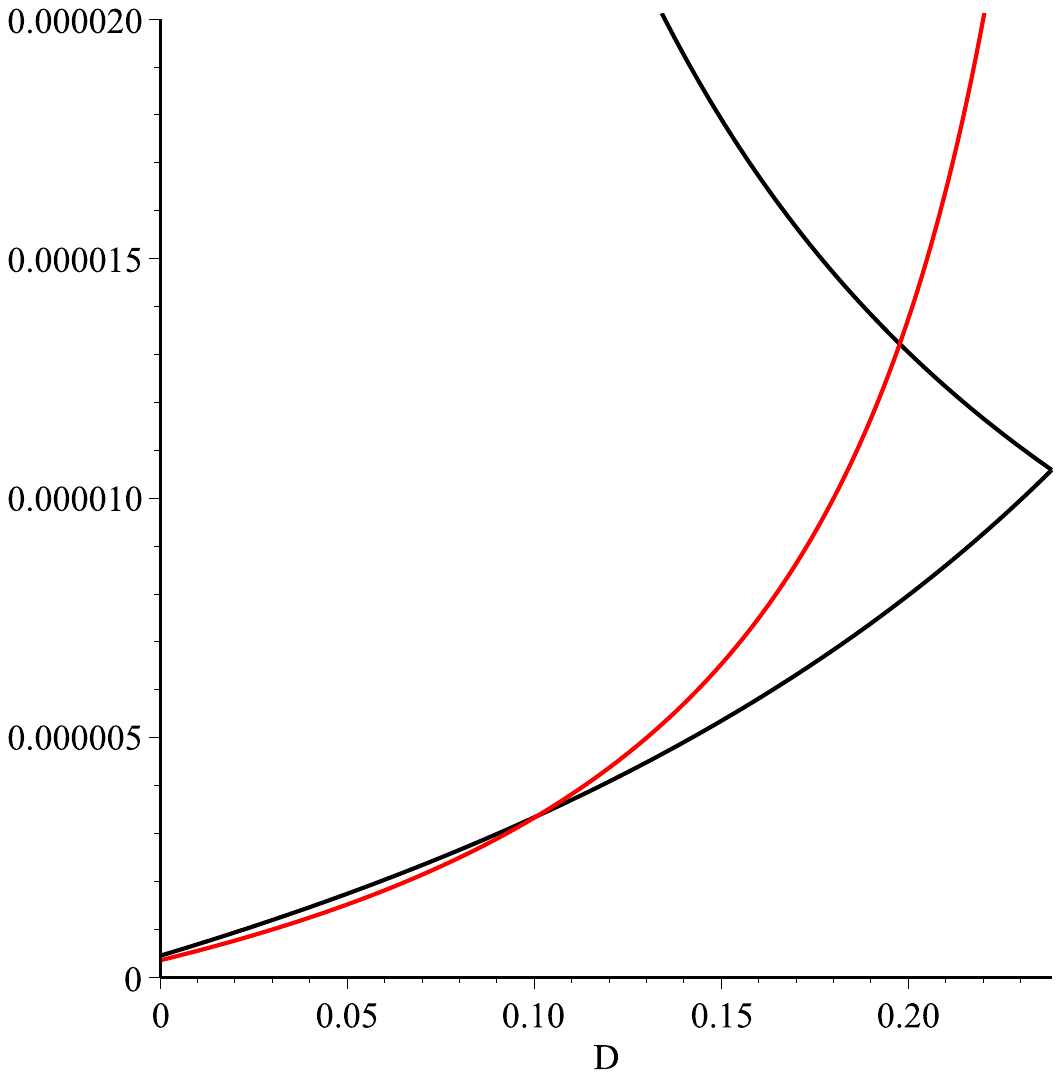}}}
%%%%%%%%%%%%
\put(-1,6){(a)}
\put(3,6){(b)}
\put(7,6){(c)}
\put(11,6){(d)}
%===============
\end{picture}
\end{center}
\caption{ 
Graphs of $s_2^0(D)$ and $s_2^1(D)$ (in black) and $M_2(D)$ (in red) and 
graphical depiction of $I_1=[0,D_1)$, where $D_1$ is the solution of $s_2^0(D)=s_2^1(D)$, and $I_2$.  
(a): $I_2=[0,D_2)$ where $D_2$ is the solution of $M_2(D)=s_2^1(D)$.  
(b) : $I_2=[0,D_2)$ where  $D_2$ is the solution of $M_2(D)=s_2^0(D)$. 
(c): $I_2$ is empty.
(d) : $I_2=(D_{2min},D_{2max})$ where $D_{2min}$ and $D_{2max}$ are the solutions of $M_2(D)=s_2^0(D)$ and $M_2(D)=s_2^1(D)$, respectively.}
\label{figI1I2}
\end{figure*}
%========================
Using the functions $M_0$, $M_1$ and $M_2$ we define the following function:
Let $\omega<1$. Let 
%==================
\begin{equation}
\psi(s_2)=M_0(D+a_0,s_2)+\frac{M_1(D+a_1,s_2)+s_2}{1-\omega}
\label{eqpsiexplicit}
\end{equation}
%=================
Notice that $\psi$ is defined if, and only if,
$$D+a_0<\mu_0(+\infty,s_2)\mbox{ and } D+a_1<\mu_1(+\infty,s_2)$$
which is equivalent to 
$$s_2^0(D)<s_2<s_2^1(D)$$
where 
$$
s_2^0(D) = \frac{L_0(D+a_0)}{m_0-D-a_0},
\quad
s_2^1(D) = \frac{K_i(m_1-D-a_1)}{D+a_1}
$$
are the solutions of equations
\begin{equation}
\mu_0(+\infty,s_2)=D+a_0,\quad \mu_1(+\infty,s_2)=D+a_1
\label{eqs201D}
\end{equation}
respectively, see Fig.~\ref{fig1} (a). 
Straightforward calculations show that
$$
\psi(s_2,D)=
\frac{K_0(D+a_0)}{m_0-D-a_0}\frac{L_0+s_2}{s_2-s_2^0(D)}
+\frac{\frac{K_1(K_i+s_2)}{s_2^1(D)-s_2}+s_2}{1-\omega}
$$
Therefore, $\psi(s_2)>0$ for $s_2^0<s_2<s_2^1$ (see Fig.~\ref{fig1} (b)), 
$$\lim_{s_2\to s_2^0}\psi(s_2)=\lim_{s_2\to s_2^1}\psi(s_2)=+\infty$$
and
{\small $$
\frac{d^2\psi}{ds_2^2}=
\frac{2K_0(D+a_0)}{m_0-D-a_0}\frac{L_0+s_2^0(D)}{\left(s_2-s_2^0(D)\right)^3}
-\frac{2K_1(K_i+s_2^1(D))}{(1-\omega)\left(s_2^1(D)-s_2\right)^3}
$$}
Hence, $\frac{d^2\psi}{ds_2^2}>0$ for all $s_2\in(s_2^0(D),s_2^1(D))$, so that the function $s_2\mapsto\psi(s_2,D)$ is convex and, thus, it has a unique minimum $\overline{s}_2(D)$, see Fig.~\ref{fig1} (b). 

Let $\omega<1$. We define the function
\begin{equation}
F_1(D)=\inf_{s_2\in(s_2^0,s_2^1)}\psi(s_2)=\psi\left(\overline{s}_2\right)
\label{eqF1explicit}
\end{equation}
as shown in Fig.~\ref{fig1} (b). The minimum $\overline{s}_2(D)$ is a solution of an algebraic equation of degree 4 in $s_2$. Although mathematical software, such as {\em Maple}, cannot give its solutions explicitly with respect to the parameters, $\overline{s}_2(D)$ could be obtained analytically since algebraic equations of degree 4 can theoretically be solved by quadratures. We do not try to obtain such an explicit formula. However, if the biological parameters are fixed, the function $\overline{s}_2(D)$ and, hence, 
$F_1(D)=\psi(\overline{s}_2(D),D)$, can be obtained numerically.

The function $F_1(D)$ is defined as long as $s_2^0(D)<s_2^1(D)$. Assuming that $s_2^0(0)<s_2^1(0)$,  $F_1(D)$ is defined 
for $0\leq D<D_1$, where $D_1$ is the positive solution of $s_2^0(D)=s_2^1(D)$ (see Fig.~\ref{figI1I2}). Therefore, $D_1$ is a solution of the 
the second order algebraic equation $\frac{L_0(D+a_0)}{m_0-D-a_0}=\frac{K_i(m_1-D-a_1)}{D+a_1}$.
We denote by 
\begin{equation}
I_1=\{D:s_2^0(D)<s_2^1(D)\}
\label{eqI1}
\end{equation}
the set on which $F_1(D)$ is defined.
%=============== ITEM8

Let $\omega<1$. We define the functions
\begin{align}
F_2(D)&=\psi\left(M_2(D+a_2)\right)	
\label{eqF2explicit}\\
F_3(D)&=\frac{d\psi}{ds_2}\left(M_2(D+a_2)\right)
\label{eqF3explicit}
\end{align}
Since $M_2$ and $\psi$ are given explicitly by Eq.~\ref{eqM2explicit} and Eq.~\ref{eqpsiexplicit}, respectively, 
the functions $F_2(D)$ and $F_3(D)$ are given explicitly with respect to the biological parameters in Eq.~\ref{monMonod}.
%=============== ITEM9
The functions $F_2(D)$ and $F_3(D)$ are defined for $D$ such that $s_2^0(D)<M_2(D)<s_2^1(D)$, that is to say, for $D$ such that  
$$\frac{L_0(D+a_0)}{m_0-D-a_0}<\frac{K_2(D+a_2)}{m_2-D-a_2}<\frac{K_i(m_1-D-a_1)}{D+a_1}$$
%=============== ITEM10
We denote by 
\begin{equation}
I_2=\{D\in I_1: s_2^0(D)<M_2(D)<s_2^1(D)\}
\label{eqI2}
\end{equation}
the subset of $I_1$ on which $F_2(D)$ and $F_3(D)$ are defined.
For all  for $D\in I_2$, $F_1(D)\leq F_2(D)$.
%=============== ITEM11
The equality $F_1(D)= F_2(D)$ holds if, and only if, $M_2(D+a_2)=\overline{s}_2(D)$ that is, 
	$\frac{d\psi}{ds_2}\left(M_2(D+a_2)\right)=0$. 
	Therefore, $F_1(D)= F_2(D)$ holds if, and only if, $F_3(D)=0$. 
	We define 
	$$I_3=\{D\in I_2 :F_3(D)<0\}$$
	
Since $D\mapsto s_2^0(D)$ is increasing and $D\mapsto s_2^1(D)$ is decreasing, and assuming $s_2^0(0)<s_2^1(0)$, the domain of definition $I_1$ of $F_1(D)$ is an interval $I_1=[0,D_1)$, where $D_1$ is the solution of $s_2^0(D)=s_2^1(D)$, see Fig. \ref{figI1I2}. A necessary condition of existence of  SS2 is $0<D<D_1$.}
	
For the domain of definition $I_2$ of $F_2(D)$, several cases can be distinguished. 
$I_2$ is an interval $I_2=[0,D_2)$, where $D_2$ is the solution of $M_2(D)=s_2^1(D)$, see Fig. \ref{figI1I2}(a), or the solution of equation  
$M_2(D)=s_2^0(D)$, see Fig. \ref{figI1I2}(b). $I_2$ is empty, see Fig. \ref{figI1I2}(c).  
$I_2$ is an interval $I_2=(D_{2min},D_{2max})$ where $D_{2min}$ and $D_{2max}$ are the solutions of $M_2(D)=s_2^0(D)$ and $M_2(D)=s_2^1(D)$ respectively, see Fig. \ref{figI1I2}(d). A necessary condition of existence of  SS3 is $D\in I_2$. 
Cases (a)--(d) are obtained with the numerical parameter values listed in Table \ref{InterestingCases} and \ref{InterestingCases(d)}.
%===================================
\begin{table}[ht]
\centering
\begin{tabular}{llllll}
\hline
& $K_{S,\mathrm{H_{2},c}}$ & $a_i$&$D_1$ & $D_2$ &$D_3$\\                                
\hline
(a)
&$1.0\times10^{-6}$
&$0.02$
&$0.432$
&$0.373$
&$0.058$
\\ 
&&
0
&$0.452$
&$0.393$
&$0.078$\\
(b)
&$4.0\times10^{-6}$ 
&$0.02$
&$0.329$
&$0.236$   
&$I_3=I_2$\\
&&
0 
&$0.349$
&$0.256$
&$I_3=I_2$        \\
(c)
&$7.0\times10^{-6}$  
&$0.02$
&$0.287$
&$I_2=\emptyset$
&\\
&&0   
& $0.303$
&$I_2=\emptyset$         \\
\hline
\end{tabular}
\caption{Parameter values for cases (a), (b) and (c) of Fig. \ref{figI1I2}.
Unspecified parameter values are as in Table \ref{table0}. 
The table gives the values of $D_1$, $D_2$ and $D_3$ where 
$I_1=[0,D_1)$, $I_2=[0,D_2)$ and $I_3=[0,D_3)$}
\label{InterestingCases}
\end{table}
%========================
\begin{table}[ht]
\centering
\begin{tabular}{llllll}
\hline
 & $a_i$&$D_1$ & $D_{2min}$ & $D_{2max}$&$D_3$\\                                
\hline
(d)
&$0.02$
&$0.238$
&$0.101$
&$0.198$
&$0.161$\\
&0
&$0.258$
&$0.121$
&$0.218$
&$0.181$\\
\hline
\end{tabular}
\caption{Parameter values for case (d) of Fig. \ref{figI1I2}:  
$K_{S,\mathrm{H_{2},c}}=1.2\times10^{-5}$, $K_{S,\mathrm{H_{2}}}=0.5\times10^{-5}$ and $k_{m,\mathrm{H_{2}}}=5$. 
Unspecified parameter values are as in Table \ref{table0}. 
The table gives the values of $D_1$, $D_{2min}$, $D_{2max}$ and $D_3$ where 
$I_1=[0,D_1)$, $I_2=(D_{2min},D_{2max})$ and $I_3=(D_{2min},D_{3})$.}
\label{InterestingCases(d)}
\end{table}

We can state now the necessary and sufficient conditions of existence of SS2 and SS3. 
%==================
\begin{lemma}\label{lemma2}
If $\omega \geq 1$ then SS2 does not exist. If $\omega <1$ then SS2 exists if, and only if, $s_0^{\rm in}\geq F_1(D)$. Therefore, a necessary condition for the existence of SS2 is that $D\in I_1$, where $I_1$ is defined by Eq.~\ref{eqI1}. 
If $s_0^{\rm in}\geq F_1(D)$ then each solution $s_2$ of equation
%=================
\begin{equation}
\psi(s_2)=s_0^{\rm in},\quad s_2\in(s_2^0,s_2^1)
\label{eqSS5}
\end{equation}
gives a steady-state
${\rm SS2}= \left(x_0,x_1,x_2=0,s_0,s_1,s_2\right)$
where
\begin{align}
&
s_0=M_0(D+a_0,s_2),\quad s_1=M_1(D+a_1,s_2)
\nonumber\\
&
x_0=\frac{D}{D+a_0}(s_0^{{\rm in}}-s_0),
\quad
x_1=\frac{D}{D+a_1}(s_0^{{\rm in}}-s_0-s_1)
\label{eqSS5x0x1}
\end{align}
\end{lemma}

\begin{lemma}\label{lemma3}
If $\omega \geq 1$ then SS3 does not exist. If $\omega <1$ then SS3 exists if, and only if, 
$s_0^{{\rm in}}>F_2(D)$. Therefore, a necessary condition of existence of SS3 is that $D\in I_2$, where $I_2$ is defined by  
Eq.~\ref{eqI2}. If $s_0^{{\rm in}}>F_2(D)$ then the steady-state 
${\rm SS3}= 
\left(x_0,x_1,x_2,
s_0,s_1,s_2\right)
$  
is given by
\begin{align}
s_0&=M_0(D+a_0,M_2(D+a_2))\nonumber \\
s_1&=M_1(D+a_1,M_2(D+a_2))\nonumber \\
s_2&=M_2(D+a_2)
\label{eqSS8}
\end{align}	
and
\begin{align}
x_0&=\frac{D}{D+a_0}(s_0^{{\rm in}}-s_0),
\quad
x_1=\frac{D}{D+a_1}(s_0^{{\rm in}}-s_0-s_1)
\nonumber\\
x_2&=\frac{D}{D+a_2}\left((1-\omega)(s_0^{{\rm in}}-s_0)-s_1-s_2\right)
\label{eqSS8x2}
\end{align}
\end{lemma}

\begin{remark}
If $s_0^{{\rm in}}>F_1(D)$ then Eq.~\ref{eqSS5} has exactly two solutions denoted by $s_2^\flat$ and $s_2^\sharp$ and such that, see Fig. \ref{fig1}(c), 
$$s_2^0<s_2^\flat<\overline{s}_2<s_2^\sharp<s_2^1$$
If $s_0^{{\rm in}}=F_0(D)$ then $s_2^0<s_2^\flat=\overline{s}_2=s_2^\sharp<s_2^1$.
\end{remark}

To these solutions, $s_2^\flat$ and $s_2^\sharp$, correspond two steady-states of SS2, which are denoted by 
${\rm SS2}^\flat$ and ${\rm SS2}^\sharp$. 
These steady-states coalesce when $s_0^{{\rm in}}=F_0(D)$.

Since $F_1(D)\leq F_2(D)$, the condition $s_0^{{\rm in}}>F_2(D)$ for the existence of the positive steady-state SS3 implies that the condition 
$s_0^{{\rm in}}>F_2(D)$ for the existence of the two steady-states ${\rm SS2}^\flat$ and ${\rm SS2}^\sharp$ is satisfied. Therefore, if SS3 exists then  ${\rm SS2}^\flat$ and ${\rm SS2}^\sharp$ exist and are distinct.
If $s_0^{\rm in}=F_2(D)$ then ${\rm SS3}$ coalesces with ${\rm SS2}^\flat$ if $F_3(D)<0$, and with ${\rm SS2}^\sharp$ if $F_3(D)>0$, respectively. 
%The results are summarised in Table \ref{SS}.
%======================================

\begin{remark}\label{remSchin}
Using Eq.~\ref{eqS0in}, the conditions $s_0^{{\rm in}}>F_1(D)$ and $s_0^{{\rm in}}>F_2(D)$ of existence of the steady-state 
SS2 and SS3 respectively  are equivalent to the conditions
$$S_{\rm{ch,in}}>\frac{F_1(D)}{Y_3Y_4}\mbox{ and }S_{\rm{ch,in}}>\frac{F_2(D)}{Y_3Y_4}$$
respectively, expressed with respect to the inflowing concentration $S_{{\rm ch,in}}$.
\end{remark}

\begin{figure}[ht]
\setlength{\unitlength}{0.6cm}
\begin{center}
\begin{picture}(12,5)(0,-0.2)
\put(-2,-8){\rotatebox{0}{\includegraphics[scale=0.3]{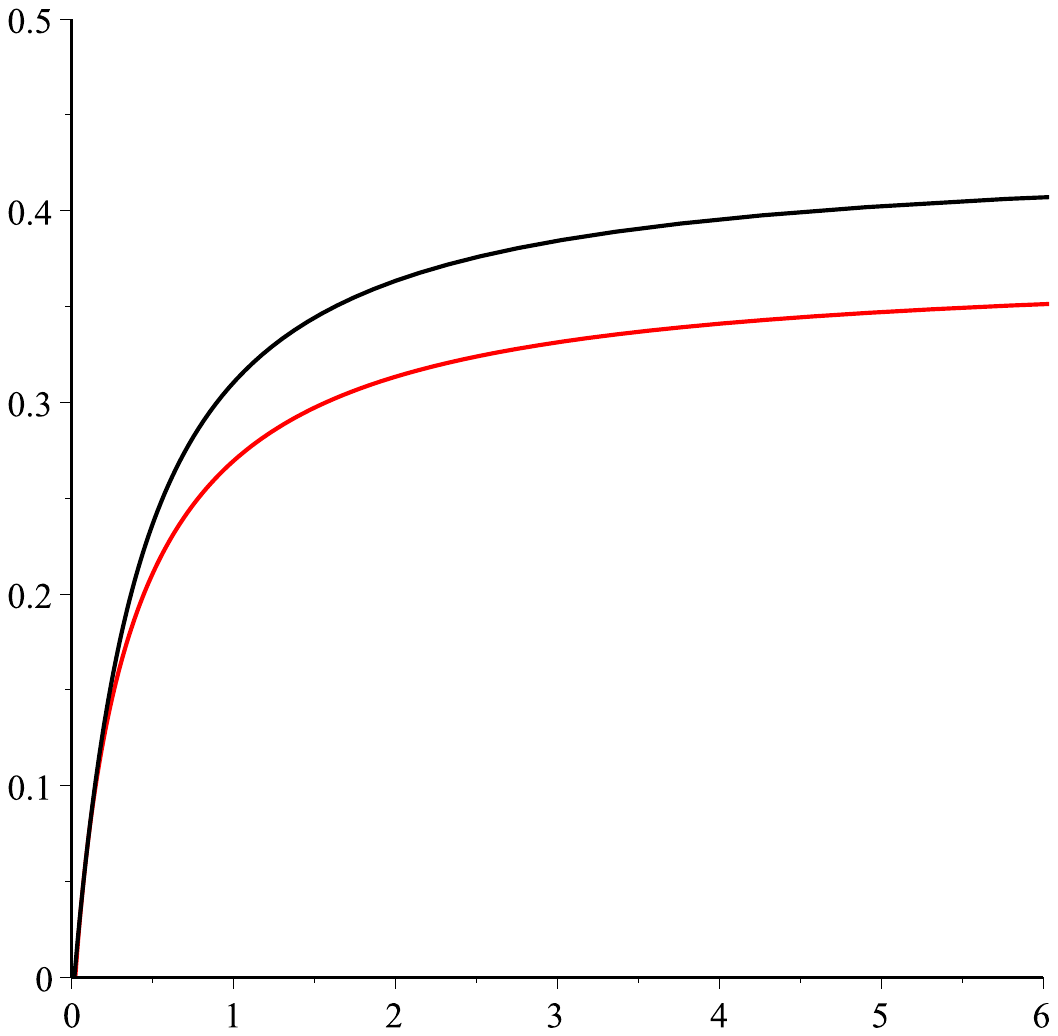}}}
\put(6,-8){\rotatebox{0}{\includegraphics[scale=0.3]{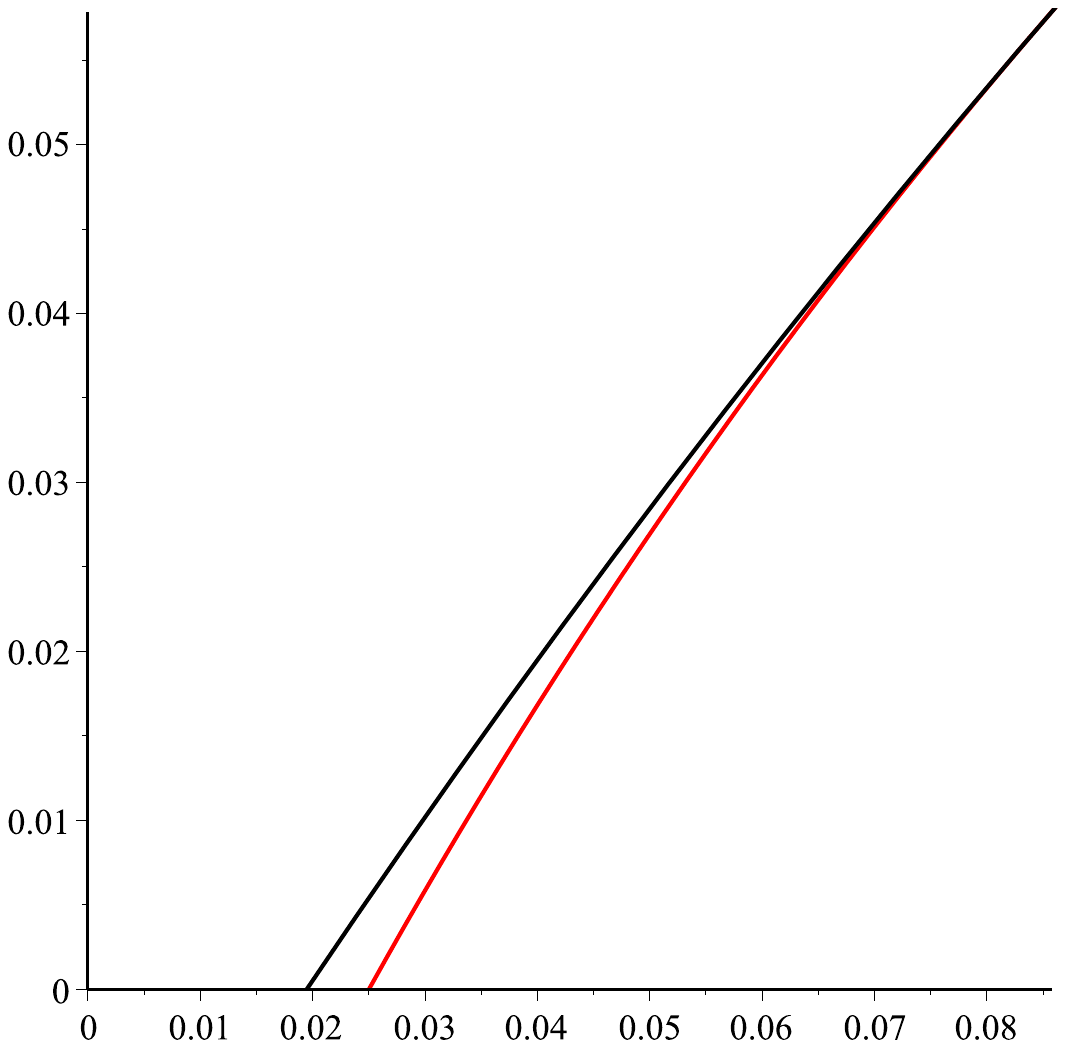}}}
%%%%%%%%%%%%
\put(-1.5,4.5){{\scriptsize{(i)}}}
\put(2,4.1){{\scriptsize{$\mathcal{J}_1$}}}
\put(4,3.5){{\scriptsize{$\mathcal{J}_2$}}}
\put(4,1){{\scriptsize{$\mathcal{J}_3$}}}
%===============================
\put(9,3.9){{\scriptsize{$\mathcal{J}_1$}}}
\put(11.5,2){{\scriptsize{$\mathcal{J}_3$}}}
\put(9,0.2){{\scriptsize{$\mathcal{J}_4$}}}
%=======================
\put(-1,2){{\rotatebox{90}{\scriptsize{$D$}}}}
\put(7,2){{\rotatebox{90}{\scriptsize{$D$}}}}
\put(2.2,-0.6){{\scriptsize{$S_{{\rm ch,in}}$}}}
\put(10.2,-0.6){{\scriptsize{$S_{{\rm ch,in}}$}}}
%====================
\put(0.6,4.2){{\color{black}\scriptsize{$\Gamma_1$}}}
\put(0.8,4.1){\vector(1,-1){0.6}}
\put(4.0,2.3){{\color{black}\scriptsize{$\Gamma_2$}}}
\put(4.1,2.65){\vector(0,1){0.6}}
\put(9,2){{\color{black}\scriptsize{$\Gamma_1$}}}
\put(9.15,1.95){\vector(1,-1){0.6}}
\put(10.5,0.4){{\color{black}\scriptsize{$\Gamma_2$}}}
\put(10.55,0.7){\vector(-1,1){0.6}}
%===============
\end{picture}
\vspace{0.5cm}

\begin{picture}(12,5)(0,-0.2)
\put(-2,-8){\rotatebox{0}{\includegraphics[scale=0.3]{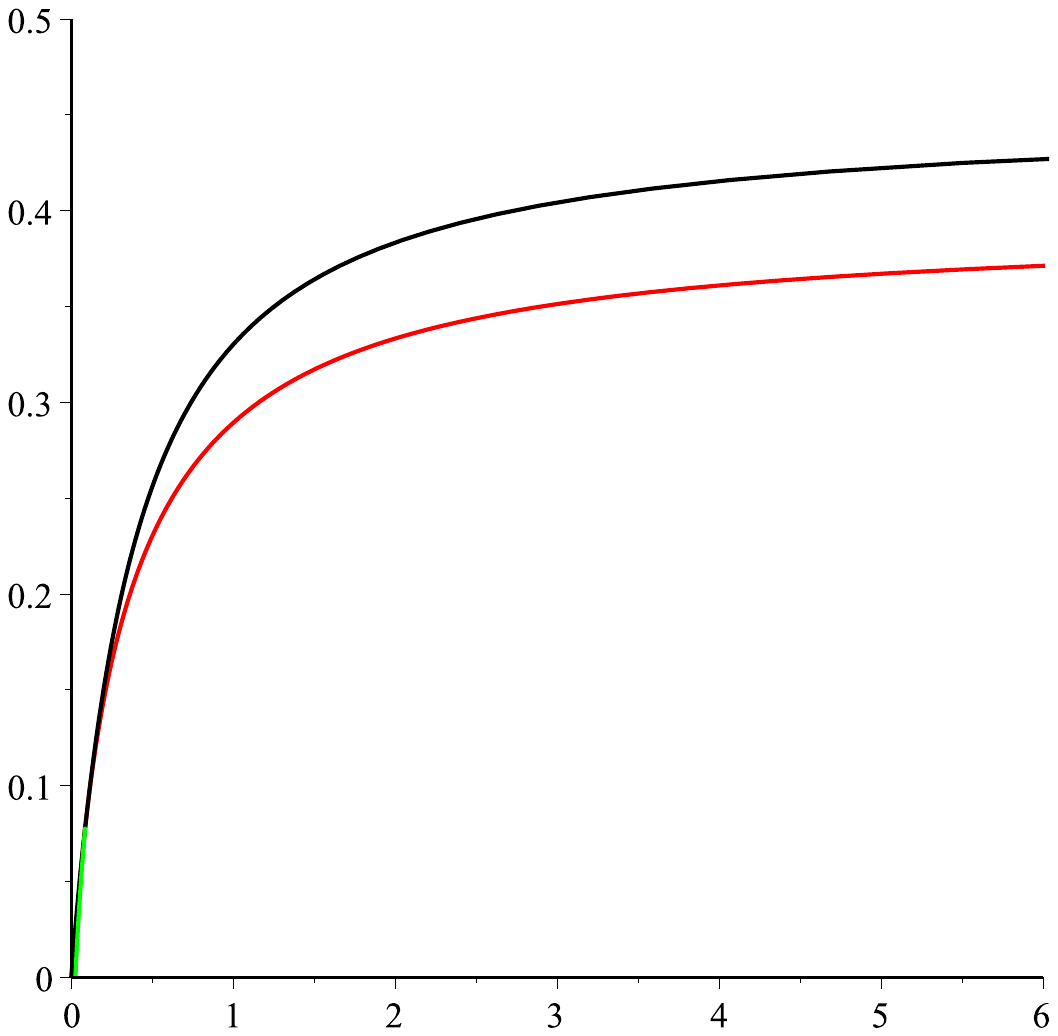}}}
\put(6,-8){\rotatebox{0}{\includegraphics[scale=0.3]{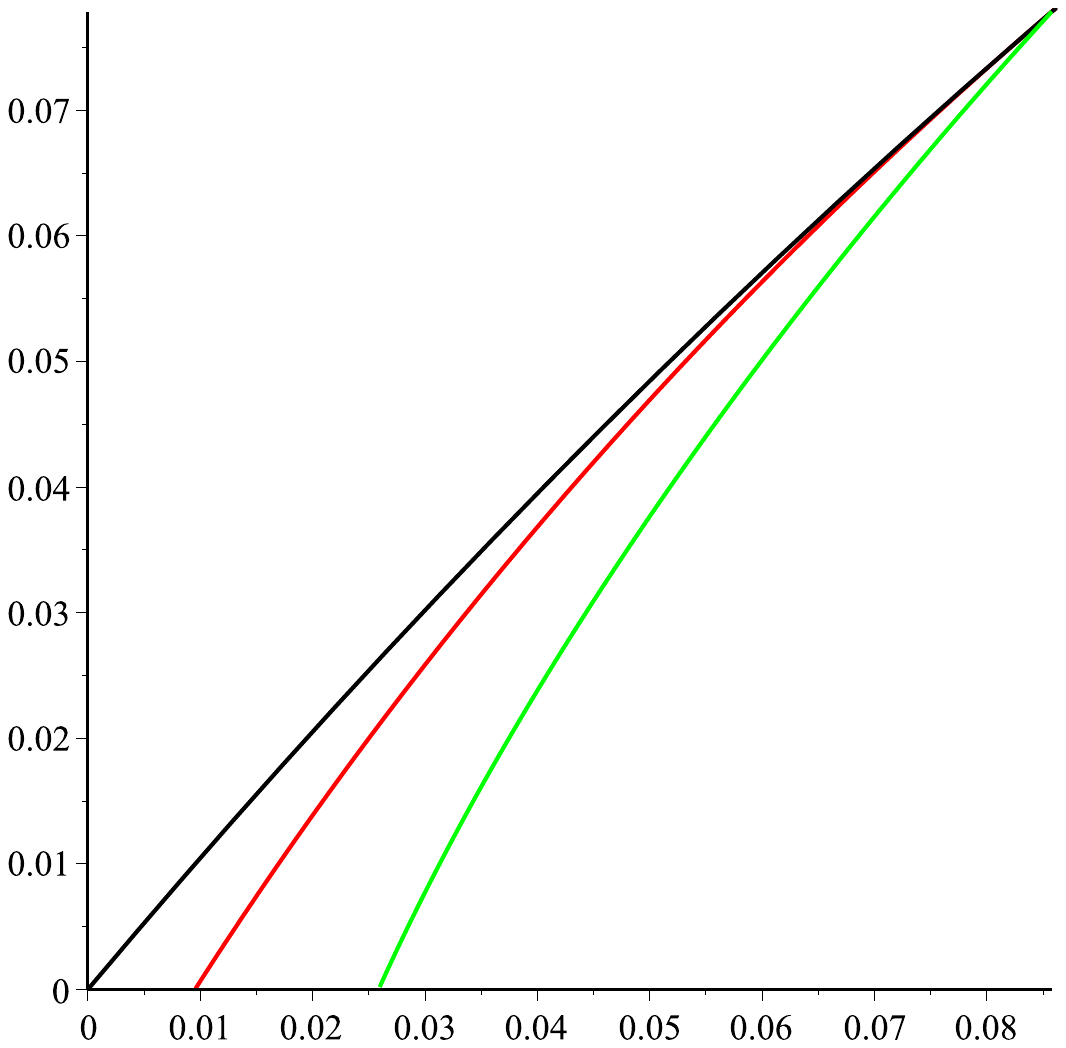}}}
%%%%%%%%%%%%
\put(-1.5,4.5){{\scriptsize{(ii)}}}
\put(2,4.2){{\scriptsize{$\mathcal{J}_1$}}}
\put(4,3.7){{\scriptsize{$\mathcal{J}_2$}}}
\put(4,1){{\scriptsize{$\mathcal{J}_3$}}}
\put(9,3.9){{\scriptsize{$\mathcal{J}_1$}}}
\put(11.3,1.5){{\scriptsize{$\mathcal{J}_3$}}}
\put(8,0.2){{\scriptsize{$\mathcal{J}_4$}}}
\put(8.7,0.2){{\scriptsize{$\mathcal{J}_5$}}}
%=======================
\put(-1,2){{\rotatebox{90}{\scriptsize{$D$}}}}
\put(7,2){{\rotatebox{90}{\scriptsize{$D$}}}}
\put(2.2,-0.6){{\scriptsize{$S_{{\rm ch,in}}$}}}
\put(10.2,-0.6){{\scriptsize{$S_{{\rm ch,in}}$}}}
%===============
\put(0.9,4.05){{\scriptsize{$\Gamma_1$}}}
\put(1.1,4){\vector(0,-1){0.4}}
\put(4.0,2.75){{\scriptsize{$\Gamma_2$}}}
\put(4.1,3.1){\vector(0,1){0.4}}
\put(0.3,0.7){{\scriptsize{$\Gamma_3$}}}
\put(0.25,0.65){\vector(-1,-1){0.6}}
\put(8.5,2.35){{\scriptsize{$\Gamma_1$}}}
\put(8.65,2.3){\vector(1,-1){0.6}}
\put(8.15,1.55){{\scriptsize{$\Gamma_2$}}}
\put(8.3,1.5){\vector(1,-1){0.6}}
\put(10.5,0.5){{\scriptsize{$\Gamma_3$}}}
\put(10.6,0.85){\vector(-1,1){0.6}}
%===========================
\end{picture}
\end{center}
\caption{The curves $\Gamma_1$ (black), $\Gamma_2$ (red) and $\Gamma_3$ (green) for case (a).
(i) : regions of steady-state existence, with maintenance. On the right,
a magnification for $0<D<D_3=0.058$ showing the region $\mathcal{J}_4$.
(ii) : regions of steady-state existence and their stability, without maintenance. On the right, 
a magnification for $0<D<D_3=0.078$ showing the regions $\mathcal{J}_4$ and $\mathcal{J}_5$.}
\label{fig3}
\end{figure}

Our aim now is to describe the operating diagram : The operating diagram shows how the system behaves when we vary the two control parameters 
$S_{\rm ch,in}$ and $D$ in Eqs.~\ref{model0a}-\ref{model0f}. 
According to Remark \ref{remSchin}, the curve $\Gamma_1$ of equation
\begin{equation}
S_{\rm ch,in}=\frac{1}{Y3Y4}F_1(D)
\label{eqGamma1}
\end{equation}
is the border to which SS2 exists, and 
the curve $\Gamma_2$ of equation
\begin{equation}
S_{\rm ch,in}=\frac{1}{Y3Y4}F_2(D)
\label{eqGamma2}
\end{equation}
is the border to which SS3 exists, see  Fig. \ref{fig3}. If we want to plot the operating diagram we must fix the values of the biological parameters. In the remainder of the Section we plot the operating diagrams corresponding to cases (a)--(d) depicted in Fig. \ref{figI1I2}.
%==============

%=============================

\subsection{Operating diagram: case (a)}
This case corresponds to the parameter values used by~\cite{wade15}. 
We have seen in Table \ref{InterestingCases} that the curves $\Gamma_1$ and $\Gamma_2$ are defined for
$D<D_1$ and $D<D_2$, respectively and that they are tangent for 
$D=D_3$, where $D_1=0.432$, $D_2=0.373$ and $D_3=0.058$. 
Therefore, they separate the operating plane $(S_{\rm ch,in},D)$ into four regions, as shown in Fig. \ref{fig3}(i), 
labelled $\mathcal{J}_1$, $\mathcal{J}_2$ and $\mathcal{J}_3$ and $\mathcal{J}_4$.
\begin{table}
\begin{center}
\begin{tabular}{ll}
\hline
Region  & Steady states\\
\hline
$\mathcal{J}_1$  & SS1 \\
$\mathcal{J}_2\cup\mathcal{J}_4$  & SS1, ${\rm SS2}^\flat$, ${\rm SS2}^\sharp$ \\
$\mathcal{J}_3$  & SS1, ${\rm SS2}^\flat$, ${\rm SS2}^\sharp$, SS3 \\
\hline
\end{tabular}
\caption{Existence of steady-states in the regions of the operating diagrams of Fig. \ref{fig3}(i) and Fig. \ref{fig3d}(i).}\label{ESS}
\end{center}
\end{table}
%======================================
\begin{figure}[ht]
\setlength{\unitlength}{0.6cm}
\begin{center}
\begin{picture}(12,5)(0,-0.2)
\put(-2,-8){\rotatebox{0}{\includegraphics[scale=0.3]{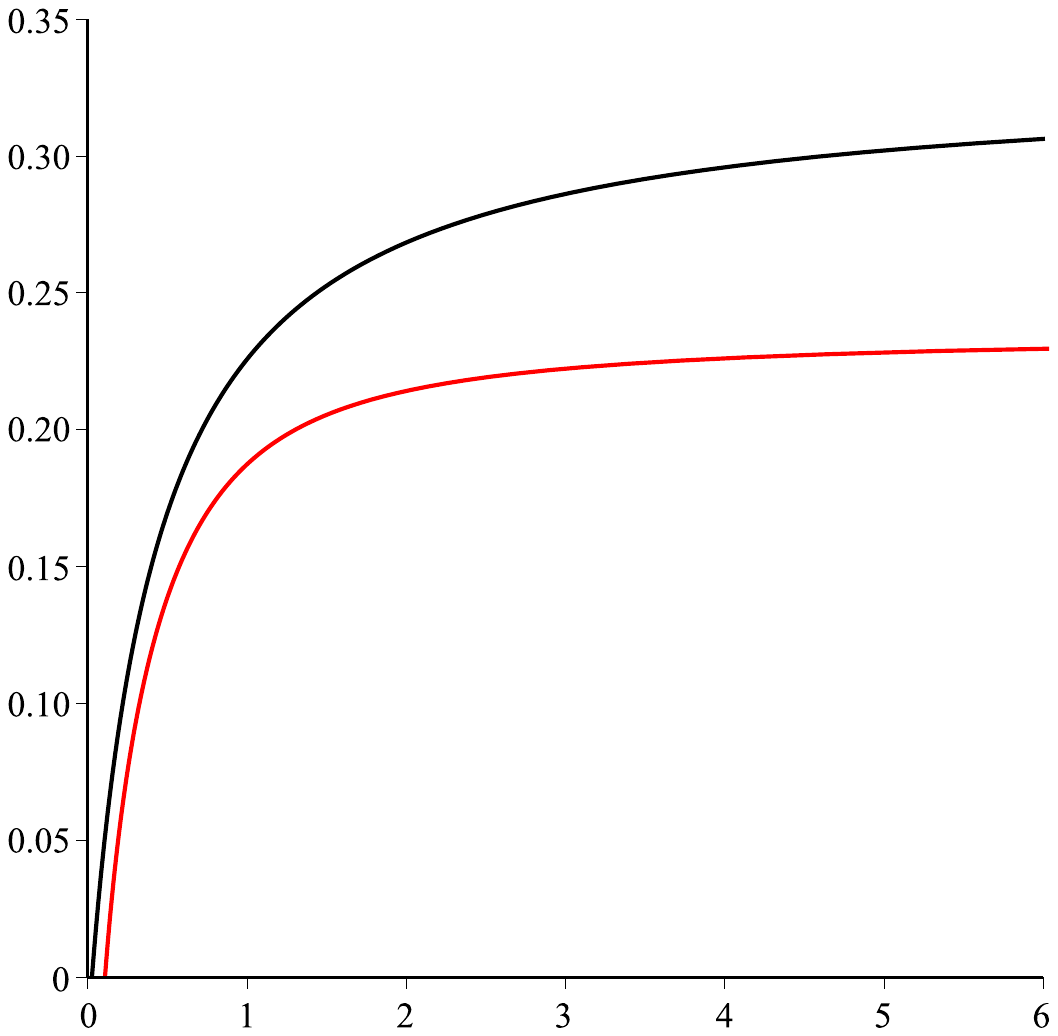}}}
\put(6,-8){\rotatebox{0}{\includegraphics[scale=0.3]{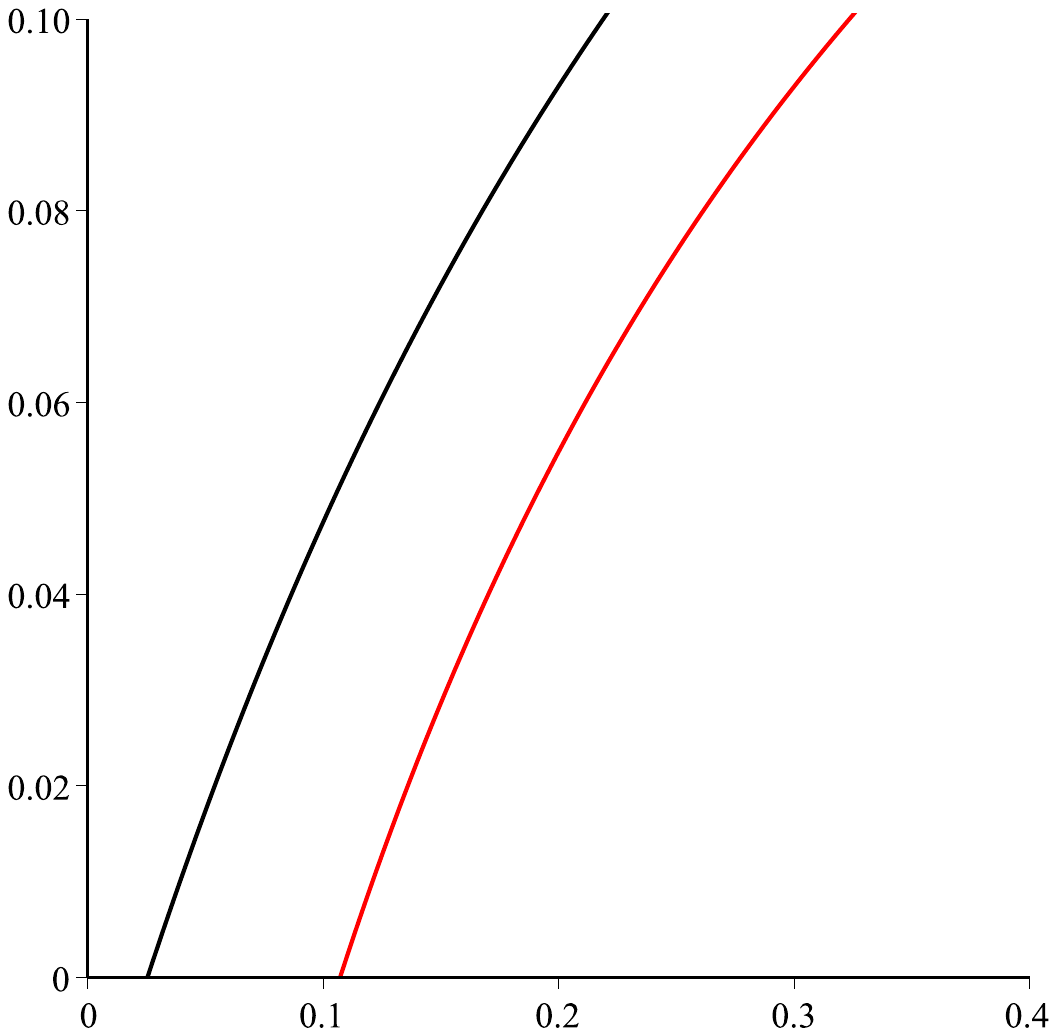}}}
%%%%%%%%%%%%
\put(-1.5,4.5){{\scriptsize{(i)}}}
\put(0.8,4.2){{\scriptsize{$\mathcal{J}_1$}}}
\put(4,3.5){{\scriptsize{$\mathcal{J}_4$}}}
\put(4,1.4){{\scriptsize{$\mathcal{J}_3$}}}
%===============================
\put(9,3.9){{\scriptsize{$\mathcal{J}_1$}}}
\put(11.5,0.6){{\scriptsize{$\mathcal{J}_3$}}}
\put(9.5,3){{\scriptsize{$\mathcal{J}_4$}}}
%=======================
\put(-1,2){{\rotatebox{90}{\scriptsize{$D$}}}}
\put(7,2){{\rotatebox{90}{\scriptsize{$D$}}}}
\put(2.2,-0.6){{\scriptsize{$S_{{\rm ch,in}}$}}}
\put(10.2,-0.6){{\scriptsize{$S_{{\rm ch,in}}$}}}
%====================
\put(4.1,4.6){{\color{black}\scriptsize{$\Gamma_1$}}}
\put(4.2,4.55){\vector(0,-1){0.4}}
\put(4.1,2.4){{\color{black}\scriptsize{$\Gamma_2$}}}
\put(4.2,2.65){\vector(0,1){0.4}}
\put(7.75,1.9){{\color{black}\scriptsize{$\Gamma_1$}}}
\put(8.3,2){\vector(1,0){0.4}}
\put(8.35,0.6){{\color{black}\scriptsize{$\Gamma_2$}}}
\put(8.85,0.7){\vector(1,0){0.4}}
%===============
\end{picture}
\vspace{0.5cm}

\begin{picture}(12,5)(0,-0.2)
\put(-2,-8){\rotatebox{0}{\includegraphics[scale=0.3]{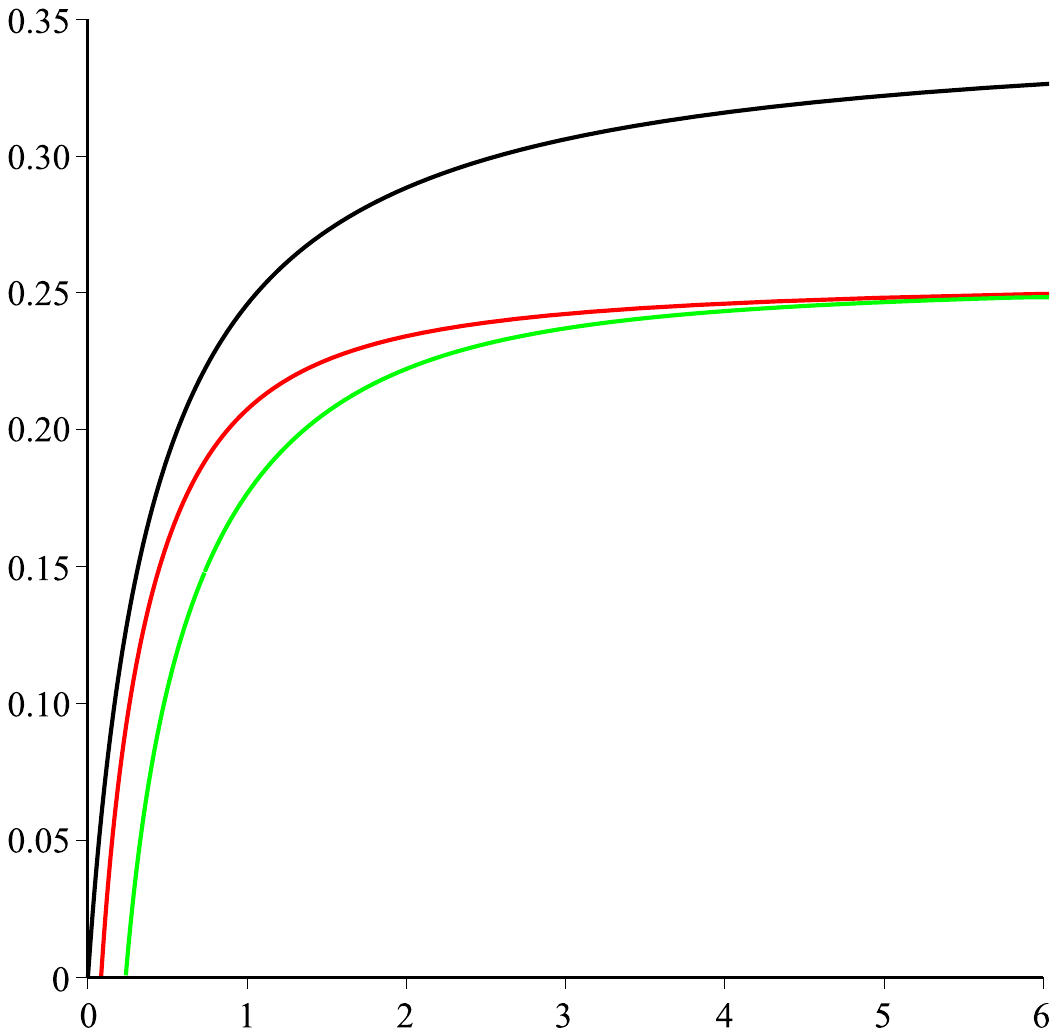}}}
\put(6,-8){\rotatebox{0}{\includegraphics[scale=0.3]{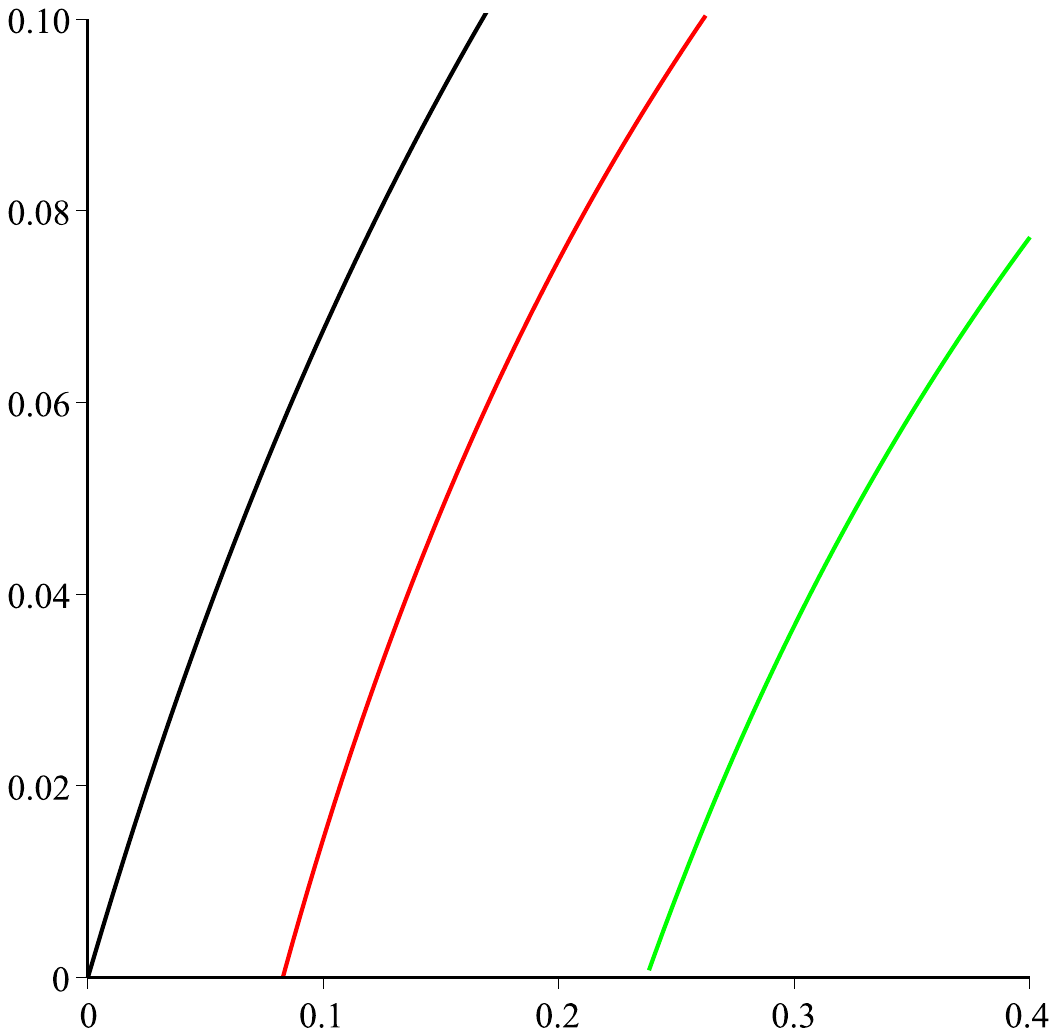}}}
%%%%%%%%%%%%
\put(-1.5,4.5){{\scriptsize{(ii)}}}
\put(0.2,4.2){{\scriptsize{$\mathcal{J}_1$}}}
\put(3.3,3.9){{\scriptsize{$\mathcal{J}_4$}}}
\put(3.5,1){{\scriptsize{$\mathcal{J}_3$}}}
\put(0.5,2.6){{\scriptsize{$\mathcal{J}_5$}}}
\put(8.2,3.9){{\scriptsize{$\mathcal{J}_1$}}}
\put(11.5,0.5){{\scriptsize{$\mathcal{J}_3$}}}
\put(8,0.5){{\scriptsize{$\mathcal{J}_4$}}}
\put(9.5,0.5){{\scriptsize{$\mathcal{J}_5$}}}
%=======================
\put(-1,2){{\rotatebox{90}{\scriptsize{$D$}}}}
\put(7,2){{\rotatebox{90}{\scriptsize{$D$}}}}
\put(2.2,-0.6){{\scriptsize{$S_{{\rm ch,in}}$}}}
\put(10.2,-0.6){{\scriptsize{$S_{{\rm ch,in}}$}}}
%===============
\put(0.8,3.4){{\scriptsize{$\Gamma_2$}}}
\put(0.8,3.4){\vector(0,-1){0.4}}
\put(4.1,3.65){{\scriptsize{$\Gamma_1$}}}
\put(4.3,4){\vector(0,1){0.4}}
\put(0.55,0.5){{\scriptsize{$\Gamma_3$}}}
\put(0.5,0.6){\vector(-1,0){0.4}}
\put(8.1,2.5){{\scriptsize{$\Gamma_1$}}}
\put(8.3,2.4){\vector(0,-1){0.4}}
\put(9.75,1.5){{\scriptsize{$\Gamma_2$}}}
\put(9.7,1.6){\vector(-1,0){0.4}}
\put(11,3){{\scriptsize{$\Gamma_3$}}}
\put(11.6,3.1){\vector(1,0){0.4}}
%===========================
\end{picture}
\end{center}
\caption{The curves $\Gamma_1$ (black), $\Gamma_2$ (red) and $\Gamma_3$ (green) for case (b). 
(i) : regions of steady-state existence, with maintenance. 
(ii) : regions of steady-state existence and their stability, without maintenance. 
On the right, a magnification for $0<D<0.1$}
\label{fig3b}
\end{figure} 
%=============================

The results are summarised in Table \ref{ESS}, which shows the existence of the steady-states SS1, SS2 and SS3 in the regions of the operating diagram in Fig. \ref{fig3}(i).

\subsection{Operating diagram: case (b)}\label{sec4.2}
This case corresponds to the parameter values used by~\cite{wade15}, except that $K_{S,\mathrm{H_{2},c}}$ is changed from $1.0\times10^{-6}$ to $4.0\times10^{-6}$. 
We have seen in Table \ref{InterestingCases} that the curves $\Gamma_1$ and $\Gamma_2$ are defined for $D<D_1$ and $D<D_2$, respectively and $F_1(D)<F_2(D)$ for all $D<D_2$, where $D_1=0.329$ and $D_2=0.236$. Therefore, they separate the operating plane $(S_{\rm ch,in},D)$ in three regions, as shown in Fig. \ref{fig3b}(i), 
labelled $\mathcal{J}_1$, $\mathcal{J}_3$ and $\mathcal{J}_4$.
\begin{table}
\begin{center}
\begin{tabular}{ll}
\hline
Region  & Steady states\\
\hline
$\mathcal{J}_1$  & SS1 \\
$\mathcal{J}_4$  & SS1, ${\rm SS2}^\flat$, ${\rm SS2}^\sharp$ \\
$\mathcal{J}_3$  & SS1, ${\rm SS2}^\flat$, ${\rm SS2}^\sharp$, SS3 \\
\hline
\end{tabular}
\caption{Existence of steady-states in the regions of the operating diagram of Fig \ref{fig3b}(i).}\label{ESSb}
\end{center}
\end{table}
%======================================

The results are summarised in Table \ref{ESSb}, which shows the existence of the steady-states SS1, SS2 and SS3 in the regions of the operating diagram in Fig. \ref{fig3b}(i). Note that the region $\mathcal{J}_2$ has disappeared.

%======================================
\begin{figure}[ht]
\setlength{\unitlength}{0.6cm}
\begin{center}
\begin{picture}(12,5)(0,-0.2)
\put(-2,-8){\rotatebox{0}{\includegraphics[scale=0.3]{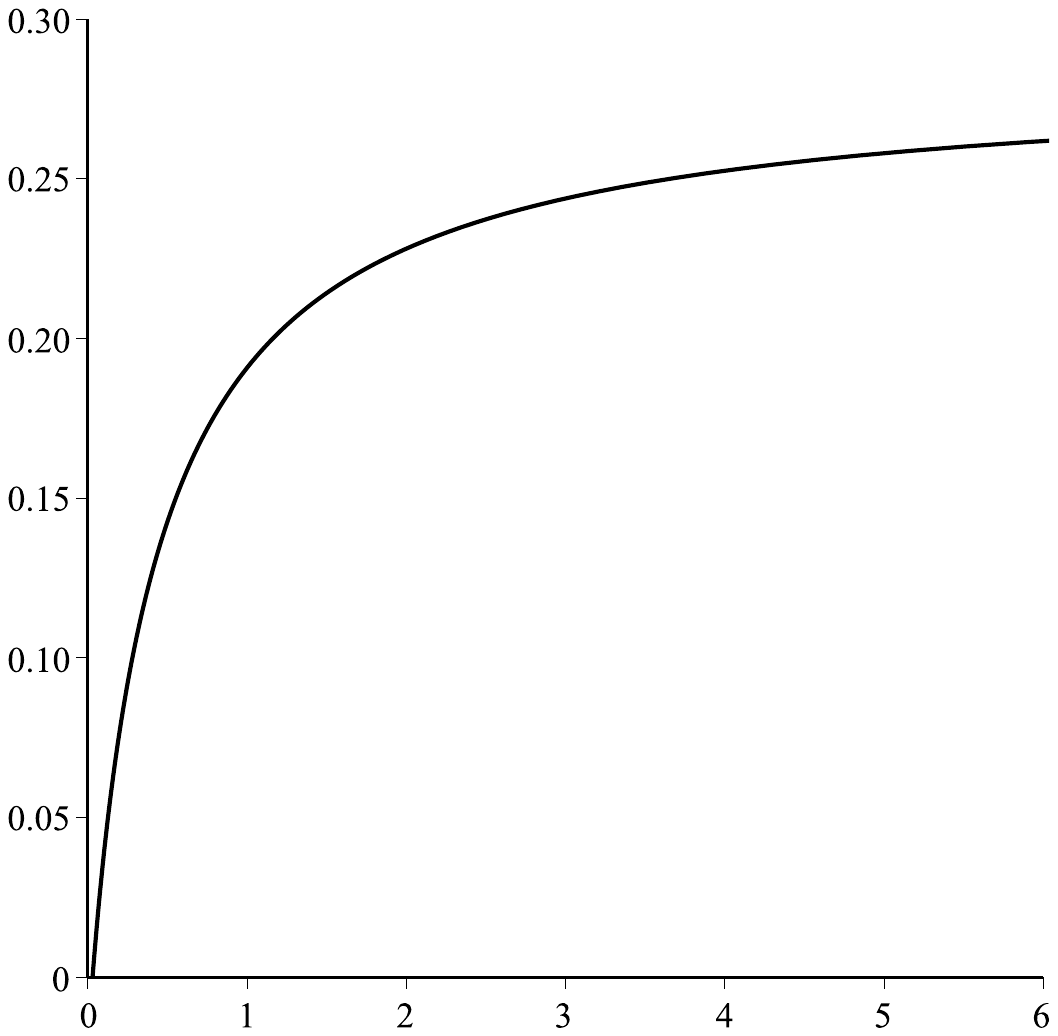}}}
\put(6,-8){\rotatebox{0}{\includegraphics[scale=0.3]{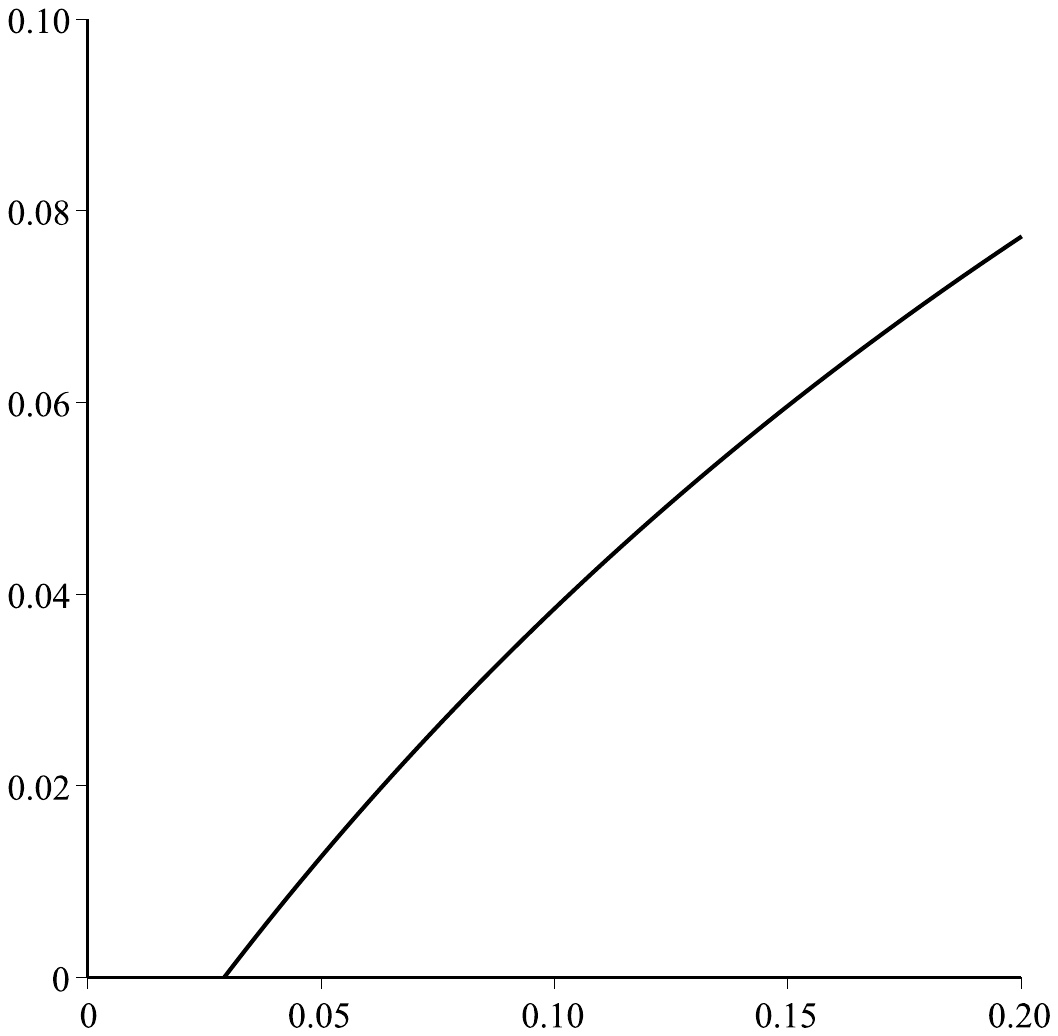}}}
%%%%%%%%%%%%
\put(-1.5,4.5){{\scriptsize{(i)}}}
\put(0.5,3.9){{\scriptsize{$\mathcal{J}_1$}}}
\put(3,2){{\scriptsize{$\mathcal{J}_4$}}}
%===============================
\put(9,3.9){{\scriptsize{$\mathcal{J}_1$}}}
\put(11,1){{\scriptsize{$\mathcal{J}_4$}}}
%=======================
\put(-1,2){{\rotatebox{90}{\scriptsize{$D$}}}}
\put(7,2){{\rotatebox{90}{\scriptsize{$D$}}}}
\put(2.2,-0.6){{\scriptsize{$S_{{\rm ch,in}}$}}}
\put(10.2,-0.6){{\scriptsize{$S_{{\rm ch,in}}$}}}
%====================
\put(4.2,3.35){{\color{black}\scriptsize{$\Gamma_1$}}}
\put(4.3,3.7){\vector(0,1){0.4}}
\put(9.2,1.9){{\color{black}\scriptsize{$\Gamma_1$}}}
\put(9.75,2){\vector(1,0){0.4}}
%===============
\end{picture}
\vspace{0.5cm}

\begin{picture}(12,5)(0,-0.2)
\put(-2,-8){\rotatebox{0}{\includegraphics[scale=0.3]{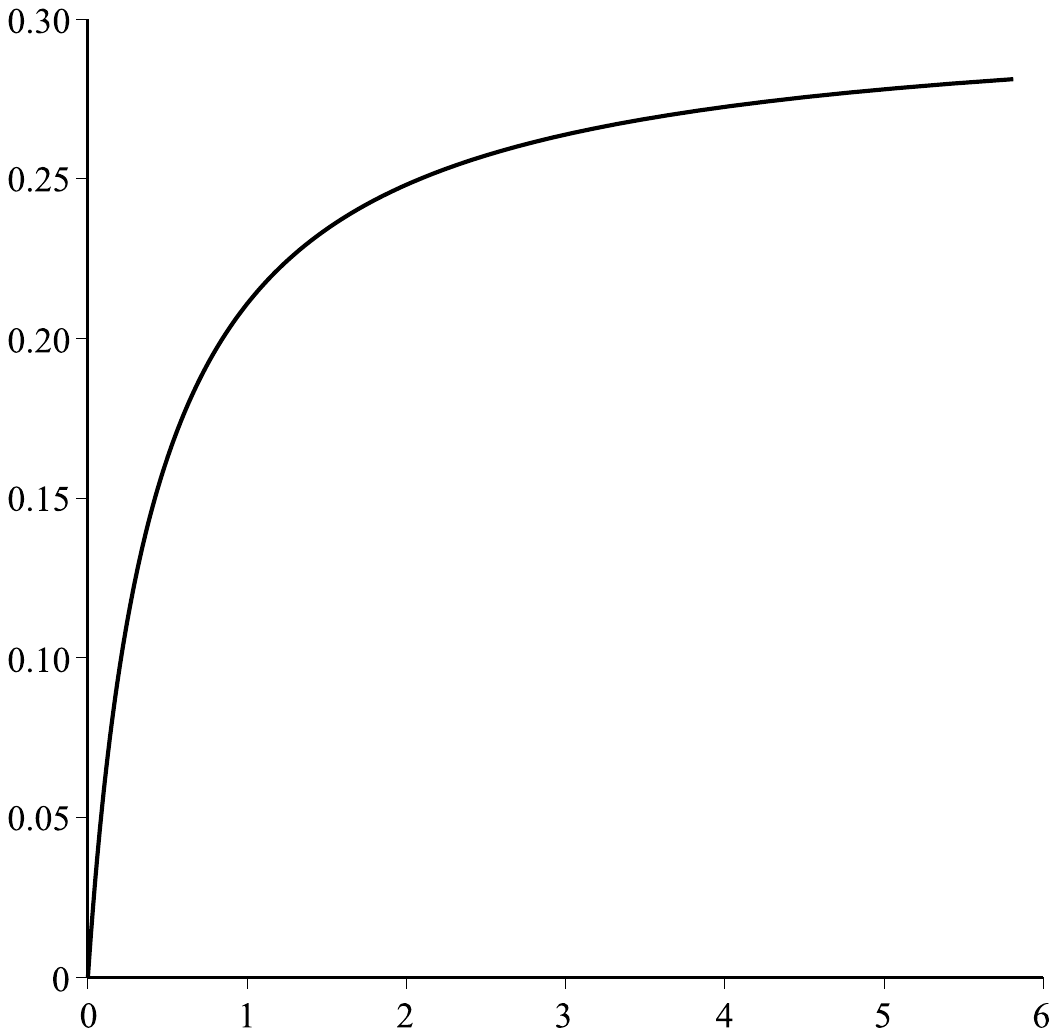}}}
\put(6,-8){\rotatebox{0}{\includegraphics[scale=0.3]{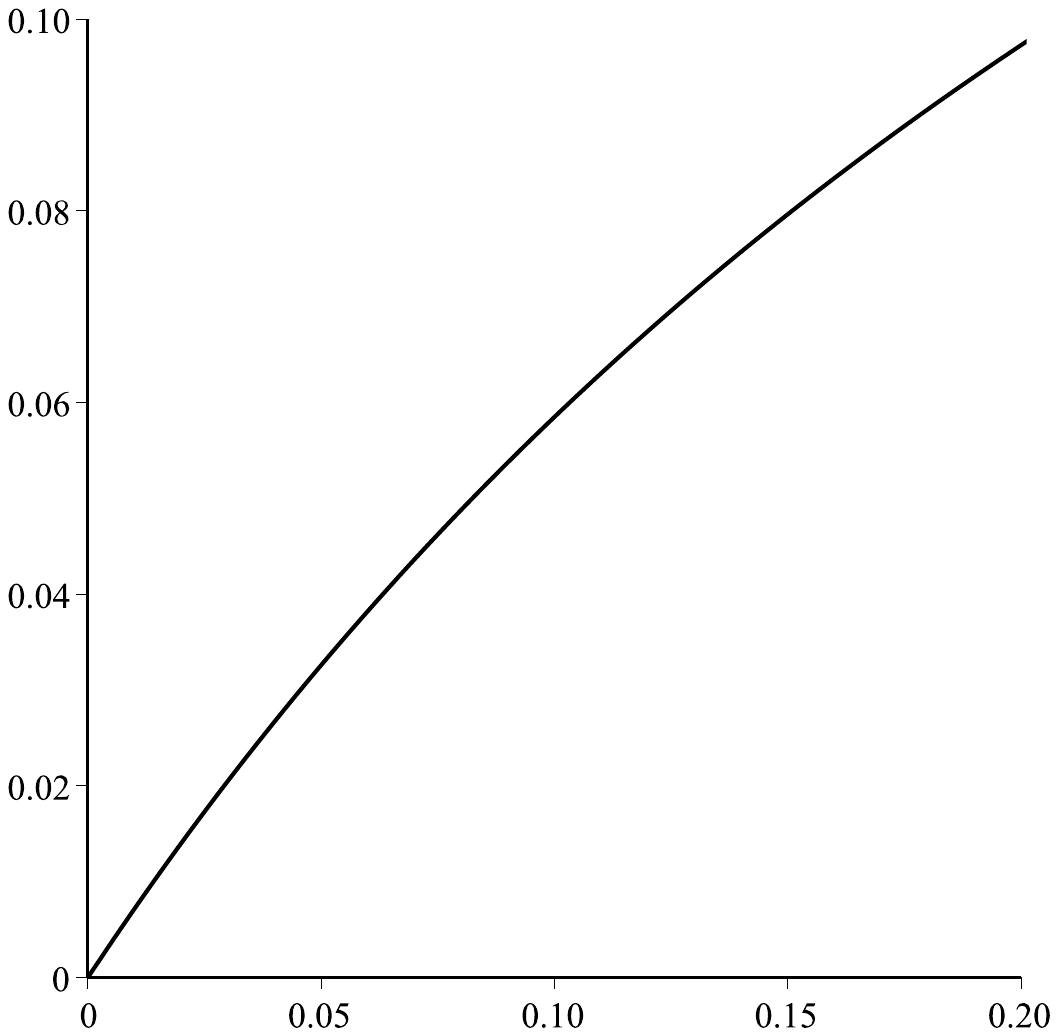}}}
%%%%%%%%%%%%
\put(-1.5,4.5){{\scriptsize{(ii)}}}
\put(0.5,3.9){{\scriptsize{$\mathcal{J}_1$}}}
\put(3,2){{\scriptsize{$\mathcal{J}_4$}}}
%===============================
\put(9,3.9){{\scriptsize{$\mathcal{J}_1$}}}
\put(11,1){{\scriptsize{$\mathcal{J}_4$}}}
%=======================
\put(-1,2){{\rotatebox{90}{\scriptsize{$D$}}}}
\put(7,2){{\rotatebox{90}{\scriptsize{$D$}}}}
\put(2.2,-0.6){{\scriptsize{$S_{{\rm ch,in}}$}}}
\put(10.2,-0.6){{\scriptsize{$S_{{\rm ch,in}}$}}}
%====================
\put(4.2,3.65){{\color{black}\scriptsize{$\Gamma_1$}}}
\put(4.3,4){\vector(0,1){0.4}}
\put(8.35,1.9){{\color{black}\scriptsize{$\Gamma_1$}}}
\put(8.9,2){\vector(1,0){0.4}}
%===========================
\end{picture}
\end{center}
\caption{The curve $\Gamma_1$ for case (c). 
(i) : regions of steady-state existence, with maintenance. 
(ii) : regions of steady-state existence. without maintenance and their stability.
On the right, a magnification for $0<D<0.1$.}
\label{fig3c}
\end{figure} 
%============================

\subsection{Operating diagram: case (c)}
This case corresponds to the parameter values used by~\cite{wade15}, except that $K_{S,\mathrm{H_{2},c}}$ is changed from $1.0\times10^{-6}$ to $7.0\times10^{-6}$. 
We have seen in Table \ref{InterestingCases} that the curve $\Gamma_1$ is defined for $D<D_1=0.287$ and that $I_2$ is empty so that SS3 does not exist. 
Therefore, $\Gamma_1$ separates the operating plane $(S_{\rm ch,in},D)$ in two regions, as shown in Fig. \ref{fig3c}(i), 
labelled $\mathcal{J}_1$  and $\mathcal{J}_4$.
\begin{table}
\begin{center}
\begin{tabular}{ll}
\hline
Region  & Steady states\\
\hline
$\mathcal{J}_1$  & SS1 \\
$\mathcal{J}_4$  & SS1, ${\rm SS2}^\flat$, ${\rm SS2}^\sharp$ \\
\hline
\end{tabular}
\caption{Existence of steady-states in the regions of the operating diagram of Fig \ref{fig3c}(i).}\label{ESSc}
\end{center}
\end{table}
%======================================

The results are summarised in Table \ref{ESSc}, which shows the existence of the steady-states SS1 and SS2 in the regions of the operating diagram in Fig. \ref{fig3c}(i). Note that the region $\mathcal{J}_3$ of existence of SS3 has disappeared.
%======================================
\begin{figure}[ht]
\setlength{\unitlength}{0.6cm}
\begin{center}
\begin{picture}(12,5)(0,-0.2)
\put(-2,-8){\rotatebox{0}{\includegraphics[scale=0.3]{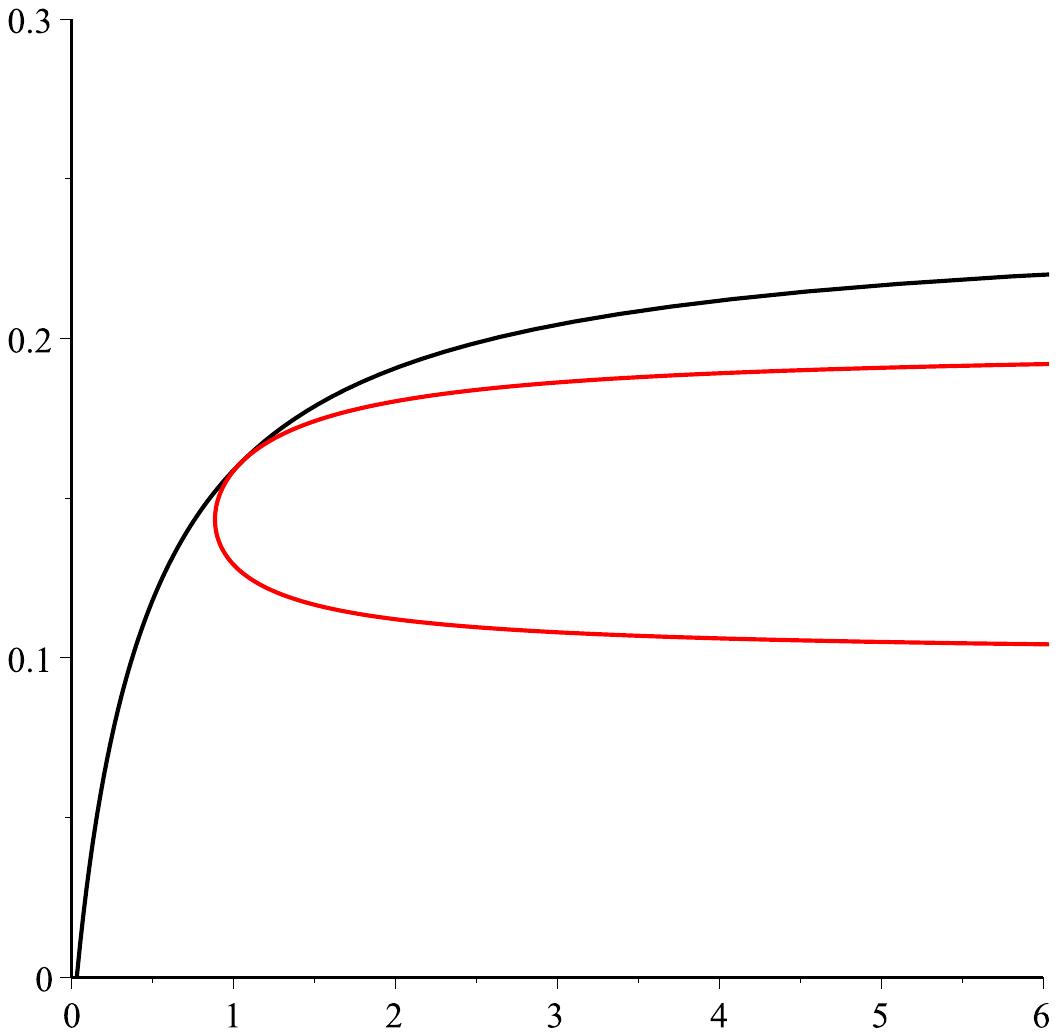}}}
\put(6,-8){\rotatebox{0}{\includegraphics[scale=0.3]{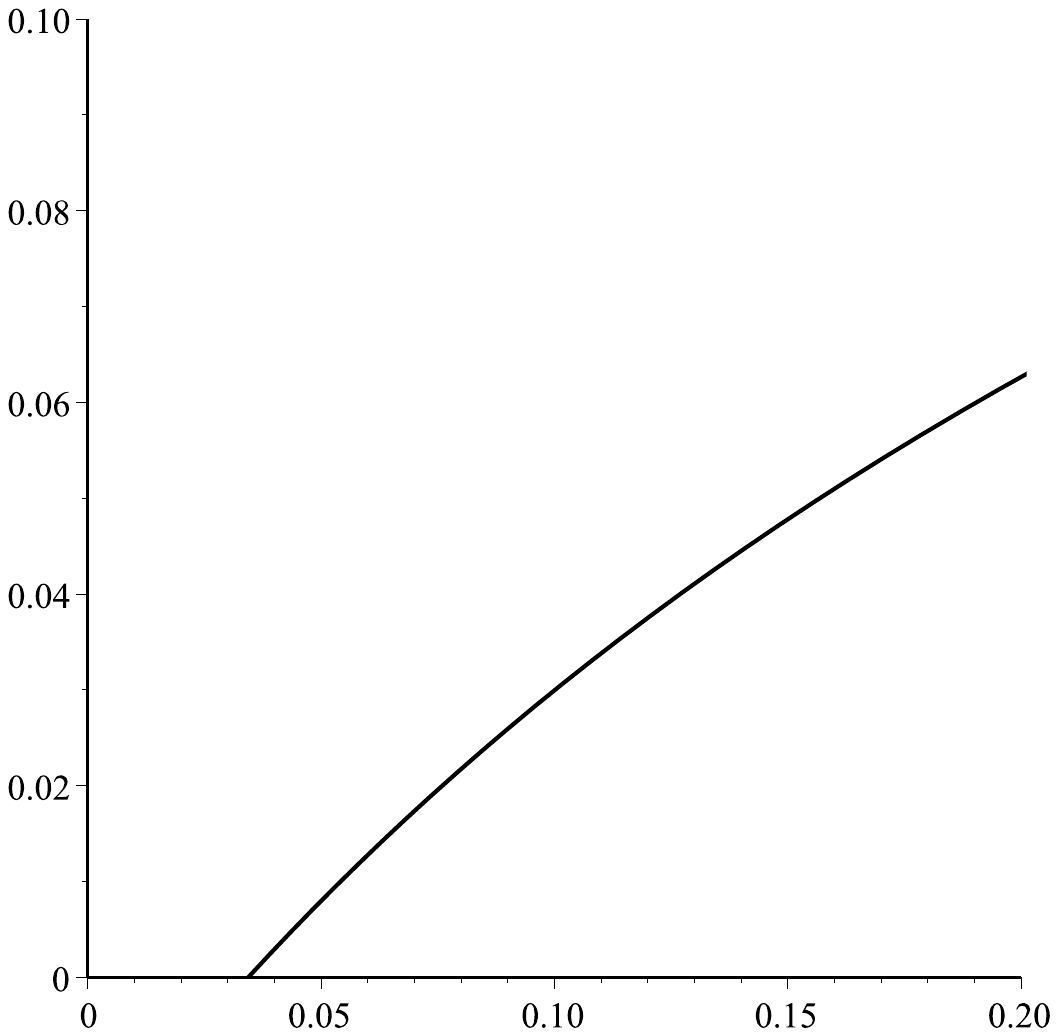}}}
%%%%%%%%%%%%
\put(-1.5,4.5){{\scriptsize{(i)}}}
\put(1,3.9){{\scriptsize{$\mathcal{J}_1$}}}
\put(4,3.2){{\scriptsize{$\mathcal{J}_2$}}}
\put(4,2){{\scriptsize{$\mathcal{J}_3$}}}
\put(4,0.7){{\scriptsize{$\mathcal{J}_4$}}}
%===============================
\put(9,3.9){{\scriptsize{$\mathcal{J}_1$}}}
\put(11.5,1){{\scriptsize{$\mathcal{J}_4$}}}
%=======================
\put(-1,2){{\rotatebox{90}{\scriptsize{$D$}}}}
\put(7,2){{\rotatebox{90}{\scriptsize{$D$}}}}
\put(2.2,-0.6){{\scriptsize{$S_{{\rm ch,in}}$}}}
\put(10.2,-0.6){{\scriptsize{$S_{{\rm ch,in}}$}}}
%====================
\put(4.2,4.1){{\color{black}\scriptsize{$\Gamma_1$}}}
\put(4.3,4){\vector(0,-1){0.4}}
\put(3.2,2.25){{\color{black}\scriptsize{$\Gamma_2$}}}
\put(3.3,2.6){\vector(0,1){0.4}}
\put(9.25,1.5){{\color{black}\scriptsize{$\Gamma_1$}}}
\put(9.8,1.6){\vector(1,0){0.4}}
%===============
\end{picture}
\vspace{0.5cm}

\begin{picture}(12,5)(0,-0.2)
\put(-2,-8){\rotatebox{0}{\includegraphics[scale=0.3]{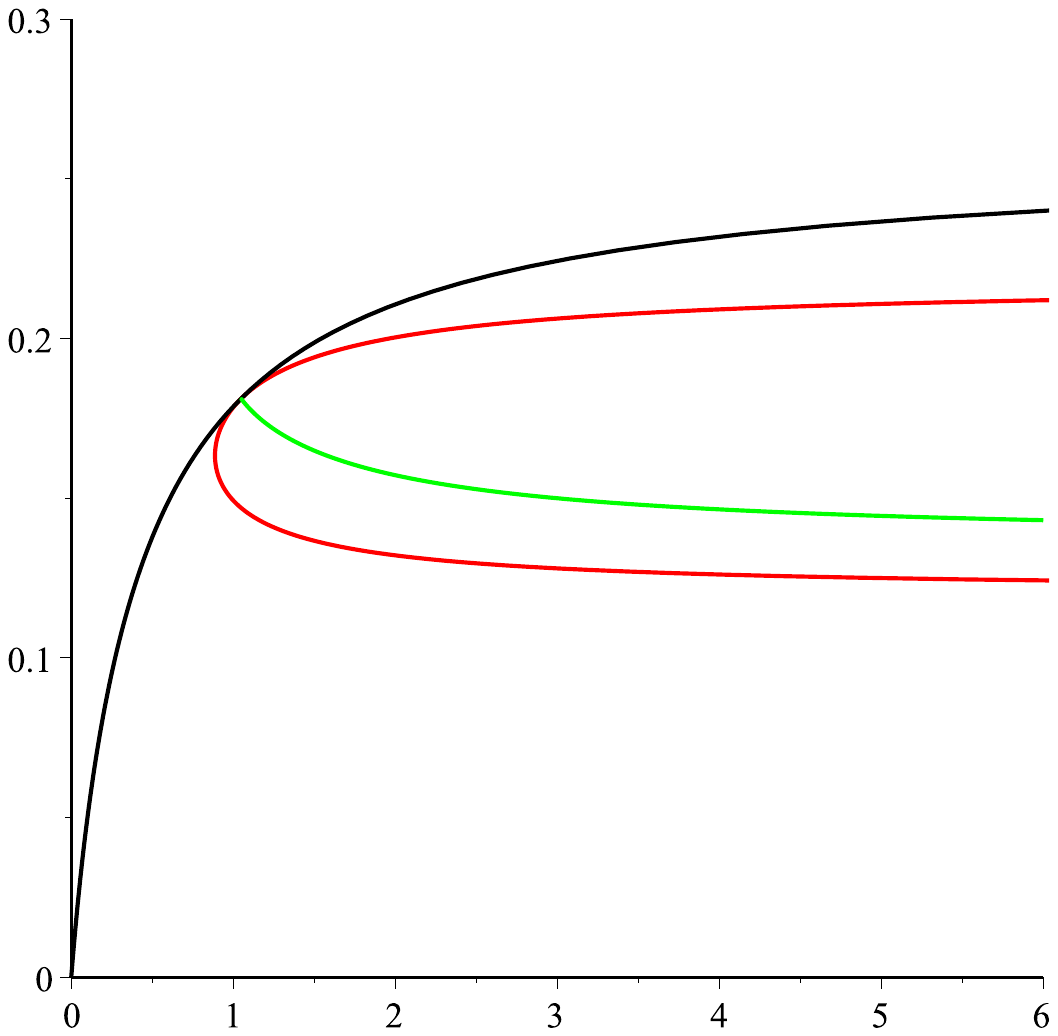}}}
\put(6,-8){\rotatebox{0}{\includegraphics[scale=0.3]{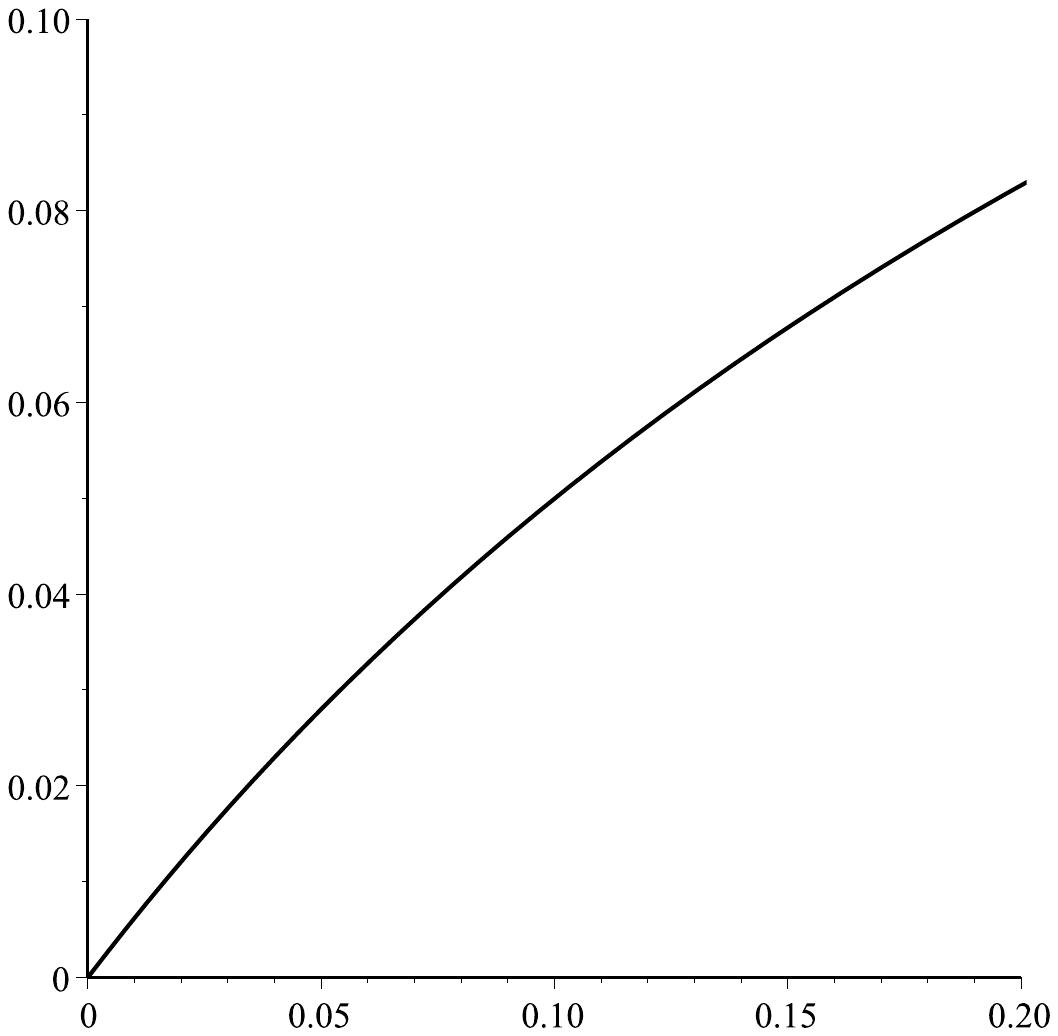}}}
%%%%%%%%%%%%
\put(-1.5,4.5){{\scriptsize{(ii)}}}
\put(1,3.9){{\scriptsize{$\mathcal{J}_1$}}}
\put(4,3.5){{\scriptsize{$\mathcal{J}_2$}}}
\put(1,2.2){{\scriptsize{$\mathcal{J}_5$}}}
\put(4,2.7){{\scriptsize{$\mathcal{J}_3$}}}
\put(4,0.7){{\scriptsize{$\mathcal{J}_4$}}}
%===============================
\put(9,3.9){{\scriptsize{$\mathcal{J}_1$}}}
\put(11,1){{\scriptsize{$\mathcal{J}_4$}}}
%=======================
\put(-1,2){{\rotatebox{90}{\scriptsize{$D$}}}}
\put(7,2){{\rotatebox{90}{\scriptsize{$D$}}}}
\put(2.2,-0.6){{\scriptsize{$S_{{\rm ch,in}}$}}}
\put(10.2,-0.6){{\scriptsize{$S_{{\rm ch,in}}$}}}
%====================
\put(4.2,4.4){{\scriptsize{$\Gamma_1$}}}
\put(4.3,4.3){\vector(0,-1){0.4}}
\put(1.2,1.35){{\scriptsize{$\Gamma_2$}}}
\put(1.3,1.7){\vector(0,1){0.4}}
\put(3.2,2.8){{\scriptsize{$\Gamma_3$}}}
\put(3.3,2.7){\vector(0,-1){0.4}}
\put(9.75,1.5){{\scriptsize{$\Gamma_1$}}}
\put(9.7,1.6){\vector(-1,0){0.4}}
%===========================
\end{picture}
\end{center}
\caption{The curves $\Gamma_1$ (black), $\Gamma_2$ (red) and $\Gamma_3$ (green) for case (d). 
(i) : regions of steady-state existence, with maintenance. 
(ii) : regions of steady-state existence and their stability, without maintenance.
On the right, a magnification for $0<D<0.1$.}
\label{fig3d}
\end{figure} 
%============================
\subsection{Operating diagram: case (d)}\label{sec4.4}
We end this discussion on the role of kinetic parameters by the presentation of this case, which presents a new behaviour that did not occur in the preceding cases: there exists $D_{2min}$ such that for $D<D_{2min}$ the system cannot have a positive steady-state SS3. This case corresponds to the parameter values used by~\cite{wade15}, except that three of them are changed as indicated in Table \ref{InterestingCases(d)}. This table shows that the curves $\Gamma_1$ and $\Gamma_2$ are defined for
$D<D_1$ and $D_{2min}<D<D_{2max}$ and that they are tangent for $D=D_3$, where $D_1=0.238$, $D_{2min}=0.101$, $D_{2max}=0.198$ and $D_3=0.161$.  
Therefore, $\Gamma_1$ and $\Gamma_2$ separate the operating plane $(S_{\rm ch,in},D)$ in four regions, as shown in Fig. \ref{fig3d}(i), labelled $\mathcal{J}_1$,  $\mathcal{J}_2$,  $\mathcal{J}_3$  and $\mathcal{J}_4$. 
The results are summarised in Table \ref{ESS}, which shows the existence of the steady-states SS1, SS2 and SS3 in the regions of the operating diagram in Fig. \ref{fig3d}(i).

\subsection{Stability of steady-states}\label{sec:ESS}

We know that SS1 is always stable. The analytical study of the stability of SS2 and SS3 is very difficult because the conditions for Routh-Hurwitz in the 6-dimensional case are intractable. For this reason we will consider in Section~\ref{sec:SS} the question of the stability only in the case without maintenance, since the system reduces to a 3-dimensional. The general case will be considered only numerically in Section~\ref{numanaly}.

%============================================
\section{Local stability without maintenance}\label{sec:SS}
%============================================
When maintenance is not considered in the model, the steady-states SS1, SS2 and SS3 are given by
\begin{enumerate}
	\item ${\rm SS1}= (0,0,0,s_0^{{\rm in}},0,0)$
	\item ${\rm SS2}= (x_0,x_1,0,s_0,s_1,s_2)$ where $s_2$ a solution of equation 
	$$s_0^{{\rm in}}=\psi(s_2)=M_0(D,s_2)+\frac{M_1(D,s_2)+s_2}{1-\omega}$$
	and 
	\begin{align}
	s_0&=M_0(D,s_2),\quad  s_1=M_1(D,s_2)\nonumber\\
	x_0&=s_0^{{\rm in}}-s_0,\quad x_1=s_0^{{\rm in}}-s_0-s_1
	\label{SS2WM}
	\end{align}
	\item ${\rm SS3}= (x_0,x_1,x_2,s_0,s_1,s_2)$ where 
	\begin{align}
	s_2&=M_2(D),\quad s_0=M_0(D,s_2),\quad s_1=M_1(D,s_2)
	\nonumber\\ 
	x_0&=s_0^{{\rm in}}-s_0,\quad
	x_1=s_0^{{\rm in}}-s_0-s_1
	\label{SS3WM}\\
	x_2&=(1-\omega)\left(s_0^{{\rm in}}-s_0\right)-s_1-s_2
	\nonumber
	\end{align}
\end{enumerate}

\begin{proposition}\label{prop1}
Let ${\rm SS2}=(x_0,x_1,0,s_0,s_1,s_2)$ be a steady-state. Then SS2 is stable if, and only if, $\mu_2(s_2)<D$ and 
$\frac{d\psi}{ds_2}>0$. 
\end{proposition}
Therefore, ${\rm SS2}^\flat$ is always unstable and ${\rm SS2}^\sharp$ is stable if, and only if, $\mu_2(s_2)<D$. This last condition is equivalent to $M_2(D)>s_2^\sharp$, which implies that $F_3(D)>0$. Hence, if SS3 exists then  ${\rm SS2}^\sharp$ is necessarily unstable. Therefore, ${\rm SS2}^\sharp$ is stable if, and only if, $F_3(D)>0$ and SS3 does not exist.

\begin{proposition}\label{prop2}
Let ${\rm SS3}=(x_0,x_1,x_2,s_0,s_1,s_2)$ be a steady-state. 
If $F_3(D)\geq 0$ then SS3 is stable as long as it exists.
If $F_3(D)<0$ then SS3 can be unstable. The instability of SS3 occurs in particular when $s_2$ is sufficiently close to $s_2^\flat$, that is to say SS3 is sufficiently close to ${\rm SS2}^\flat$.
\end{proposition}

The condition $F_3(D)\geq 0$ is equivalent to $\frac{d\psi}{ds_2}(M_2(D))\geq 0$, that is to say 
$s_2=M_2(D)\in[\overline{s}_2,s_2^\sharp)$.
If $\frac{d\psi}{ds_2}<0$, that is to say $s_2\in(s_2^\flat,\overline{s}_2)$, then SS3 can be unstable.  

When $D$ is such that $F_3(D)<0$, the determination of the boundary between the regions of stability and instability of SS3 needs to examine the Routh-Hurwitz condition of stability for SS3. For this purpose we define the following functions. 
Let ${\rm SS3}=(x_0,x_1,x_2,s_0,s_1,s_2)$ be a steady-state. Let
$$\scriptstyle{
E=\frac{\partial \mu_0}{\partial s_0},\quad
F=\frac{\partial \mu_0}{\partial s_2},\quad
G=\frac{\partial \mu_1}{\partial s_1},\quad
H=-\frac{\partial \mu_1}{\partial s_2},\quad
I=\frac{ d\mu_2}{d s_2}}
$$
evaluated at the steady-state SS3 defined by (\ref{SS3WM}), that is to say, for 
$$s_2=M_2(D),\quad s_0=M_0(D,s_2),\quad s_1=M_1(D,s_2)$$  
For $D\in I_3$ and $s_0^{\rm in}>F_2(D)$, we define 
{\small
\begin{align}
F_4\left(D,s_0^{\rm in}\right)& =
(EIx_0x_2+\left[E(G+H)-(1-\omega)FG\right] x_0x_1)f_2
\nonumber\\
&+(Ix_2+(G+H)x_1+\omega Fx_0)GIx_1x_2
\label{eqF4}
\end{align}}%
where $f_2=Ix_2+(G+H)x_1+(E+\omega F)x_0$.
Notice that to compute $F_4\left(D,s_0^{\rm in}\right)$, we must replace $x_0$, $x_1$, $x_2$, $s_0$, $s_1$ and $s_2$ by their values at SS3, given by (\ref{SS3WM}). Hence, this function depends on the operating parameters $D$ and $s_0^{\rm in}$. For each fixed $D\in I_3$, $F\left(D,s_0^{\rm in}\right)$ is polynomial in $s_0^{\rm in}$ of degree 3 and tends to $+\infty$ when 
$s_0^{\rm in}$ tends to $+\infty$. Therefore, it is necessarily positive for large enough $s_0^{\rm in}$. The values of the operating parameters $D$ and $s_0^{\rm in}$ for which $F\left(D,s_0^{\rm in}\right)$ is positive correspond to the stability of SS3 as shown in the following proposition.

\begin{proposition}\label{prop3}
Let ${\rm SS3}=(x_0,x_1,x_2,s_0,s_1,s_2)$ be a steady-state. 
If $F_3(D)<0$ then SS3 is stable if, and only if, $F_4\left(D,s_0^{\rm in}\right)>0$.
\end{proposition} 

The results on the existence of steady states (with or without maintenance) of Lemma \ref{lemma1}, Lemma \ref{lemma2} and Lemma \ref{lemma3}, and their stability (without maintenance) of Prop \ref{prop1}, Prop \ref{prop2} and Prop \ref{prop3}, are summarised in Table \ref{SSstability}.
\begin{table}
\begin{center}
\begin{tabular}{lll}
\hline
&Existence & Stability \\
\hline
SS1
&Always exists
&Always stable\\
${\rm SS2}^\flat$
&$s_0^{\rm in}>F_1(D)$
&Always unstable\\
${\rm SS2}^\sharp$    
&$s_0^{\rm in}>F_1(D)$
& $F_3(D)>0$ and $s_0^{in}<F_2(D)$\\
SS3          
&$s_0^{\rm in}>F_2(D)$
& $F_3(D)\geq 0$ or \\
&&$F_3(D)<0$ and $F_4\left(D,s_0^{\rm in}\right)>0$\\
\hline
\end{tabular}
\caption{Existence (with or without maintenance) and stability (without maintenance) of steady-states.}\label{SSstability}
\end{center}
\end{table}

\subsection{Operating diagram: case (a)}
This case corresponds to the parameter values used by~\cite{wade15} but without maintenance. 
We see from Table \ref{InterestingCases} that 
the curves $\Gamma_1$ and $\Gamma_2$ of the operating diagram, 
given by Eq. \ref{eqGamma1} and Eq. \ref{eqGamma2}, respectively, 
are defined now for $D<D_1=0.452$ and $D<D_2=0.393$, respectively and that they are tangent for $D=D_3=0.078$. 
Beside these curves, we plot also on the operating diagram of Fig.~\ref{fig3}(ii), 
the curve $\Gamma_3$ of equation
\begin{equation}
F_4\left(D,Y_3Y_4S_{{\rm ch,in}}\right)=0
\label{eqGamma3}
\end{equation}
According to Prop.~\ref{prop3}, this curve is defined for $D<D_3=0.078$ and it separates the region of existence of SS3 into two subregions labelled $\mathcal{J}_3$ and  $\mathcal{J}_5$, such that SS3 is stable in $\mathcal{J}_3$ and unstable in $\mathcal{J}_5$.
The other regions $\mathcal{J}_1$, 
$\mathcal{J}_2$ and  $\mathcal{J}_4$ are defined as in the previous section. The operating diagram is shown Fig. \ref{fig3}(ii). It looks very similar to Fig. \ref{fig3}(i), except near the origin, as it is indicated in the magnification for $0<D<D_3=0.078$.
From Table \ref{SSstability}, we deduce the following result 
\begin{proposition}\label{propODa}
Table \ref{SSS} shows the existence and stability of the steady-states SS1, SS2 and SS3 in the regions of the operating diagram 
in Fig. \ref{fig3}(ii).
\end{proposition}
%=============================================
\begin{table}
\begin{center}
\begin{tabular}{lcccc}
\hline
Region  & SS1 & ${\rm SS2}^\flat$ & ${\rm SS2}^\sharp$ & SS3\\
\hline
$\mathcal{J}_1$  & S &   &   &  \\
$\mathcal{J}_2$  & S & U & S &  \\
$\mathcal{J}_3$  & S & U & U & S \\
$\mathcal{J}_4$  & S & U & U &  \\
$\mathcal{J}_5$  & S & U & U & U\\
\hline
\end{tabular}
\caption{Existence and stability of steady-states in the regions of the operating diagrams of Fig. \ref{fig3}(ii) 
and Fig. \ref{fig3d}(ii).}\label{SSS}
\end{center}
\end{table}

%==================================  

\subsection{Operating diagram: case (b)} \label{sec5.2}
We see from Table \ref{InterestingCases} that the curves $\Gamma_1$ and $\Gamma_2$ are defined now for $D<D_1=0.349$ and $D<D_2=0.256$, respectively and that $F_1(D)<F_2(D)$ for all $D$. Beside these curves, we plot also on the operating diagram of Fig.  \ref{fig3b}(ii), the curve $\Gamma_3$ of equation (Eq.~\ref{eqGamma3})
which separates the region of existence of SS3  into two subregions labelled $\mathcal{J}_3$ and  $\mathcal{J}_5$, such that SS3 is stable in $\mathcal{J}_3$ and unstable in $\mathcal{J}_5$. Therefore, the curves $\Gamma_1$, $\Gamma_2$ and $\Gamma_3$ separate the operating plane $(S_{\rm ch,in},D)$ into four regions, as shown in Fig. \ref{fig3b}(ii), labelled $\mathcal{J}_1$,  $\mathcal{J}_3$,  $\mathcal{J}_4$  and $\mathcal{J}_5$.
\begin{proposition}\label{propODb}
Table \ref{SSSb} shows the existence and stability of the steady-states SS1, SS2 and SS3 in the regions of the operating diagram 
in Fig. \ref{fig3b}(ii)
\end{proposition}

%=============================================
\begin{table}
\begin{center}
\begin{tabular}{lcccc}
\hline
Region  & SS1 & ${\rm SS2}^\flat$ & ${\rm SS2}^\sharp$ & SS3\\
\hline
$\mathcal{J}_1$  & S &   &   &  \\
$\mathcal{J}_3$  & S & U & U & S \\
$\mathcal{J}_4$  & S & U & U &  \\
$\mathcal{J}_5$  & S & U & U & U\\
\hline
\end{tabular}
\caption{Existence and stability of steady-states in the regions of the operating diagram of Fig. \ref{fig3b}(ii).}\label{SSSb}
\end{center}
\end{table}
%===============================

\subsection{Operating diagram: case (c)}
We see from Table \ref{InterestingCases} that $\Gamma_1$ is defined for $D<D_1=0.303$ and that $I_2$ is empty so that SS3 does not exist. Therefore, $\Gamma_1$ separates the operating plane $(S_{\rm ch,in},D)$ into two regions, as shown in Fig. \ref{fig3c}(ii), labelled $\mathcal{J}_1$  and $\mathcal{J}_4$.
\begin{proposition}\label{propODc}
Table \ref{SSSc} shows the existence and stability of the steady-states SS1, SS2 and SS3 in the regions of the operating diagram 
in Fig. \ref{fig3c}(ii).
\end{proposition}

%============================================

%=============================================
\begin{table}
\begin{center}
\begin{tabular}{lcccc}
\hline
Region  & SS1 & ${\rm SS2}^\flat$ & ${\rm SS2}^\sharp$ & SS3\\
\hline
$\mathcal{J}_1$  & S &   &  \\
$\mathcal{J}_4$  & S & U & U &  \\
\hline
\end{tabular}
\caption{Existence and stability of steady-states in the regions of the operating diagram of Fig. \ref{fig3c}(ii).}\label{SSSc}
\end{center}
\end{table}
%===============================

\subsection{Operating diagram: case (d)} \label{sec5.4}

We see in Table \ref{InterestingCases(d)} that the curves $\Gamma_1$ and $\Gamma_2$ are defined for $D<D_1$ and $D_{2min}<D<D_{2max}$ and that they are tangent for $D=D_3$, where $D_1=0.258$ and $D_{2min}=0.121$, $D_{2max}=0.218$ and $D_3=0.181$. Beside these curves, we plot also on the operating diagram of Fig.  \ref{fig3d}(ii), the curve $\Gamma_3$ defined by Eq.~\ref{eqGamma3}, which separates the region of existence of SS3 into two subregions labelled $\mathcal{J}_3$ and  $\mathcal{J}_5$, such that SS3 is stable in $\mathcal{J}_3$ and unstable in $\mathcal{J}_5$. Therefore, the curves $\Gamma_1$, $\Gamma_2$ and $\Gamma_3$ separate the operating plane $(S_{\rm ch,in},D)$ into five regions, as shown in Fig. \ref{fig3d}(ii), labelled $\mathcal{J}_1$,  $\mathcal{J}_2$,  $\mathcal{J}_3$, $\mathcal{J}_4$  and $\mathcal{J}_5$. 
\begin{proposition}\label{propODd}
Table \ref{SSS} shows the existence and stability of the steady-states SS1, SS2 and SS3 in the regions of the operating diagram 
in Fig. \ref{fig3d}(ii).
\end{proposition}

%=============================
\section{Numerical analysis to confirm and extend the analytical results}\label{numanaly}
%=============================
The aim of this section is to study numerically (the method is explained in Appendix~\ref{nummeth}) the existence and stability of the steady-states SS2 and SS3. We obtain numerically the operating diagrams that were described in Sections \ref{sec:ES} and \ref{sec:SS}. The results in this section confirm the results on existence of the steady-states obtained in Section \ref{sec:ES} in the case with or without maintenance and the results of stability obtained in Section \ref{sec:SS} in the case without maintenance. These results permit also to elucidate the problem of the local stability of SS2 and SS3, which was left open in Section \ref{sec:ESS}.

\subsection{Operating diagram: case (a)}

\begin{figure}[ht!]
\setlength{\unitlength}{1cm}
\begin{picture}(6,4.2)(0.2,0.2)
\includegraphics[scale=0.14]{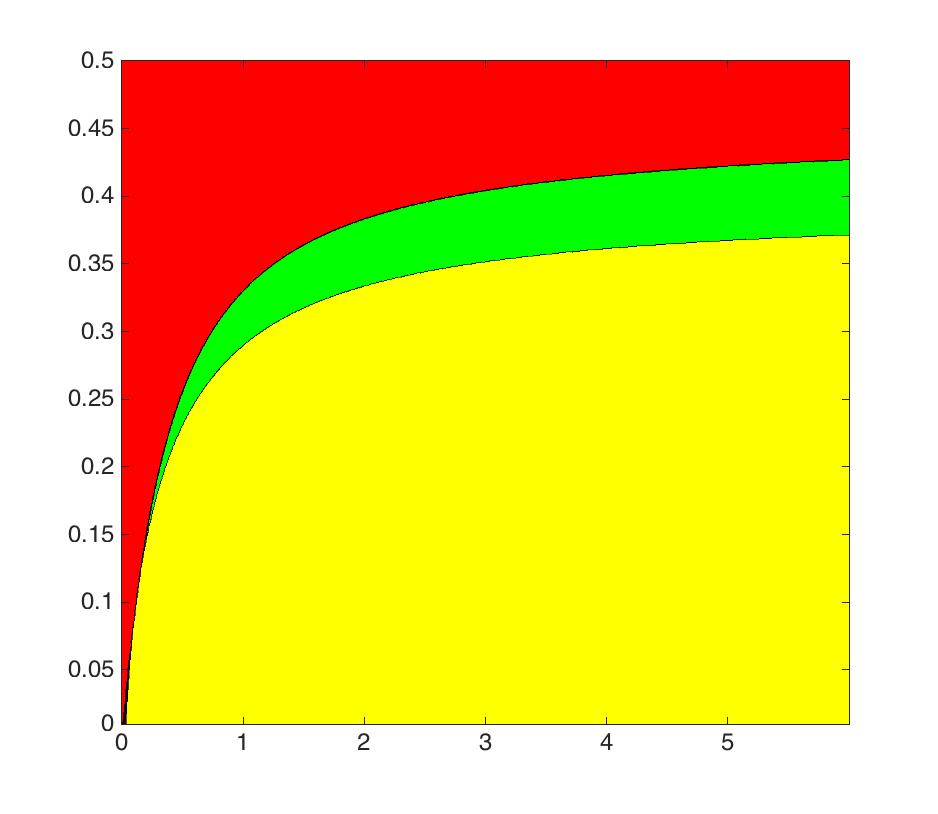}
%%%%%%%%%%%%
\put(-3.8,3){{\color{white}{$\mathcal{J}_1$}}}
\put(-2.2,2.85){{\color{black}{$\mathcal{J}_2$}}}
\put(-2,1.5){{\color{black}{$\mathcal{J}_3$}}}
\put(-3.6,1.2){{\color{black}{$\mathcal{J}_5$}}}
\put(-3.6,0.85){{\color{black}{$\mathcal{J}_4$}}}
\put(-3.6,1.1){\vector(-1,-1){0.4}}
\put(-3.65,0.85){\vector(-1,-1){0.4}}
\put(-4.6,3.6){{\color{black}{\scriptsize(i)}}}
\put(-4.6,2){{\rotatebox{90}{\scriptsize{$D$}}}}
\put(-2.5,0.1){{\scriptsize{$S_{{\rm ch,in}}$}}}
%===============
\end{picture}
\begin{picture}(0,4.2)(2,0.2)
\includegraphics[scale=0.14]{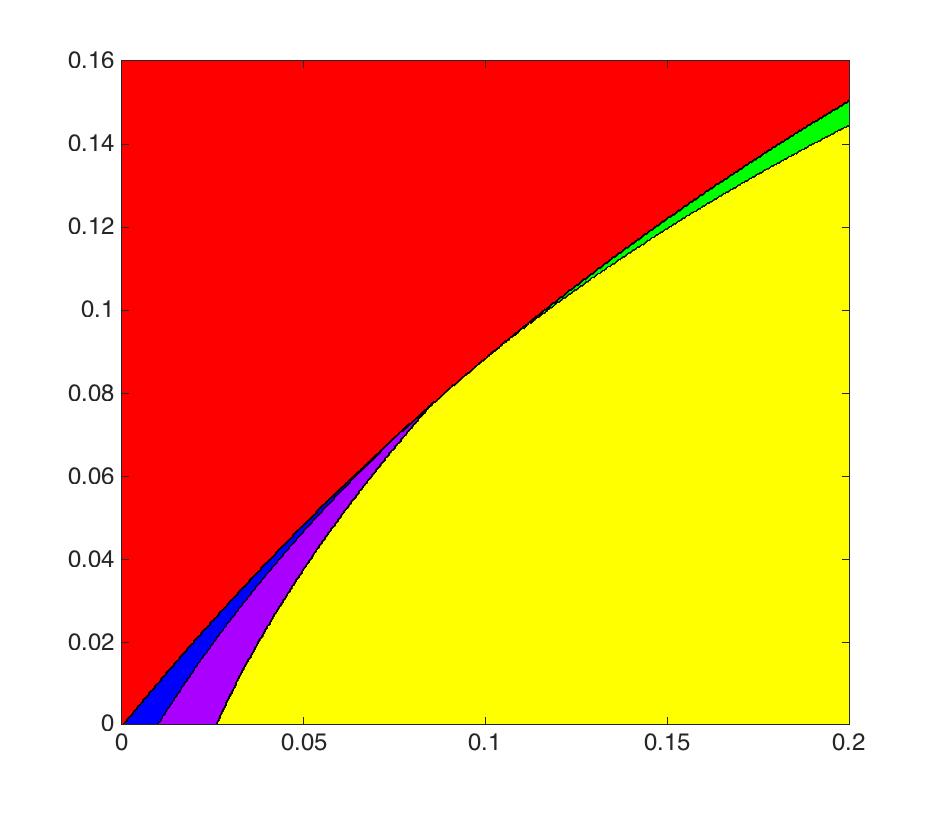}
%%%%%%%%%%%%
\put(-3.4,2.5){{\color{white}{$\mathcal{J}_1$}}}
\put(-1,2.5){{\color{black}{$\mathcal{J}_2$}}}
\put(-1.8,2){{\color{black}{$\mathcal{J}_3$}}}
\put(-4,1.1){{\color{white}{$\mathcal{J}_4$}}}
\put(-3,0.7){{\color{black}{$\mathcal{J}_5$}}}
\put(-4.7,3.6){{\color{black}{\scriptsize(ii)}}}
\put(-3,0.8){\vector(-1,0){0.4}}
\put(-0.8,2.8){\vector(0,0){0.4}}
\put(-3.8,1){{\color{white}{\vector(0,-1){0.4}}}}
\put(-4.6,2){{\rotatebox{90}{\scriptsize{$D$}}}}
\put(-2.5,0.1){{\scriptsize{$S_{{\rm ch,in}}$}}}
%===============
\end{picture}
\caption{Numerical analysis for the existence and stability of steady-states for case (a), without maintenance. On the right, a magnification for $0<D<0.16$.}
\label{od_ins_full}
\end{figure}
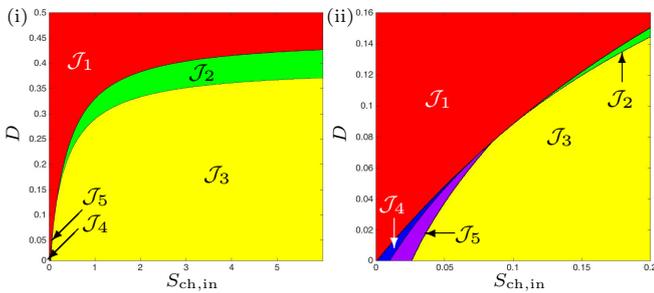

\begin{figure}[ht]
\setlength{\unitlength}{1.6cm}
\begin{center}
\begin{picture}(6,4.2)(0,0.3)
\includegraphics[scale=0.27]{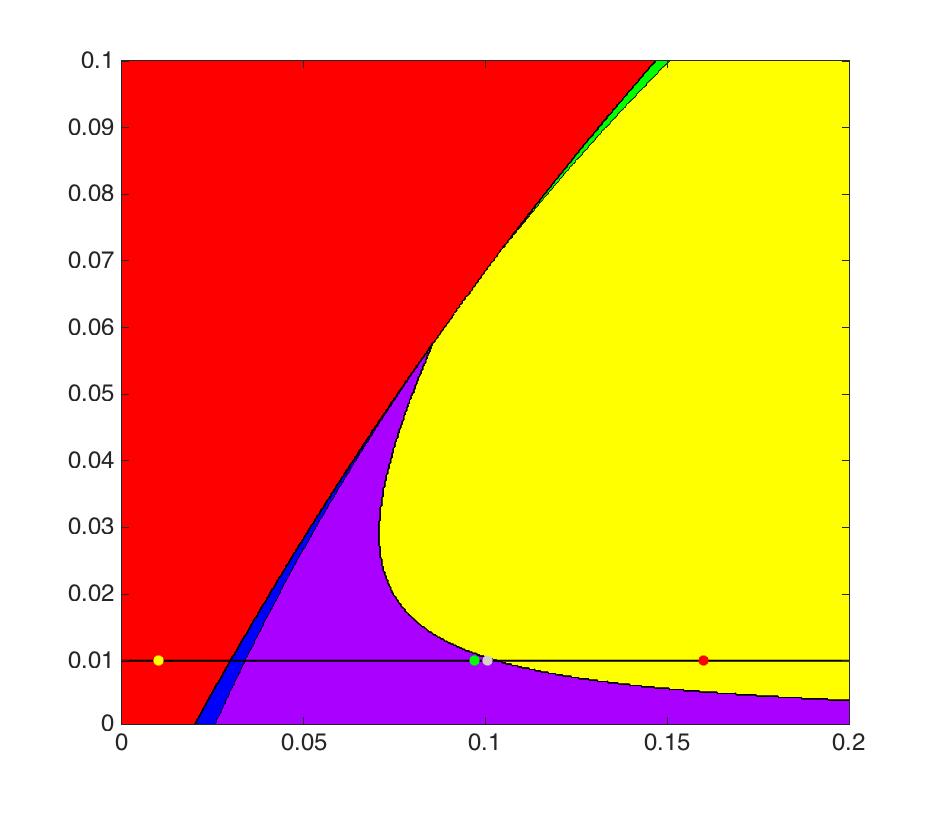}
%%%%%%%%%%%%
\put(-4.8,0.7){{\color{yellow}{($\alpha$)}}}
\put(-3,0.7){{\color{green}{($\beta$)}}}
\put(-2.7,0.7){{\color{white}{($\gamma$)}}}
\put(-1.53,0.7){{\color{red}{($\delta$)}}}
\put(-4,3){{\color{white}{$\mathcal{J}_1$}}}
\put(-1.8,4.2){{\color{black}{$\mathcal{J}_2$}}}
\put(-1.6,2.5){{\color{black}{$\mathcal{J}_3$}}}
\put(-4.4,0.6){{\color{white}{$\mathcal{J}_4$}}}
\put(-3.7,1.1){{\color{white}{$\mathcal{J}_5$}}}
\put(-5.5,2.5){{\rotatebox{90}{{$D$}}}}
\put(-3,0.2){{{$S_{{\rm ch,in}}$}}}
%===============
\end{picture}
\end{center}
\caption{Numerical analysis for the existence and stability of steady-states for case (a), with maintenance. This is a magnification for $0<D<0.1$, showing the presence and extent of region $\mathcal{J}_5$ undetectable by the analytical method. The coordinates labelled $(\alpha) - (\delta)$ are subsequently used to simulate the system dynamics, as shown in the proceeding Fig.~\ref{oscil}.}
\label{od_ins}
\end{figure}
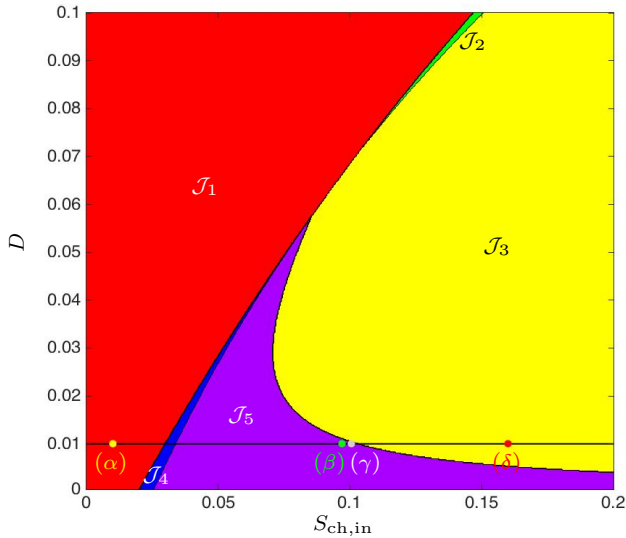

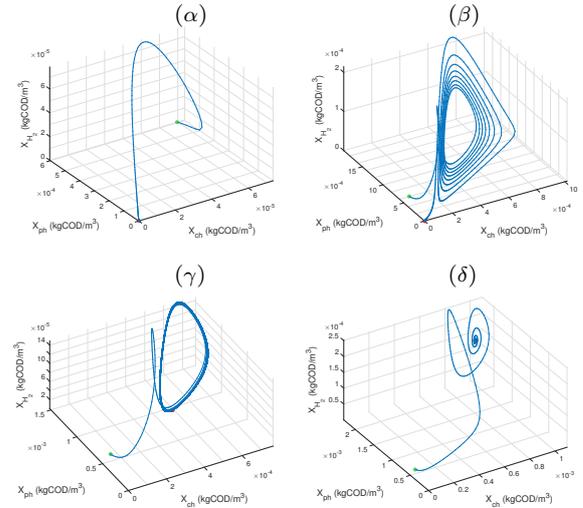
\begin{figure}[ht!]
\setlength{\unitlength}{1.5cm}
\begin{center}
\begin{picture}(5,4.5)(0.2,0.6)
\includegraphics[scale=0.17]{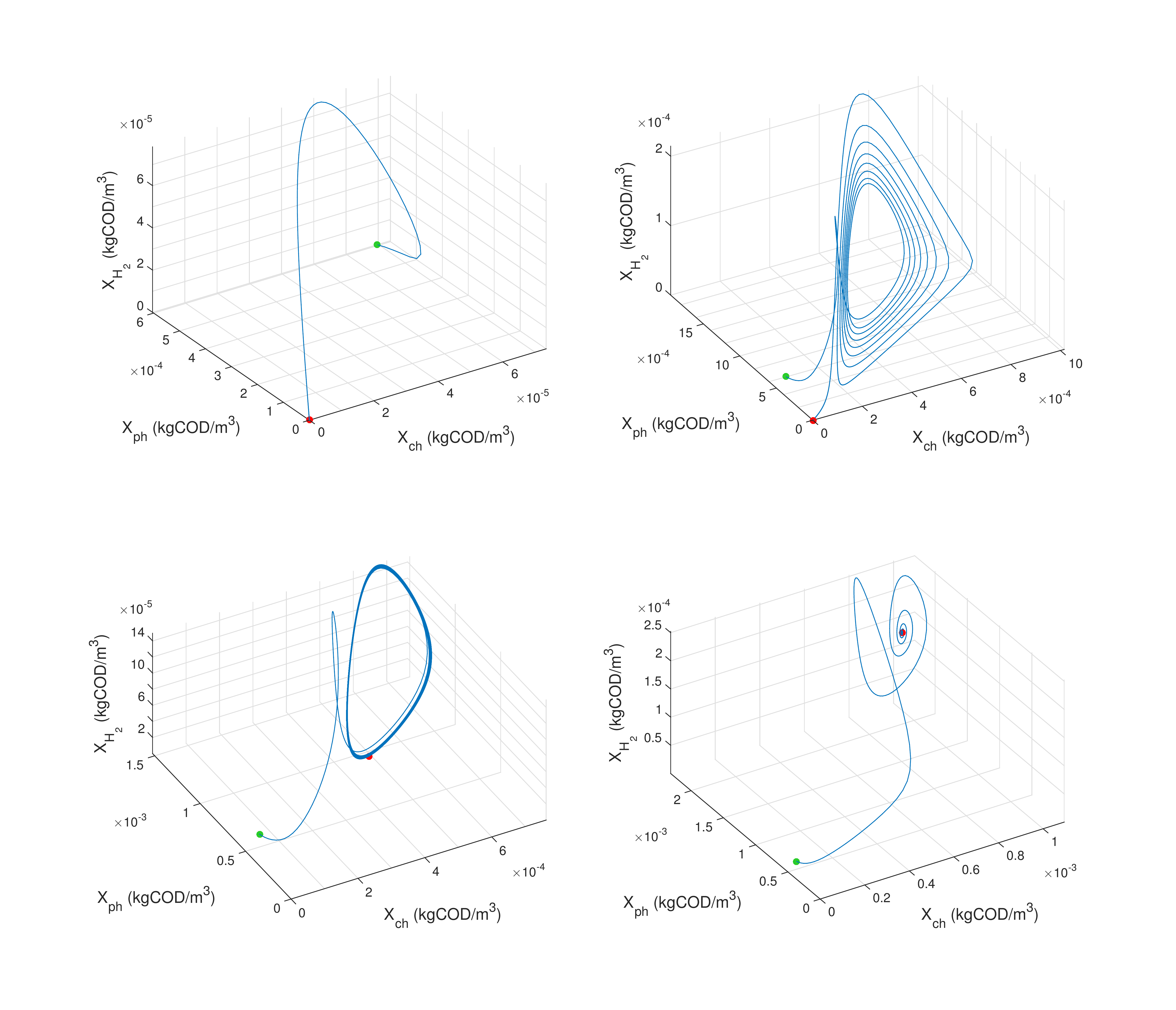}
%%%%%%%%%%%%
\put(-4,4.7){($\alpha$)}
\put(-1.6,4.7){($\beta$)}
\put(-4,2.4){($\gamma$)}
\put(-1.6,2.4){($\delta$)}
%===============
\end{picture}
\end{center}
\caption{Three-dimensional phase plane diagrams of the biomass dynamics for $t = 10000~d$, showing initial (green dot) and final (red dot) conditions for a dilution rate, $D = 0.01~d^{-1}$ and chlorophenol input, $S_{{\rm ch},in}$ ($kgCOD/m^3$) of $\alpha$) 0.01 - the system converges to SS1, $\beta$) 0.097 - the system enters a periodic orbit of increasing amplitude, ultimately converging to SS1, $\gamma$) 0.10052 - the system is close to a stable limit cycle, $\delta$) 0.16 - the system undergoes damped oscillations and converges to SS3.}
\label{oscil}
\end{figure}

We endeavoured to find numerically the operating conditions under which SS3 is unstable, previously unreported by~\cite{wade15}. Given that we have determined analytically in Proposition~\ref{prop2} that when SS3 is close to $SS2^{\flat}$ it becomes unstable, we performed numerical simulations with the parameters defined in Table~\ref{table0} over an operating region similar to that shown in Fig.~2 from~\cite{wade15} whilst also satisfying our conditions. In Fig.~\ref{od_ins_full} we show the case when maintenance is excluded. When magnified, we observe more clearly that region $\mathcal{J}_5$ does exist for the conditions described above, and also note that the region $\mathcal{J}_4$ occurs in a small area between $\mathcal{J}_1$ and $\mathcal{J}_5$, which corresponds to the results shown in Fig.~\ref{fig3}(ii), and is in agreement with Proposition~\ref{propODa}. In Fig.~\ref{od_ins} we confirm that region $\mathcal{J}_5$ does exist for the conditions described above, when maintenance is included, but could not be determined analytically, the curve $\Gamma_3$ is absent in Fig.~\ref{fig3}(i).  Furthermore, we demonstrate that a Hopf bifurcation occurs along the boundary of $F_3(D)$ for values of $D<D_3$ by selecting values of $S_{\rm{ch},in}$ (indicated by $(\alpha)-(\delta)$ in Fig.~\ref{od_ins}) at a fixed dilution rate of $0.01~d^{-1}$, and running dynamic simulations for $10000~d$. The three-dimensional phase plots, with the axes representing biomass concentrations, are shown in Fig.~\ref{oscil}, and show that as $S_{\rm{ch},in}$ approaches $\mathcal{J}_3$ from $\mathcal{J}_5$, emergent periodic orbits are shown to diminish to a stable limit cycle at the boundary (see Appendix~\ref{Hopf} for proof). Subsequently, increasing $S_{\rm{ch},in}$ to $\mathcal{J}_3$ results in the orbit reducing to a fixed point equilibrium at SS3.

\subsection{Operating diagram: case (b)}
\begin{figure}[ht!]
\setlength{\unitlength}{1cm}
\begin{picture}(6,4.2)(0.2,0.2)
\includegraphics[scale=0.138]{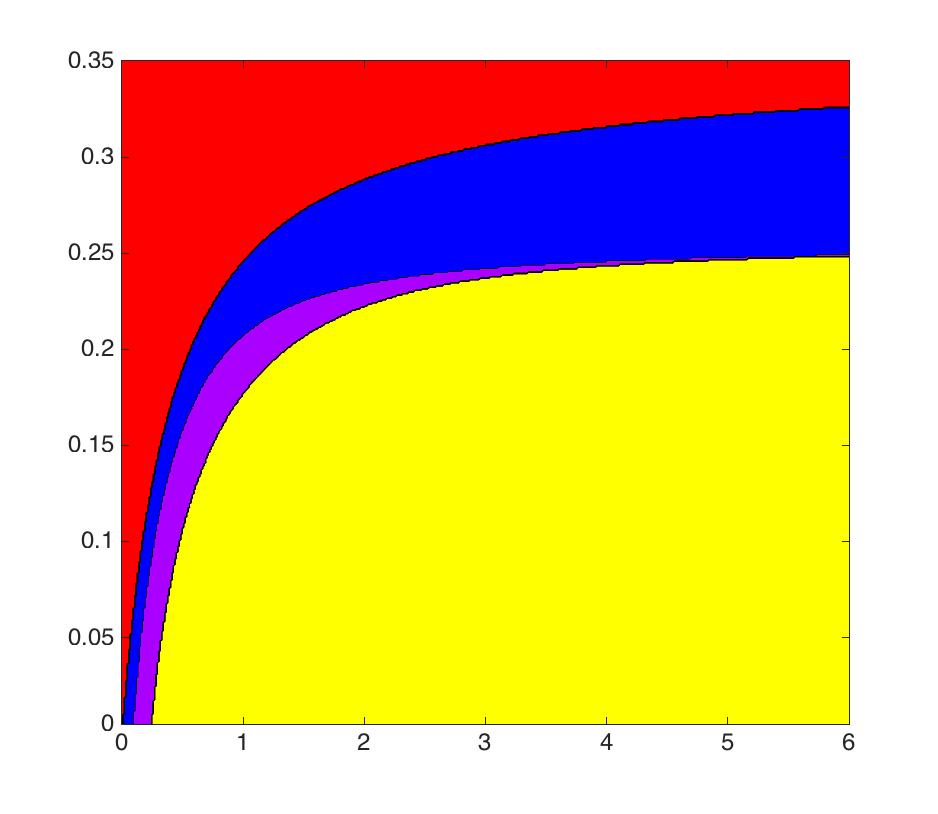}
\put(-3.7,3.1){{\color{white}{$\mathcal{J}_1$}}}
\put(-1.8,3){{\color{white}{$\mathcal{J}_4$}}}
\put(-1.8,1.5){{\color{black}{$\mathcal{J}_3$}}}
\put(-2.8,1.8){{\color{black}{$\mathcal{J}_5$}}}
\put(-2.8,1.9){\vector(-1,1){0.4}}
\put(-4.5,3.4){\color{black}{\scriptsize(i)}}
\put(-4.5,2){{\rotatebox{90}{\scriptsize{$D$}}}}
\put(-2.5,0.1){{\scriptsize{$S_{{\rm ch,in}}$}}}
\end{picture}
\begin{picture}(0,4.2)(1.9,0.2)
\includegraphics[scale=0.138]{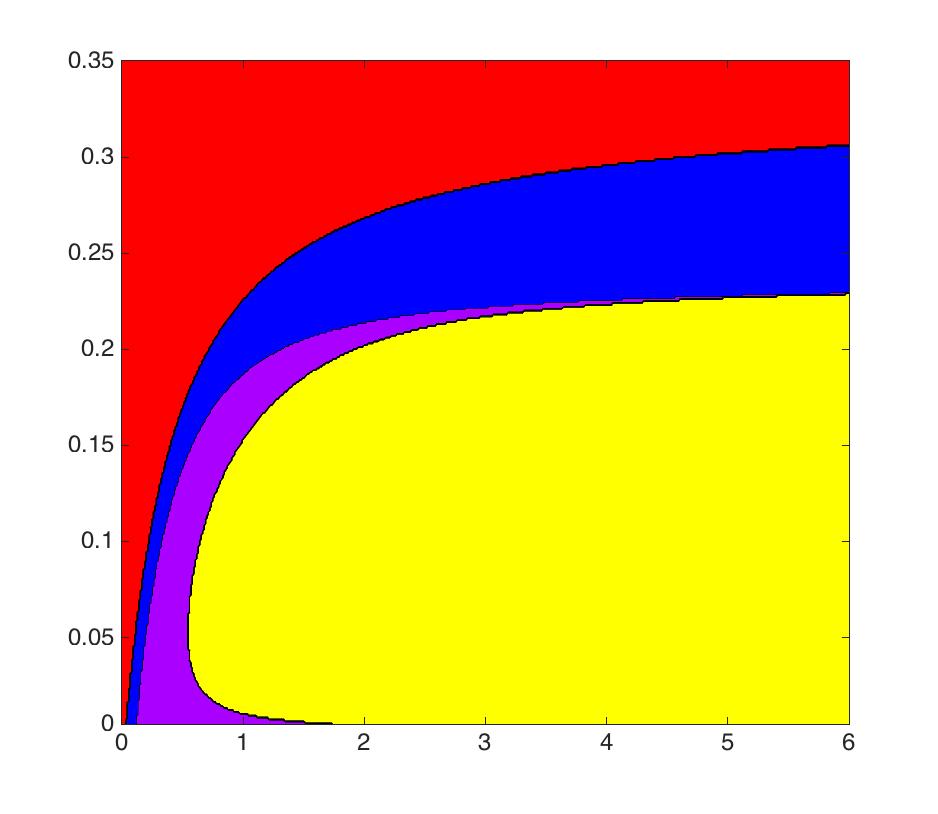}
\put(-3.7,3.1){{\color{white}{$\mathcal{J}_1$}}}
\put(-1.8,2.7){{\color{white}{$\mathcal{J}_4$}}}
\put(-1.8,1.5){{\color{black}{$\mathcal{J}_3$}}}
\put(-2.8,1.6){{\color{black}{$\mathcal{J}_5$}}}
\put(-2.8,1.7){\vector(-1,1){0.4}}
\put(-4.5,3.4){\color{black}{\scriptsize(ii)}}
\put(-4.5,2){{\rotatebox{90}{\scriptsize{$D$}}}}
\put(-2.5,0.1){{\scriptsize{$S_{{\rm ch,in}}$}}}
\end{picture}
\caption{Numerical analysis for the existence and stability of steady-states for case (b). (i) : without maintenance. (ii) : with maintenance.}
\label{caseb_full}
\end{figure}
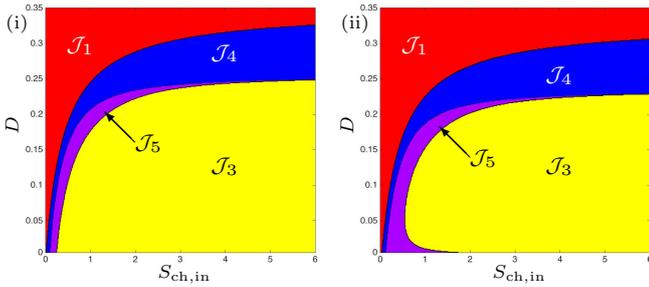

Whilst the numerical parameters chosen for this work are taken from the original study~\cite{wade15}, there somewhat arbitrary nature leaves room to explore the impact of the parameters on the existence and stability of the steady-states. Case (b), discussed in Sections~\ref{sec4.2} and~\ref{sec5.2}, involves a small increase to the half-saturation constant (or inverse of substrate affinity), $K_{S,\rm {H_2,c}}$, of the chlorophenol degrader on hydrogen. Following the same approach as with the preceding case, we confirm in Fig.~\ref{caseb_full}(i) the Proposition~\ref{propODb} in the scenario without maintenance. Furthermore, the extension of this proposition with maintenance included, corresponding to the existence and stability of all three steady-states given in Table~\ref{SSSb}, is show in Fig~\ref{caseb_full}(ii). It shows the region $\mathcal{J}_5$ that cannot be obtained analytically (cf. Fig.~\ref{fig3b}(i)). In both cases, region $\mathcal{J}_2$ has disappeared, as observed analytically. Additionally, the ideal condition $\mathcal{J}_3$, where all organisms are present and stable, diminishes. 

\subsection{Operating diagram: case (c)}
\begin{figure}[ht!]
\setlength{\unitlength}{1cm}
\begin{picture}(6,4.2)(0.1,0.4)
\includegraphics[scale=0.138]{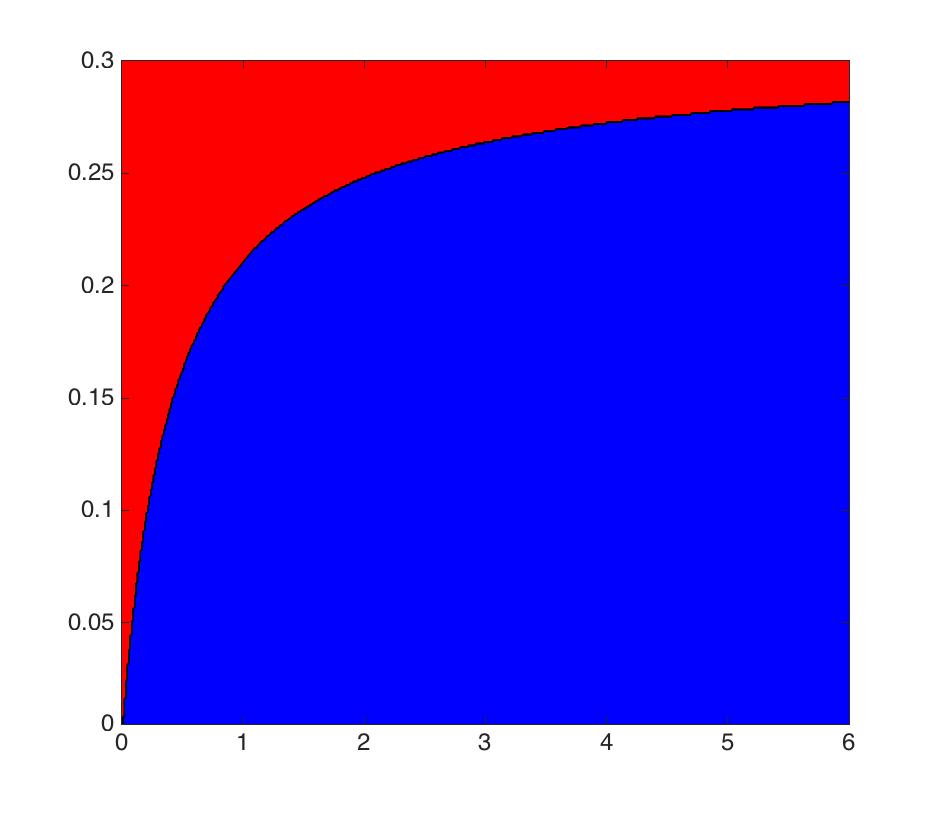}
\put(-3.6,3){{\color{white}{$\mathcal{J}_1$}}}
\put(-1.8,1.6){{\color{white}{$\mathcal{J}_4$}}}
\put(-4.5,3.8){{\color{black}{\scriptsize(i)}}}
\put(-4.5,2){{\rotatebox{90}{\scriptsize{$D$}}}}
\put(-2.5,0.1){{\scriptsize{$S_{{\rm ch,in}}$}}}
\end{picture}
\begin{picture}(0,4.2)(4.2,0.4)
\includegraphics[scale=0.138]{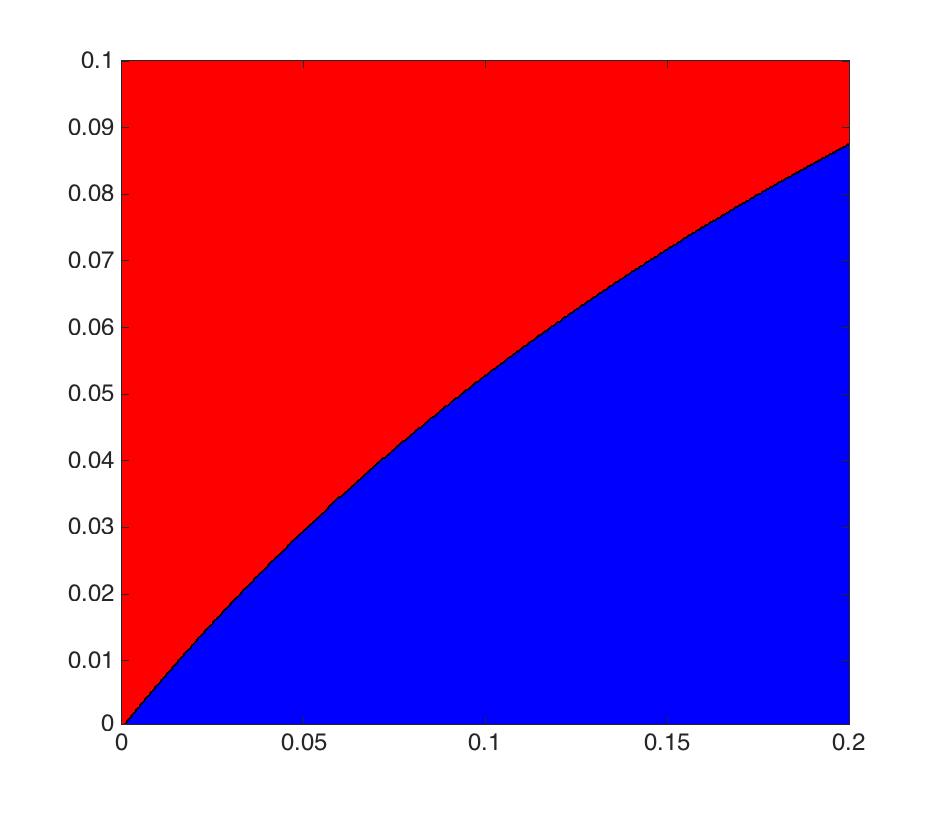}
\put(-3.6,3){{\color{white}{$\mathcal{J}_1$}}}
\put(-1.8,1.6){{\color{white}{$\mathcal{J}_4$}}}
\put(-4.5,2){{\rotatebox{90}{\scriptsize{$D$}}}}
\put(-2.5,0.1){{\scriptsize{$S_{{\rm ch,in}}$}}}
\end{picture}
\begin{picture}(6,4.2)(0.1,0.2)
\includegraphics[scale=0.138]{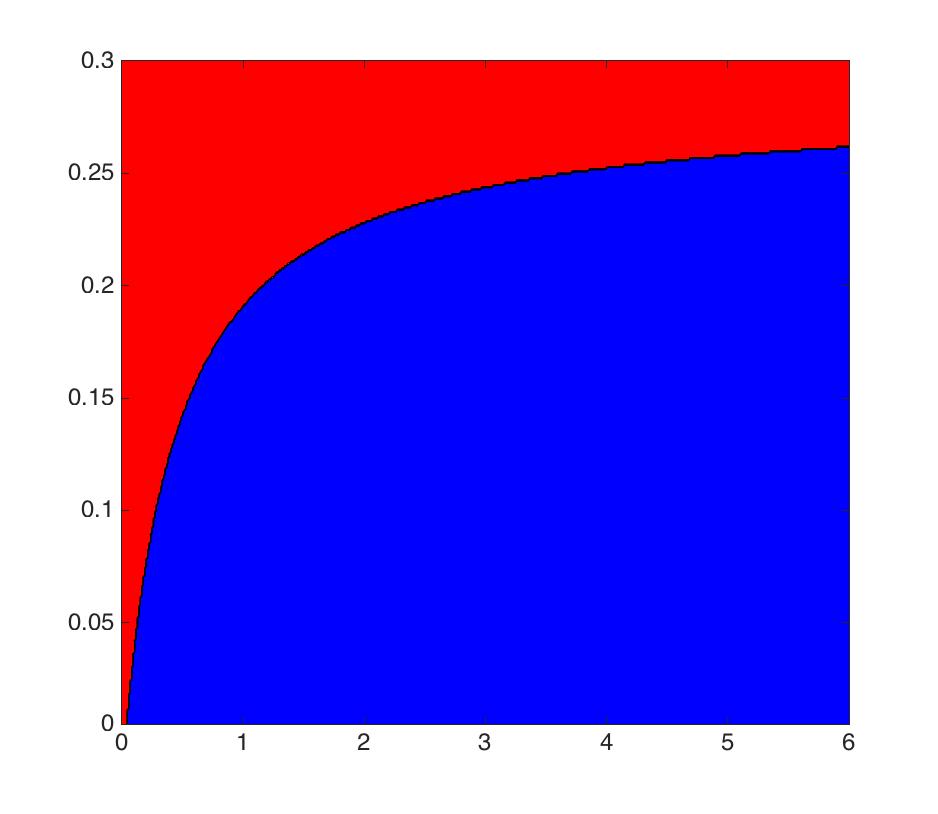}
\put(-3.6,3){{\color{white}{$\mathcal{J}_1$}}}
\put(-1.8,1.6){{\color{white}{$\mathcal{J}_4$}}}
\put(-4.5,3.8){{\color{black}{\scriptsize(ii)}}}
\put(-4.5,2){{\rotatebox{90}{\scriptsize{$D$}}}}
\put(-2.5,0.1){{\scriptsize{$S_{{\rm ch,in}}$}}}
\end{picture}
\begin{picture}(0,4.2)(1.9,0.2)
\includegraphics[scale=0.138]{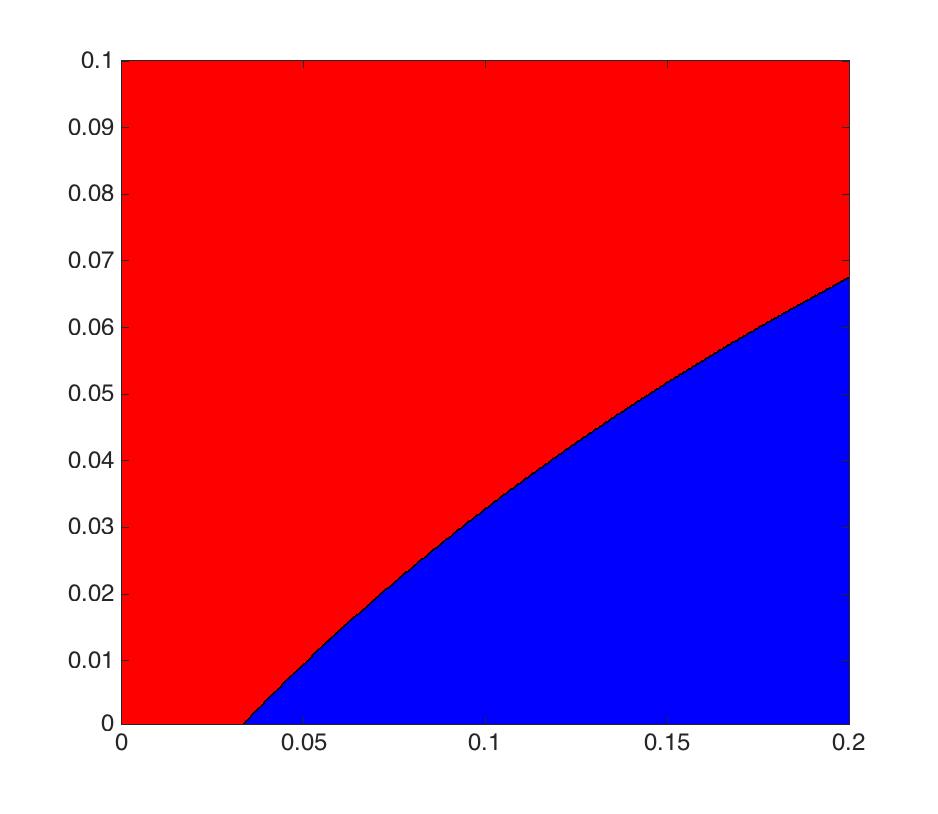}
\put(-3.6,3){{\color{white}{$\mathcal{J}_1$}}}
\put(-1.8,1.6){{\color{white}{$\mathcal{J}_4$}}}
\put(-4.5,2){{\rotatebox{90}{\scriptsize{$D$}}}}
\put(-2.5,0.1){{\scriptsize{$S_{{\rm ch,in}}$}}}
\end{picture}
\caption{Numerical analysis for the existence and stability of steady-states for case (c). (i) : without maintenance. (ii) : with maintenance. On the right, a magnification for $0<D<0.1$.}
\label{par_full}
\end{figure}
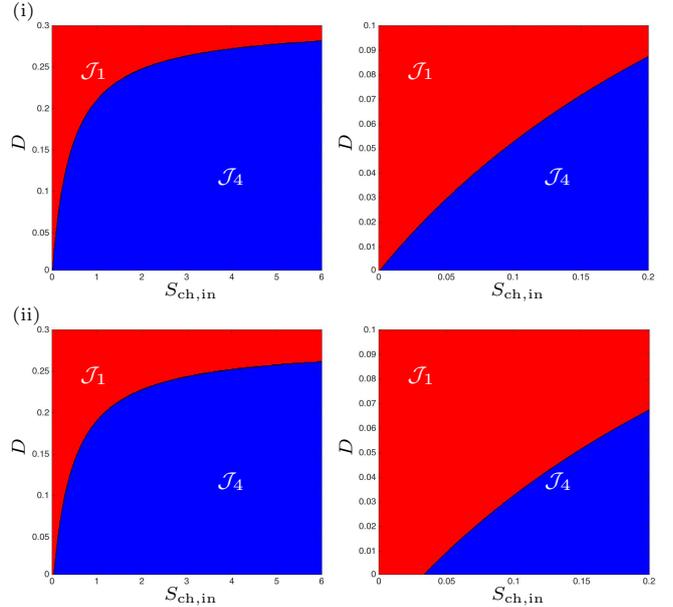

Here, $K_{S,H_2,c}$, was further increased and confirm the Proposition~\ref{propODc}, where the function SS3 never exist and SS2 never stable for the case without maintenance. The extension of this proposition to the case with maintenance, shown in 
Table~\ref{SSSc}, produce similar results as shown in the comparison of Figs.~\ref{par_full}(i) and (ii). 

\subsection{Operating diagram: case (d)}
\begin{figure}[ht!]
\setlength{\unitlength}{1cm}

\begin{picture}(6,4.2)(0.2,0.2)
\includegraphics[scale=0.12]{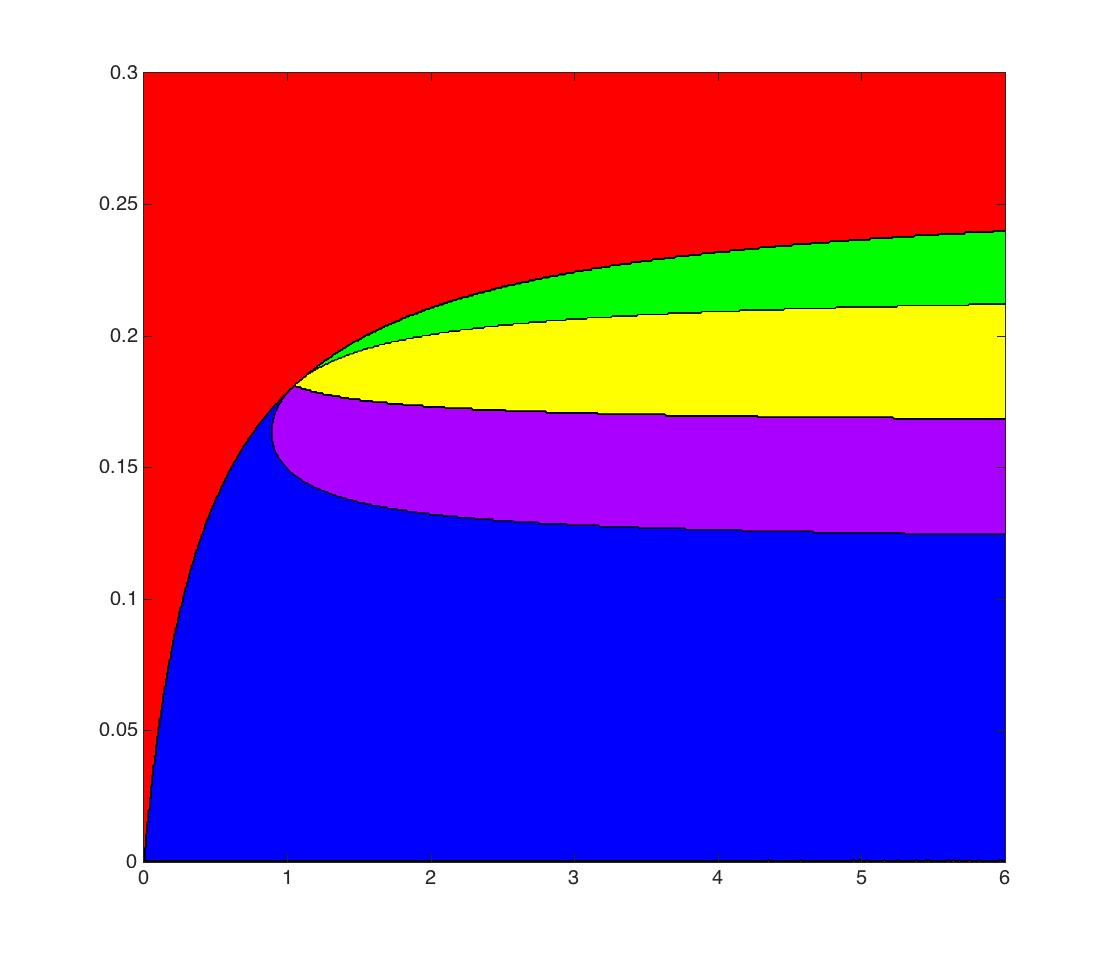}
\put(-3.7,3.1){{\color{white}{$\mathcal{J}_1$}}}
\put(-1.4,2.87){{\color{black}{$\mathcal{J}_2$}}}
\put(-1.4,2.5){{\color{black}{$\mathcal{J}_3$}}}
\put(-1.4,1){{\color{white}{$\mathcal{J}_4$}}}
\put(-1.4,2){{\color{white}{$\mathcal{J}_5$}}}
\put(-4.55,3.6){\color{black}{\scriptsize(i)}}
\put(-4.55,2){{\rotatebox{90}{\scriptsize{$D$}}}}
\put(-2.5,0.1){{\scriptsize{$S_{{\rm ch,in}}$}}}
\end{picture}
\begin{picture}(0,4.2)(2,0.2)
\includegraphics[scale=0.12]{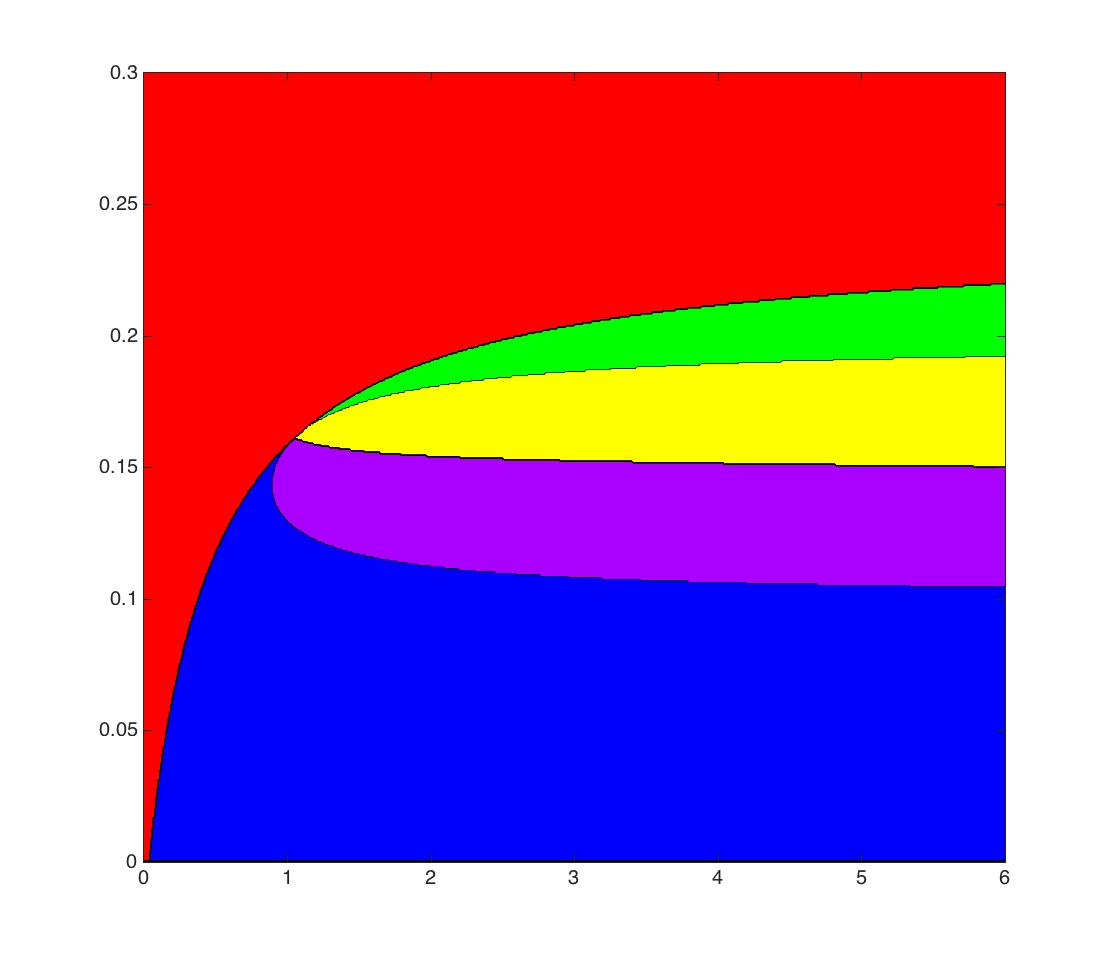}
\put(-3.7,3.1){{\color{white}{$\mathcal{J}_1$}}}
\put(-1.4,2.65){{\color{black}{$\mathcal{J}_2$}}}
\put(-1.4,2.3){{\color{black}{$\mathcal{J}_3$}}}
\put(-1.4,1){{\color{white}{$\mathcal{J}_4$}}}
\put(-1.4,1.85){{\color{white}{$\mathcal{J}_5$}}}
\put(-4.65,3.6){\color{black}{\scriptsize(ii)}}
\put(-4.55,2){{\rotatebox{90}{\scriptsize{$D$}}}}
\put(-2.5,0.1){{\scriptsize{$S_{{\rm ch,in}}$}}}
\end{picture}
\caption{Numerical analysis for the existence and stability of steady-states for case (d). (i) : without maintenance. (ii) : with maintenance.}
\label{cased_full}
\end{figure}
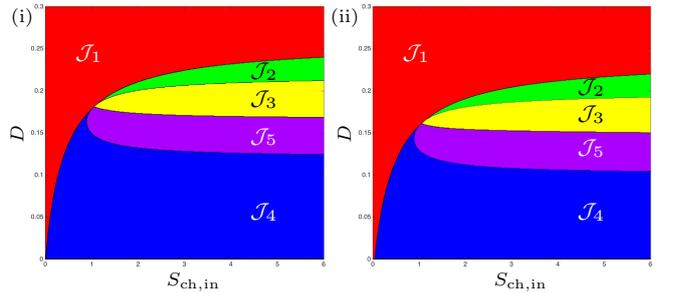

With the final investigated scenario, where $k_{m \rm,H_2} < k_{m, \rm ch}$ and $K_{S,\rm H_2} < K_{S, \rm {H_2,c}}$,  we observe once again the presence of all operating regions, $\mathcal{J}_1 - \mathcal{J}_5$, without and with maintenance, as shown in Fig.~\ref{cased_full}. It can be seen that regions $\mathcal{J}_4$ and $\mathcal{J}_5$ increase at low dilution rates across a much larger range of $S_{\rm ch,in}$ than in the default case (a), and the desired condition (stable SS3) is restricted to a much narrower set of $D$. 

As with the previous cases, the numerical analysis for case (d) confirms the Proposition~\ref{propODd} without maintenance and its extension to the case with maintenance, indicated in Table~\ref{SSS}. 

%=============================
\section{The role of kinetic parameters}\label{kinetic}
%=============================
Finally, we give brief consideration to the characterisation of the four cases discussed in the preceding sections.
%=============================
The main difference between cases (a) or (b) and cases (c) or (d) is that, for small values of $D$, the coexistence steady-state SS3 can exist for cases (a) and (b), but cannot exist for cases (c) or (d). The cases (a) or (b) occur if and only if $s_2^0(0)< M_2(0)$ holds or 
$s_2^0(0)= M_2(0)$ and $\frac{ds_2^0}{dD}(0)< \frac{dM_2}{dD}(0)$ hold, that is to say
\begin{align}
\frac{L_0a_0}{m_0-a_0}<\frac{K_2a_2}{m_2-a_2}\mbox{ or }\label{case(a)or(b)1}\\
\frac{L_0a_0}{m_0-a_0}=\frac{K_2a_2}{m_2-a_2}\mbox{ and }
\frac{L_0m_0}{(m_0-a_0)^2}<\frac{K_2m_2}{(m_2-a_2)^2}
\label{case(a)or(b)2}
\end{align}
The cases (c) or (d) occur if and only if $s_2^0(0)> M_2(0)$ holds or 
$s_2^0(0)= M_2(0)$ and $\frac{ds_2^0}{dD}(0)> \frac{dM_2}{dD}(0)$ hold, that is to say
\begin{align}
\frac{L_0a_0}{m_0-a_0}>\frac{K_2a_2}{m_2-a_2}\mbox{ or }\label{case(c)or(d)2}\\
\frac{L_0a_0}{m_0-a_0}=\frac{K_2a_2}{m_2-a_2}\mbox{ and }
\frac{L_0m_0}{(m_0-a_0)^2}>\frac{K_2m_2}{(m_2-a_2)^2}
\label{case(c)or(d)2}
\end{align}

Notice that it is easy to make the difference between case (c) and case (d): the first occurs when $M_2(D_1)<s_2^0(D_1)$ and the second when
$M_2(D_1)>s_2^0(D_1)$. Since $D_1$ is the positive solution of the algebraic quadratic equation $s_2^0(D)=s_2^1(D)$, it is possible to have an expression for $D_1$ with respect to the biological parameters.
However, this is a complicated expression involving many parameters and the preceding conditions 
$M_2(D_1)<s_2^0(D_1)$ or $M_2(D_1)>s_2^0(D_1)$ have no biological interpretation.
We simply remark here that the function $s_2^0(D)$ has a vertical asymptote for $D=m_0-a_0$ and the function $M_2(D)$ has a vertical asymptote for $D=m_2-a_2$. Therefore, if $m_0-a_0<m_2-a_2$ then case (c) occurs, so that a necessary (but not sufficient) condition for case (d) to occur is $m_0-a_0>m_2-a_2$. If $m_2$ is sufficiently small then case (d) can occur.

The observations from the numerical analysis suggest that the role of the chlorophenol degrader as a secondary hydrogen scavenger is critical in maintaining full chlorophenol mineralisation and system stability, particularly at higher dilution rates, as shown by comparing cases (c) and (d) . More significantly, the results coupled with the parameter relationships shown in Eqs.~\ref{case(a)or(b)1}-\ref{case(c)or(d)2}, highlight the necessary conditions under which the ideal case ($SS3$ stable) is achieved and, in general, this is a coupling of the two key parameters describing the half-saturation constant and maximum specific growth rates between the two hydrogen competitors. 

%============================= 

\section{Conclusions}
In this work we have generalised a simplified mechanistic model describing the anaerobic mineralisation of chlorophenol in a two-step food-web. We are able to show complete analytical solutions describing the existence and stability of the steady-states in the case that maintenance is excluded from the system, whilst with a decay term present, purely analytical determination of stability is not possible. 

We confirm the findings of previous numerical analysis by~\cite{wade15} that with chlorophenol as the sole input substrate, three steady-states are possible. However, the analysis goes further and we determine that under certain operating conditions, two of these steady-states (SS2 and SS3) can become stable, whilst SS1 always exists and is always stable. Furthermore, without maintenance we can explicitly determine the stability of the system, and form analytical expressions of the boundaries between the different stability regions.  

As the boundary of $\mathcal{J}_3$ is not open to analytical determination in the case with maintenance, we determined numerically (substituting the general growth function with the classical Monod-type growth kinetics) the existence and stability of the system over a range of practical operating conditions (dilution rate and chlorophenol input). For comparison and confirmation, we also performed this for the case without maintenance and found the same regions in both cases, with variations only in their shape and extent. For example, whilst the boundary between $\mathcal{J}_1$ and $\mathcal{J}_4$ terminates at the origin without maintenance, with maintenance it is located at $F_1(0)/Y_3Y_4\approx 0.0195$. More interestingly, the addition of a decay term results in an extension of the SS3 unstable steady-state, reducing the potential for successful chlorophenol demineralisation at relatively low dilution rates and substrate input concentrations. Additionally, we show that at the boundary between $\mathcal{J}_3$ and $\mathcal{J}_5$, a Hopf bifurcation occurs and a limit cycle in SS3 emerges.

Finally, we gave an example of how the model could be used to probe the system to answer specific questions regarding model parameterisation. Here we have indicated that a switch in dominance between two organisms competing for hydrogen results in the system becoming unstable and a loss in viability. This is perhaps intuitive to microbiologists, but here it has been proven using mathematical analysis, and could be used to determine critical limits of the theoretical parameter values in shifting between a stable and unstable system. Whilst parameters are not arbitrary in real organisms, the potential for microbial engineering or synthetic biology to manipulate the properties of organisms makes this observation all the more pertinent.

%=============================

\begin{appendix}
\section{Numerical methods}\label{nummeth}
We consider sets of operating parameters ($D$ and $S_{\mathrm{ch,in}}$) for each of the three steady-states, and using Matlab, the complex polynomials for each steady-state can be solved by substitution of parameter values (see Table~\ref{table0}) into the explicit solution. By investigating the signs of the solutions and the eigenvalues, respectively, we determine which steady-states are meaningful and stable. By exploring a localised region of suitable operating parameters, we then generate a phase plot showing where each steady-state is stable, bistable or unstable. 

\section{Proof for Hopf Bifurcation}\label{Hopf}
\begin{figure}[ht!]
\setlength{\unitlength}{1.5cm}
\begin{center}
\begin{picture}(8,6)(0.4,0.5)
\includegraphics[scale=0.255]{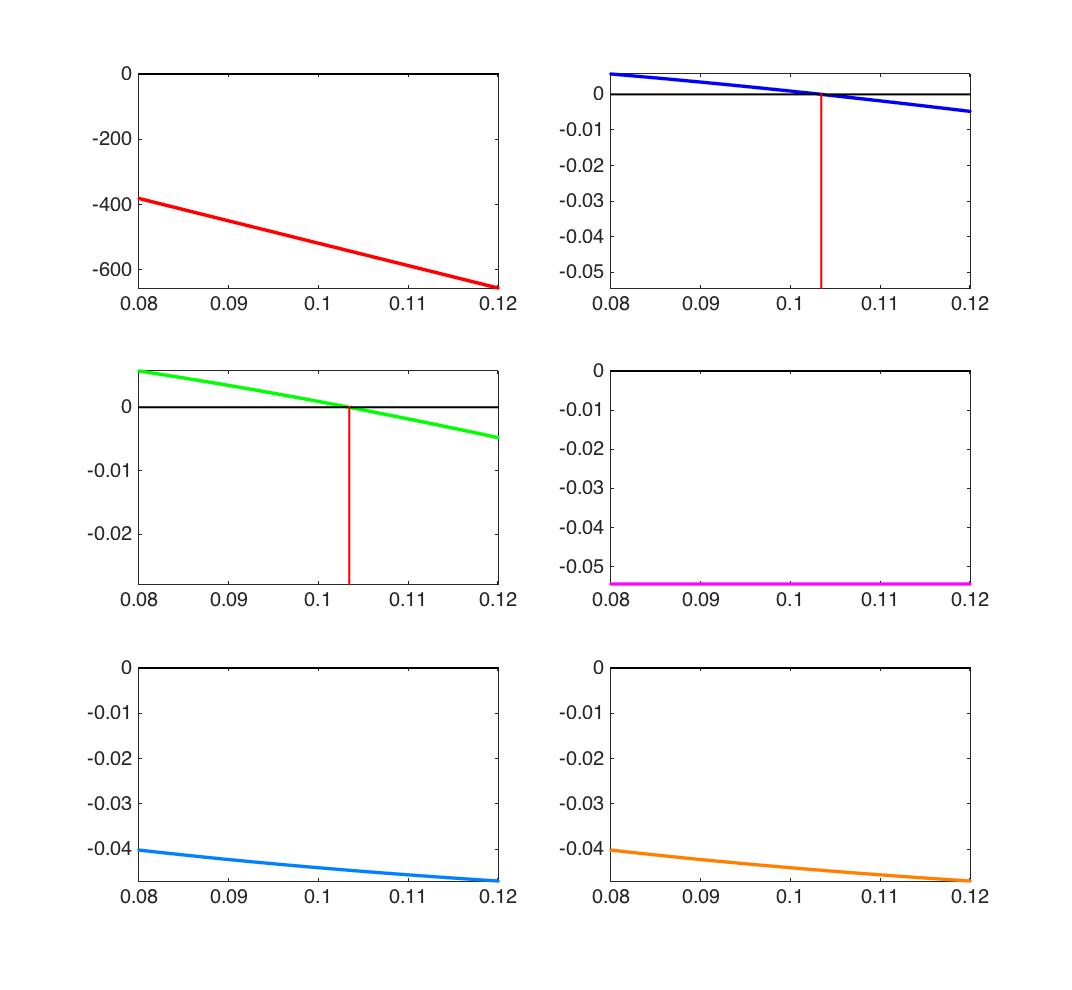}
\put(-6.1,4.4){{\rotatebox{90}{\scriptsize{Eigenvalue 1}}}}
\put(-3.3,4.4){{\rotatebox{90}{\color{red}{\scriptsize{Eigenvalue 2}}}}}
\put(-6.1,2.5){{\rotatebox{90}{\color{red}{\scriptsize{Eigenvalue 3}}}}}
\put(-3.3,2.5){{\rotatebox{90}{\scriptsize{Eigenvalue 4}}}}
\put(-6.1,0.75){{\rotatebox{90}{\scriptsize{Eigenvalue 5}}}}
\put(-3.3,0.75){{\rotatebox{90}{\scriptsize{Eigenvalue 6}}}}
\put(-3.5,0.3){{{$S_{{\rm ch,in}}$}}}
\end{picture}
\end{center}
\caption{Real parts of the eigenvalues determined at $D = 0.01$ and $S_{\rm ch,in} = [0.08,0.12]$, in the case with maintenance. The red vertical lines indicate the location where the eigenvalue crosses zero.}
\label{hopfproof}
\end{figure}
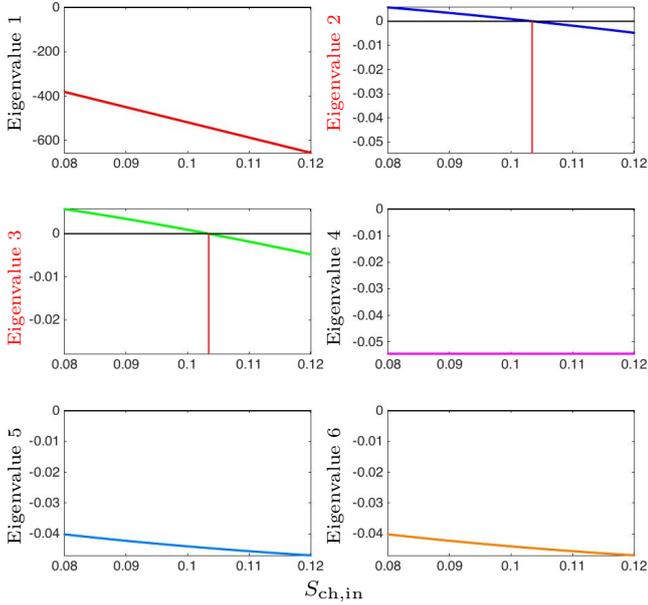

In Section~\ref{numanaly}, we show the operating diagrams with the parameters given in Table~\ref{table0}, and determine numerically that as the parameter $S_{\rm ch,in}$ increases at a fixed dilution rate ($D = 0.01~d^{-1}$), the system bifurcates through several stability domains. We claim that as we cross the boundary between regions $\mathcal{J}_5$ and $\mathcal{J}_3$, we observe a Hopf bifurcation, and, in $\mathcal{J}_5$, close to the boundary with $\mathcal{J}_3$, a limit cycle appears. In order to test this numerically, we checked the real parts of the six eigenvalues at each point along the transect shown in Fig~\ref{od_ins} (10000 points in total), and plotted their values. Fig~\ref{hopfproof} indicates the conditions for a Hopf bifurcation are satisfied as eigenvalues 2 and 3 both change their sign when passing through the coordinate (0, 0.1034) and the real part of all eigenvalues 1, 4 and 6 remain negative. 

\section{General case}\label{general}
As mentioned at the end of Section \ref{sec:model} our study does not require that growth functions are of Monod type (Eq.~\ref{monMonod}). Actually, the results are valid for a more general class of growth functions satisfying the following conditions, which concur with those given by Eq.~\ref{monMonod}:
\begin{description}
\item[H1] For all $s_0>0$ and $s_2> 0$ then $0<\mu_0\left(s_0,s_2\right)<+\infty$ and 
$\mu_0\left(0,s_{2}\right) = 0$, 
$\mu_0\left(s_0,0\right) = 0$. 
\item[H2] For all $s_1>0$ and $s_2\geq 0$ then $0<\mu_1\left(s_1,s_2\right)<+\infty$ and 
$\mu_1\left(0,s_{2}\right) = 0$.
\item[H3] For all $s_2>0$ then $0<\mu_2\left(s_2\right)<+\infty$ and 
$\mu_2(0) = 0$.
\item[H4] For all $s_{0}>0$ and $s_2>0$,
$$
\displaystyle \frac{\partial \mu_0}{\partial s_0}\left(s_0,s_2\right) > 0,
\quad
\displaystyle \frac{\partial \mu_0}{\partial s_2}\left(s_0,s_2\right) > 0.
$$
\item[H5] For  all $s_{1}>0$ and $s_2>0$,
$$
\displaystyle \frac{\partial \mu_1}{\partial s_1}\left(s_1,s_2\right) > 0,
\quad
\displaystyle \frac{\partial \mu_1}{\partial s_2}\left(s_1,s_2\right) < 0.
$$
\item[H6] For all $s_2>0$,
$\displaystyle \frac{d\mu_2}{ds_2}\left(s_2\right)> 0$.
\item[H7] The function $s_2\mapsto \mu_0(+\infty,s_2)$ is monotonically increasing and the function
$s_2\mapsto \mu_1(+\infty,s_2)$ is monotonically decreasing.
\end{description}

We use Eq. \ref{eqM0explicit}, Eq. \ref{eqM1explicit} and Eq. \ref{eqM2explicit} to define $M_0(y,s_2)$, $M_1(y,s_2)$ and $M_2(y)$, respectively.
\begin{lemma}\label{lemmaM0}
Let $s_2\geq 0$ be fixed. There exists a unique function  
%===========================
$$y\in[0,\mu_0(+\infty,s_2))\mapsto M_0(y,s_2)\in[0,+\infty),$$
%===========================
such that for $s_0\geq 0$, $s_2\geq 0$ and $y\in[0,\mu_0(+\infty,s_2))$, we have
%================
\begin{equation}
s_0=M_0(y,s_2)\Longleftrightarrow y=\mu_0(s_0,s_2)
\label{eqM0}
\end{equation}
\end{lemma}

\begin{lemma}\label{lemmaM1}
Let $s_2\geq 0$ be fixed. There exists a unique function
%========================
$$y\in[0,\mu_1(+\infty,s_2))\mapsto M_1(y,s_2)\in[0,+\infty),$$ 
%========================
such that for $s_1\geq 0$, $s_2\geq 0$ and $y\in[0,\in[0,\mu_1(+\infty,s_2))$, we have
%===========================
\begin{equation}
s_1=M_1(y,s_2)\Longleftrightarrow y=\mu_1(s_1,s_2)
\label{eqM1}
\end{equation}
\end{lemma}

\begin{lemma}\label{lemmaM2}
There exists a unique function
%=====================
$$y\in[0,\mu_2(+\infty))\mapsto M_2(y)\in[0,+\infty),$$
%=========================
such that, for
$s_2\geq 0$ and $y\in[0,\mu_2(+\infty))$ we have
%==================
\begin{equation}
s_2=M_2(y)\Longleftrightarrow y=\mu_2(s_2)
\label{eqM2}
\end{equation}
\end{lemma}

We use Eq. \ref{eqs201D} to define the functions $s_2^1(D)$ and $s_2^1(D)$
\begin{lemma}\label{lemmas201D}
For  
$D+a_0<\mu_0(+\infty,+\infty)$ and $D+a_1<\mu_1(+\infty,0)$ there exist unique values $s_2^0$ and $s_2^1$ 
such that
%====================
%========================================================
\begin{equation}
\mu_0(+\infty,s_2^0)=D+a_0,\quad\mu_1(+\infty,s_2^1)=D+a_1
\label{eqs201}
\end{equation}
\end{lemma}

Let $\omega<1$. We use Eq.~\ref{eqpsiexplicit} to define $\psi(s_2,D)$ in the general case: we let $\psi:(s_2^0,s_2^1)\longrightarrow\mathbb{R}$ defined by
%==================
\begin{equation}
\psi(s_2)=M_0(D+a_0,s_2)+\frac{M_1(D+a_1,s_2)+s_2}{1-\omega},
\label{eqpsi}
\end{equation}
%=================
It should be noted that $\psi(s_2)>0$ for $s_2^0<s_2<s_2^1$.
From Eq.~\ref{eqM0}, Eq.~\ref{eqM1} and Eq.~\ref{eqs201} we deduce that
$$M_0(D+a_0,s_2^0)=+\infty,\quad M_1(D+a_1,s_2^1)=+\infty$$
Therefore, we have
 $$\lim_{s_2\to s_2^0}\psi(s_2)=\lim_{s_2\to s_2^1}\psi(s_2)=+\infty$$ 
Hence, the function $\psi(s_2)$, which is positive and tends to $+\infty$ at the extremities of the interval $(s_2^0,s_2^1)$, has a minimum value on this interval. We add the following assumption:
\begin{description}
\item[H8] The function $\psi$ has a unique minimum $\overline{s}_2$ on the interval $\left(s_2^0,s_2^1\right)$ and 
$\frac{d\psi}{ds_2}(s_2)$ is negative on $\left(s_2^0,\overline{s}_2\right)$ and positive on $\left(\overline{s}_2,s_2^1\right)$, respectively. 
\end{description}

The function $\psi$ together with the values $s_2^0$, $s_2^1$ and $\overline{s}_2$ all depend on $D$.
However, to avoid cumbersome notations we will use the more precise notations 
$\psi(s_2,D)$, $s_2^0(D)$, $s_2^1(D)$ and $\overline{s}_2(D)$ only if necessary.

We use Eq.~\ref{eqF1explicit}, Eq.~\ref{eqF2explicit} and Eq.~\ref{eqF3explicit} to define $F_1(D)$, $F_2(D)$ and $F_3(D)$ in the general case:
\begin{align}
F_1(D)&=\inf_{s_2\in(s_2^0,s_2^1)}\psi(s_2)=\psi\left(\overline{s}_2\right)
\label{eqF1}\\
F_2(D)&=\psi\left(M_2(D+a_2)\right)	
\label{eqF2}\\
F_3(D)&=\frac{d\psi}{ds_2}\left(M_2(D+a_2)\right)
\label{eqF3}
\end{align}
%============================
The function $F_1(D)$ is defined for 
	$$D\in I_1=\{D\geq 0:s_2^0(D)<s_2^1(D)\}$$
The function $F_2(D)$ and $F_3(D)$ are defined for 
	$$D\in I_2=\{D\in I_1:s_2^0(D)<M_2(D+a_2)<s_2^1(D)\}$$
For all  for $D\in I_2$, $F_1(D)\leq F_2(D)$. 
The equality $F_1(D)= F_2(D)$ holds if, and only if, $M_2(D+a_2)=\overline{s}_2(D)$ that is, 
	$\frac{d\psi}{ds_2}\left(M_2(D+a_2)\right)=0$, that is if, and only if, $F_3(D)=0$.
 
As it will be shown in Appendix \ref{proofs}, the Lemmas \ref{lemma2} and \ref{lemma3}, stated in Section \ref{sec:SS} in the particular case of the Monod type growth functions (Eq.~\ref{monMonod}), are true in the general case of growth functions satisfying assumptions {\bf H1}--{\bf H8}. 
 
\section{Proofs}\label{proofs}
In this Section we give the proofs of the results. 
In these proofs, we do not assume that the growth function are of Monod type (Eq.~\ref{monMonod}). We only assume that the growth functions satisfy {\bf H1}--{\bf H8}. 
\subsection{Existence of steady-states}
\begin{proof} {[Lemma \ref{lemma1}]}
Assume first that $x_0=0$. Then, as a consequence of Eq.~\ref{m4}, we have $s_0=s_0^{{\rm in}}$ and, as a consequence of Eq.~\ref{m5}, we have 
$$D s_1+\mu_1(s_1,s_2)x_1=0$$ 
which implies $s_1=0$ and $\mu_1(s_1,s_2)x_1=0$.  Therefore, as a consequence of Eq.~\ref{m2} we have $x_1=0$. 
Replacing $x_0=0$ and $x_1=0$ in Eq.~\ref{m6}, we have 
$$D s_2+\mu_2(s_2)x_2=0$$
which implies $s_2=0$ and $\mu_2(s_2)x_2=0$.  
Therefore, as a consequence of Eq.~\ref{m3} we have $x_2=0$. Hence, the steady-state is SS1.\\
Assume now that $x_1=0$. Then, as a consequence of Eq.~\ref{m6}, we have
$$D s_2+\omega\mu_0(s_0,s_2)x_0+\mu_2(s_2)x_2=0$$
which implies $s_2=0$, $\mu_0(s_0,s_2)x_0=0$ and $\mu_2(s_2)x_2=0$. Therefore, as a consequence of Eq.~\ref{m1}, we have $x_0=0$. As shown previously this implies that the steady-state is SS1. 
Evaluated at SS1 the Jacobian matrix of Eqs.~\ref{modela}-\ref{modelf} is 
$$
\left[
\begin{array}{cccccc}
-D-a_0  & 0    & 0    & 0  & 0  & 0\\
0     & -D-a_1 & 0    & 0  & 0  & 0\\
0     & 0    & -D-a_2 & 0  & 0  & 0\\
0     & 0    & 0    & -D & 0  & 0\\
0     & 0    & 0    & 0  & -D & 0\\
0     & 0    & 0    & 0  & 0  & -D
\end{array}
\right]
$$
Thus, SS1 is stable.
\qed
\end{proof}

\begin{proof} {[Lemma \ref{lemma2}]}
Since $x_0>0$ and $x_1>0$, then, as a consequence of Eq.~\ref{m1} and Eq.~\ref{m2}, we have 
$$\mu_0(s_0,s_2)=D+a_0,\quad \mu_1(s_1,s_2)=D+a_1$$
Hence, we have  
\begin{equation}
s_0=M_0(D+a_0,s_2),\quad s_1=M_1(D+a_1,s_2)
\label{eqs0s1}
\end{equation}
$$$$
Using Eq.~\ref{m4} and Eq.~\ref{m5}, we have Eq.~\ref{eqSS5x0x1}.
Using Eq.~\ref{m6} we have
\begin{equation}
-s_2+(s_0^{\rm in}-s_0-s_1)-\omega(s_0^{{\rm in}}-s_0)=0
\label{eqs2}
\end{equation}
If $\omega \geq 1$ this equation has no solution. If $\omega < 1$ this equation is equivalent to 
$$s_0^{{\rm in}}=s_0+\frac{s_1+s_2}{1-\omega}.$$
Using Eq.~\ref{eqs0s1} we see that $s_2$ must be a solution of Eq.~\ref{eqSS5}. 
Since $s_1>0$ and $s_2>0$ then, from Eq.~\ref{eqs2} we have necessarily
$$s_1+s_2=(1-\omega)(s_0^{{\rm in}}-s_0)>0$$
so that $s_0^{{\rm in}}-s_0>0$.  
From Eq.~\ref{eqSS5x0x1} we deduce that $x_0>0$.
Since $s_0^{{\rm in}}-s_0>0$ and $s_2>0$ then, from Eq.~\ref{eqs2} we have necessarily
$$\omega(s_0^{{\rm in}}-s_0)+s_2=s_0^{{\rm in}}-s_0-s_1>0$$
so that $s_0^{{\rm in}}-s_0-s_1>0$
From Eq.~\ref{eqSS5x0x1} we deduce that $x_1>0$. 
\qed
\end{proof}

\begin{proof}{[Lemma \ref{lemma3}]}
Since $x_0>0$, $x_1>0$ and $x_2>0$, then, as a consequence of Eq.~\ref{m1}, Eq.~\ref{m2}) and Eq.~\ref{m3}, 
we have
$$\mu_0(s_0,s_2)=D+a,\quad \mu_1(s_1,s_2)=D+b,\quad \mu_2(s_2)=D+c$$
Hence, $s_0$, $s_1$ and $s_2$ are given by Eq.~\ref{eqSS8}. 
Using Eq.~\ref{m4}, Eq.~\ref{m5} and Eq.~\ref{m6} we have Eq.~\ref{eqSS8x2}.
For $x_2$ to be positive it is necessary that $s_0$, $s_1$ and $s_2$ satisfy the condition 
\begin{equation}
(1-\omega)(s_0^{{\rm in}}-s_0)>s_1+s_2,
\label{eqSS8existence}
\end{equation}
If $\omega \geq 1$ this equation has no solution. If $\omega < 1$ this equation is equivalent to the condition 
$$s_0^{{\rm in}}>s_0+\frac{s_1+s_2}{1-\omega}.$$
Using Eq.~\ref{eqSS8}, this condition is the same as 
$$s_0^{{\rm in}}>\psi\left(M_2(D+a_2)\right)=F_2(D)$$
Therefore, from Eq.~\ref{eqSS8existence} we have $s_0^{{\rm in}}-s_0>0$ and $s_0^{{\rm in}}-s_0-s_1>0$, so that $x_0>0$ and $x_1>0$. \qed 
\end{proof}
\subsection{Stability of steady-states}
We use the change of variables
\begin{equation}
z_0=s_0+x_0,\quad z_1=s_1+x_1-x_0,\quad z_2=s_2+x_2+\omega x_0-x_1
\label{eq:chvar}
\end{equation}
Therefore, Eqs.~\ref{modela}-\ref{modelf}, with $a_0=a_1=a_2=0$, become
\begin{align}
\frac{{\rm d}x_{0}}{{\rm d}t} &=  -D x_{0} + \mu_{0}\left(z_{0}-x_0,z_{2}-\omega x_0+x_1-x_2\right) x_{0}\label{modelza}\\%[4mm]
\frac{{\rm d}x_{1}}{{\rm d}t} &=  -D x_{1} + \mu_{1}\left(z_{1}+x_0-x_1,z_{2}-\omega x_0+x_1-x_2\right) x_{1}\label{modelzb}\\%[4mm]
\frac{{\rm d}x_{2}}{{\rm d}t} &=  -D x_{2} + \mu_{2}\left(z_{2}-\omega x_0+x_1-x_2\right) x_{2}\label{modelzc}\\%[4mm]
\frac{{\rm d}z_{0}}{{\rm d}t} &= D \left(s_{0}^{{\rm in}} - z_{0}\right)\label{modelzd}\\%[4mm]
\frac{{\rm d}z_{1}}{{\rm d}t} &= -Dz_{1}\label{modelze}\\%[4mm]
\frac{{\rm d}z_{2}}{{\rm d}t} &=  -Dz_{2} \label{modelzf}
\end{align}
In the  variables $(x_0,x_1,x_2,z_0,z_1,z_2)$ where $z_1$, $z_2$ and $z_3$ are defined by Eq.~\ref{eq:chvar}, 
the steady-states SS1, SS2 and SS3 are given by
\begin{enumerate}
	\item ${\rm SS1}= (0,0,0,s_0^{{\rm in}},0,0)$
	\item ${\rm SS2}= (x_0,x_1,0,s_0^{{\rm in}},0,0)$, where $x_0$ and $x_1$ are defined by Eq. \ref{SS2WM}.
	\item ${\rm SS3}= (x_0,x_1,x_2,s_0^{{\rm in}},0,0)$, where $x_0$, $x_1$ and $x_2$ are defined by Eq. \ref{SS3WM}.
	\end{enumerate}
	
	Let $(x_0,x_1,x_2,s_0^{{\rm in}},0,0)$ be a steady-state. The Jacobian matrix of Eqs.~\ref{modelza}-\ref{modelzf} has the block triangular form  
	$$\mathbf{J}=\left[
	\begin{array}{cc}
	\mathbf{J_{1}}&\mathbf{J_{2}}\\
	0&\mathbf{J_3}	

	\end{array}
	\right]$$
	where 
	$$\mathbf{J_1}=
	\left[
\scriptsize{\begin{array}{ccc}
\mu_0-D-(E+\omega F)x_0    & Fx_0         & -Fx_0    \\
(G+\omega H)x_1          & \mu_1-D-(G+H)x_1 & Hx_1 \\
-\omega Ix_2          & Ix_2        & \mu_2-D-Ix_2
\end{array}}
\right]
	$$
	$$
\mathbf{J_2}=
\left[
\begin{array}{ccc}
 Ex_0         & 0       &  Fx_0\\
 0            & Gx_1    & -Hx_1\\
 0            & 0       &  Ix_2
\end{array}
\right],
\quad
\mathbf{J_3}=
\left[
\begin{array}{ccc}
 -D         & 0       &  0\\
 0            & -D    & 0\\
 0            & 0       &  -D
\end{array}
\right]
$$
and 
$$\small{
E=\frac{\partial \mu_0}{\partial s_0},\quad
F=\frac{\partial \mu_0}{\partial s_2},\quad
G=\frac{\partial \mu_1}{\partial s_1},\quad
H=-\frac{\partial \mu_1}{\partial s_2},\quad
I=\frac{ d\mu_2}{d s_2}}
$$
are evaluated at the steady-state. 

Since $\mathbf{J}$ is a block triangular matrix, its eigenvalues are $-D$ (with multiplicity 3) together with the eigenvalues of the
$3\times 3$ upper-left matrix $\mathbf{J_1}$. 
Note that we have used the opposite sign for the partial derivative 
$H=-\partial \mu_1/\partial s_2$, 
so that all constants involved in the computation become positive, which will simplify the analysis of the characteristic polynomial of $\mathbf{J_1}$.
 
\begin{proof}{[Proposition \ref{prop1}]}
Evaluated at SS2, the matrix $\mathbf{J_1}$ is 
$$
\mathbf{J_1}=\left[
\begin{array}{ccc}
 -(E+\omega F)x_0 & Fx_0  & -Fx_0\\
 (G+\omega H)x_1  & -(G+H)x_1 & Hx_1\\
 0  & 0  & \mu_2-D
\end{array}
\right]
$$
Since $\mathbf{J_1}$ is a block triangular matrix, its eigenvalues are simply $\mu_2-D$, together with the eigenvalues of the $2\times 2$ upper-left matrix. Note that the trace of this $2\times 2$ matrix is negative. Hence, its eigenvalues are of negative real part if, and only if, its determinant is positive, that is if, and only if,
\begin{align} \label{ss5_det}
E(G+H)-(1-\omega)FG > 0
 \end{align}
Using
$$\frac{\partial M_0}{\partial s_2}=-\frac{\partial \mu_0}{\partial s_2}\left[\frac{\partial \mu_0}{\partial s_0}\right]^{-1}=-F/E$$
$$\frac{\partial M_1}{\partial s_2}=-\frac{\partial \mu_1}{\partial s_2}\left[\frac{\partial \mu_1}{\partial s_0}\right]^{-1}=H/G$$
we deduce from
$$\psi(s_2)=M_0(D,s_2)+\frac{M_1(D,s_2)+s_2}{1-\omega}$$
that 
$$
\frac{d\psi}{ds_2}=\frac{\partial M_0}{\partial s_2}+\frac{\frac{\partial M_0}{\partial s_2}+1}{1-\omega}
=-\frac{F}{E}+\frac{\frac{H}{G}+1}{1-\omega}
$$
Hence,
\begin{equation}
\frac{d\psi}{ds_2}=\frac{E(G+H)-(1-\omega)FG}{(1-\omega)EG}
\label{Delta}
\end{equation}
Therefore, the condition of stability, (Eq.~\ref{ss5_det}), is equivalent to $\frac{d\psi}{ds_2}>0$. Hence, we have proved that SS2 is stable if, and only if, $\mu_2(s_2)<D$ and $\frac{d\psi}{ds_2}>0$. \qed
\end{proof}

%=========================
\begin{proof}{[Proposition \ref{prop2}]}
Evaluated at SS3, the matrix $\mathbf{J_1}$ is
$$
\mathbf{J_1}=\left[
\begin{array}{ccc}
 -(E+\omega F)x_0 & Fx_0  & -Fx_0\\
 (G+\omega H)x_1  & -(G+H)x_1 & Hx_1\\
 -\omega Ix_2  & Ix_2 & -Ix_2
\end{array}
\right]
$$
The characteristic polynomial is given by
\begin{align}
\lambda^3 + f_2\lambda^2 + f_1\lambda + f_0 = 0
\end{align}
where
\begin{align}
f_2 &= Ix_2+(G+H)x_1+(E+\omega F)x_0\\
f_1 &= \Delta x_0x_1+EIx_0x_2+GIx_1x_2\\
f_0 &= EGIx_0x_1x_2
\end{align}
and
$\Delta=E(G+H)-(1-\omega)FG$.

To satisfy the Routh-Hurwitz criteria, we require $f_i >0$, for $i=0,1,2$ and $f_1f_2-f_0>0$. Notice that
\begin{align}
f_1f_2-f_0 & =
(EIx_0x_2+\Delta x_0x_1)f_2
\nonumber\\
&+(Ix_2+(G+H)x_1+\omega Fx_0)GIx_1x_2
\label{eqF4bis}
\end{align}
We always have $f_0 > 0$ and $f_2 > 0$.

From Eq.~\ref{Delta} we deduce that $\Delta=(1-\omega)EG\frac{d\psi}{ds_2}$ . 
Therefore, if $F_3(D)\geq 0$, that is to say $\frac{d\psi}{ds_2}\geq 0$, then $\Delta>0$. 
Hence, $f_1>0$ and $f_1f_2-f_0>0$, so that SS3 is stable

On the other hand, if $\frac{d\psi}{ds_2}<0$ and $x_2$ is very small, which occurs when SS3 is very close to ${\rm SS2}^{\flat}$, then $f_2$ has the sign of $\Delta$ since the term with $x_2$ is negligible compared to the term $\Delta x_0x_1$:
$$f_2=\Delta x_0x_1+(EIx_0+xGIx_1)x_2<0$$
Thus, SS3 is unstable. \qed
\end{proof}

%=========================
\begin{proof}{[Proposition \ref{prop3}]}
Since we always have $f_0 > 0$ and $f_2 > 0$, from the previous proof it follows that SS3 is stable if, and only if, 
$f_1f_2-f_0>0$. Indeed, this condition implies that we have also $f_1>f_0/f_2>0$. 
Using Eq.~\ref{eqF4} and Eq.~\ref{eqF4bis}, we see that
$$f_1f_2-f_0=F_4\left(D,s_0^{\rm in}\right)$$
Therefore, the condition $f_1f_2-f_0>0$ is equivalent to $F_4\left(D,s_0^{\rm in}\right)>0$. \qed
\end{proof}
%===========================================
\subsection{Operating diagrams}
\begin{proof}{[Proposition \ref{propODa}]} 
We know that SS1 always exist and is stable. We know that ${\rm SS2}^\flat$ is unstable if it exists. Using Table \ref{SSstability} and Remark \ref{remSchin}, we obtain the following results
\begin{itemize}
\item $\mathcal{J}_1$ is defined by $D\geq D_1$ or $0< D<D_1$ and $S_{\rm ch,in}<F_1(D)/Y_3Y_4$. Therefore, SS1 is the only existing steady state in this region.
\item $\mathcal{J}_2$ if defined by $D_3< D< D_1$ and  $F_1(D)/Y_3Y_4<S_{\rm ch,in}<F_2(D)/Y_3Y_4$. Therefore, 
both steady state SS2 exist and ${\rm SS2}^\sharp$ is stable since $F_3(D)>0$.  
\item $\mathcal{J}_3$ if defined by $0< D< D_2$ and  
		$F_2(D)/Y_3Y_4<S_{\rm ch,in}$ and $F_4\left(D,S_{\rm ch,in}/Y_3Y_4\right)>0$ when $0<D<D_3$. Therefore, SS3 exists and is stable, both steady state SS2 exist and ${\rm SS2}^\sharp$ is unstable since $F_3(D)<0$..
\item $\mathcal{J}_4$ if defined by $0< D< D_3$ and  
		$F_1(D)/Y_3Y_4<S_{\rm ch,in}<F_2(D)/Y_3Y_4$. Therefore, 
both steady state SS2 exist and ${\rm SS2}^\sharp$ is unstable since $F_3(D)<0$.         
\item $\mathcal{J}_5$ if defined by $0< D< D_3$,   
		$F_2(D)/Y_3Y_4<S_{\rm ch,in}$ and $F_4\left(D,S_{\rm ch,in}/Y_3Y_4\right)<0$. Therefore, SS3 exists and is unstable and both steady state SS2 exist and ${\rm SS2}^\sharp$ is unstable since $F_3(D)<0$..
\end{itemize} \qed
\end{proof}

\begin{proof}{[Propositions \ref{propODb}, \ref{propODc} and \ref{propODd}]} The result follows from Table \ref{SSstability} and Remark \ref{remSchin}. The details are as in the proof of Proposition \ref{propODb}. \qed
\end{proof}

\subsection{General case}
\begin{proof}{[Lemma \ref{lemmaM0}]}
Let $s_2\geq 0$ be fixed. 
By {\bf H4}, the function 
%==================
$$s_0\in[0,+\infty)\mapsto \mu_0(s_0,s_2)\in[0, \mu_0(+\infty,s_2))$$ 
%==================
is {{monotonically}} increasing. Hence, it has an inverse function denoted by 
%===========================
$$y\in[0,\mu_0(+\infty,s_2))\mapsto M_0(y,s_2)\in[0,+\infty),$$
%===========================
such that for all $s_0\geq 0$, $s_2\geq 0$ and $y\in[0,\mu_0(+\infty,s_2))$ 
(\ref{eqM0}) holds. \qed
\end{proof}

\begin{proof}{[Lemma \ref{lemmaM1}]}
Let $s_2\geq 0$ be fixed. 
By {\bf H5}, the function 
%===========================
$$s_1\in[0,+\infty)\mapsto \mu_1(s_1,s_2)\in[0, \mu_1(+\infty,s_2))$$ 
%===========================
is {{monotonically}} increasing. Hence, it has an inverse function denoted by
%========================
$$y\in[0,\mu_1(+\infty,s_2))\mapsto M_1(y,s_2)\in[0,+\infty),$$ 
%========================
such that for all $s_1\geq 0$, $s_2\geq 0$ and $y\in[0,\in[0,\mu_1(+\infty,s_2))$
%===========================
(\ref{eqM1}) holds. \qed
\end{proof}

\begin{proof}{[Lemma \ref{lemmaM2}]}
By {\bf H6}, the function 
%=====================
$s_2\in[0,+\infty)\mapsto \mu_2(s_2)\in[0,\mu_2(+\infty))$
%======================
is {{monotonically}} increasing. Hence, it has an inverse function denoted by
%=====================
$$y\in[0,\mu_2(+\infty))\mapsto M_2(y)\in[0,+\infty),$$
%=========================
such that, for all 
$s_2\geq 0$ and $y\in[0,\mu_2(+\infty))$
(\ref{eqM2}) holds. \qed
\end{proof}

\begin{proof}{[Lemma \ref{lemmas201D}]}
By {\bf H7}, for 
$D+a_0<\mu_0(+\infty,+\infty)$ and $D+a_1<\mu_1(+\infty,0)$ there exist unique values $s_2^0$ and $s_2^1$ 
such that Eq.~\ref{eqs201} holds, see Fig. \ref{fig1}(a). \qed
\end{proof}
\end{appendix}
\bibliographystyle{elsarticle-num} 
 \bibliography{Sari_Wade2015_08_20}

\begin{thebibliography}{10}
\expandafter\ifx\csname url\endcsname\relax
  \def\url#1{\texttt{#1}}\fi
\expandafter\ifx\csname urlprefix\endcsname\relax\def\urlprefix{URL }\fi
\expandafter\ifx\csname href\endcsname\relax
  \def\href#1#2{#2} \def\path#1{#1}\fi

\bibitem{henze87}
M.~Henze, C.~P. L.~J. Grady, W.~Gujer, G.~v.~R. Marais, T.~Matsuo, Activated
  {S}ludge {M}odel {No. 1}, Tech. Rep.~1, {AWPRC} Scientific and Technical
  Reports, London, UK (1987).

\bibitem{henze99}
M.~Henze, W.~Gujer, T.~Mino, T.~Matsuo, M.~C. Wentzel, G.~v.~R. Marais, M.~C.
  van Loosdrecht, {A}ctivated {S}ludge {M}odel {No. 2D}, {ASM2D}, Wat. Sci.
  Technol. 39~(1) (1999) 165--182.
\newblock \href {http://dx.doi.org/10.1016/S0273-1223(98)00829-4}
  {\path{doi:10.1016/S0273-1223(98)00829-4}}.

\bibitem{batstone02}
D.~J. Batstone, J.~Keller, I.~Angelidaki, S.~V. Kalyuzhnyi, S.~G. Pavlostathis,
  A.~Rozzi, W.~T.~M. Sanders, H.~Siegrist, V.~A. Vavilin, Anaerobic {D}igestion
  {M}odel {N}o. 1, Tech. Rep. Report No. 13, IWA Publishing, London, UK (2002).

\bibitem{jeppsson13}
U.~Jeppsson, J.~Alex, D.~J. Batstone, L.~Benedetti, J.~Comas, J.~B. Copp,
  L.~Corominas, X.~Flores-Alsina, K.~V. Gernaey, I.~Nopens, M.-N. Pons,
  I.~Rodr{\'\i}guez-Roda, C.~Rosen, J.~P. Steyer, P.~A. Vanrolleghem, E.~I.~P.
  Volcke, D.~Vrecko, Benchmark simulation models, \textit{quo vadis}?, Wat.
  Sci. Technol. 68~(1) (2013) 1--15.
\newblock \href {http://dx.doi.org/10.2166/wst.2013.246}
  {\path{doi:10.2166/wst.2013.246}}.

\bibitem{garcia13}
C.~Garc\'{i}a-Di\'{e}guez, O.~Bernard, E.~Roca, Reducing the {A}naerobic
  {D}igestion {M}odel {No. 1} for its application to an industrial wastewater
  treatment plant treating winery effluent wastewater, Biores. Technol. 132
  (2013) 244--253.
\newblock \href {http://dx.doi.org/10.1016/j.biortech.2012.12.166}
  {\path{doi:10.1016/j.biortech.2012.12.166}}.

\bibitem{bornhoft13}
A.~Bornh\"{o}ft, R.~Hanke-Rauschenbach, K.~Sundmacher, Steady-state analysis of
  the {A}naerobic {D}igestion {M}odel {N}o. 1 ({ADM1}), Nonlinear Dynam.
  73~(1-2) (2013) 535--549.
\newblock \href {http://dx.doi.org/10.1007/s110710130807x}
  {\path{doi:10.1007/s110710130807x}}.

\bibitem{weedermann15}
M.~Weedermann, G.~S.~K. Wolkowicz, J.~Sasara, Optimal biogas production in a
  model for anaerobic digestion, Nonlinear Dynam. 81 (2015) 1097--1112.
\newblock \href {http://dx.doi.org/10.1007/s11071-015-2051-z}
  {\path{doi:10.1007/s11071-015-2051-z}}.

\bibitem{xu11}
A.~Xu, J.~Dolfing, T.~P. Curtis, G.~Montague, E.~Martin, Maintenance affects
  the stability of a two-tiered microbial `food chain'?, J. Theor. Biol.
  276~(1) (2011) 35--41.
\newblock \href {http://dx.doi.org/10.1016/j.jtbi.2011.01.026}
  {\path{doi:10.1016/j.jtbi.2011.01.026}}.

\bibitem{sari14b}
T.~Sari, J.~Harmand, Maintenance does not affect the stability of a two-tiered
  microbial 'food chain', {HAL ID}: hal-01026149, version 1 [math.DS] (July
  2014).

\bibitem{wade15}
M.~J. Wade, R.~W. Pattinson, N.~G. Parker, J.~Dolfing, Emergent behaviour in a
  chlorophenol-mineralising three-tiered microbial `food web', arXiv:1503.01580
  [q-bio.PE] (March 2015).

\end{thebibliography}

\end{document}